\documentclass{article}%
\usepackage{amssymb}
\usepackage{amsfonts}
\usepackage{amsmath}
\usepackage{graphicx}%
\providecommand{\U}[1]{\protect\rule{.1in}{.1in}}

\newtheorem{theorem}{Theorem}[section]

\newtheorem{corollary}[theorem]{Corollary}

\newtheorem{definition}[theorem]{Definition}
\newtheorem{example}[theorem]{Example}

\newtheorem{lemma}[theorem]{Lemma}

\newtheorem{proposition}[theorem]{Proposition}
\newtheorem{remark}[theorem]{Remark}

\begin{document}

\title{A Hamilton-Jacobi Theory for general dynamical systems and integrability by
quadratures in symplectic and Poisson manifolds}
\author{Sergio \ Grillo\\{\small Instituto Balseiro, Universidad Nacional de Cuyo and CONICET}\\[-5pt] {\small Av. Bustillo 9500, San Carlos de Bariloche}\\[-5pt] {\small R8402AGP, Rep\'{u}blica Argentina}\\[-5pt] {\small sergiog@cab.cnea.gov.ar}\\[5pt] Edith \ Padr\'on\\[-5pt]{\small ULL-CSIC Geometr\'{\i}a Diferencial y Mec\'anica Geom\'etrica}\\[-5pt] {\small Departamento de Matem\'aticas, Estad{\'\i}stica e IO}\\[-5pt] {\small Universidad de la Laguna, La Laguna, Tenerife, Canary Islands,
Spain}\\[-5pt] {\small mepadron@ull.edu.es}}
\maketitle

\begin{abstract}
In this paper we develope, in a geometric framework, a Hamilton-Jacobi Theory
for general dynamical systems. Such a theory contains the classical theory for
Hamiltonian systems on a cotangent bundle and recent developments in the
framework of general symplectic, Poisson and almost-Poisson manifolds
(including some approaches to a Hamilton-Jacobi theory for nonholonomic
systems). Given a dynamical system, we show that every complete solution of
its related Hamilton-Jacobi Equation (HJE) gives rise to a set of first
integrals, and vice versa. From that, and in the context of symplectic and
Poisson manifolds, a deep connection between the HJE and the (non)commutative
integrability notion, and consequently the integrability by quadratures, is
stablished. Moreover, in the same context, we find conditions on the complete
solutions of the HJE that also ensures integrability by quadratures, but they
are weaker than those related to the (non)commutative integrability. Examples
are developed along all the paper in order to illustrate the theoretical results.

\end{abstract}

\section{Introduction}

In the context of Classical Mechanics, the \emph{standard} or \emph{classical}
(time-independent) Hamilton-Jacobi Theory was designed to construct, for a
given Hamiltonian system on a cotangent bundle $T^{\ast}Q$, local coordinates
such that the Hamilton equations expressed on these coordinates adopt a very
simple form. Here, by \emph{simple} we mean that such equations can be
integrated by \emph{quadratures }(i.e. its solutions can be given in terms of
primitives and inverses of known functions). The fundamental tool of the
theory is the so-called (time-independent) \emph{Hamilton-Jacobi equation
(HJE)} for a Hamiltonian function $H:T^{\ast}Q\rightarrow\mathbb{R}$. The
problem is to find a function $W$ on $Q$, known as the \emph{characteristic
Hamilton function}, such that the equation (the classical HJE)%
\[
H\left(  q,\frac{\partial W}{\partial q}\right)  =\mbox{constant}
\]
is satisfied (see for example \cite{am}). When we are given a family of
solutions $\left\{  W_{\lambda}\right\}  _{\lambda\in\Lambda}$ such that the
square matrix $\left.  \partial^{2}W_{\lambda}\right/  \partial q\partial
\lambda$ is non-degenerated, with $\Lambda$ an open subset of $\mathbb{R}%
^{\dim Q}$, the above mentioned coordinates can be constructed. More
precisely, a type $2$ canonical transformation (see \cite{gold}) can be
defined from the functions $W_{\lambda}$'s such that the equations of motions
of the system become, under such a transformation, simple enough to be solved
by quadratures.

In modern terms (see Ref. \cite{pepin-holo}), the classical HJE reads
$d\left(  H\circ\sigma\right)  =0$, and its unknown is a closed $1$-form
$\sigma:Q\rightarrow T^{\ast}Q$. If a solution $\sigma$ is given for a
Hamiltonian function $H$, the problem of finding the integral curves of its
Hamiltonian vector field $X_{H}$ (w.r.t. the canonical symplectic structure on
$T^{\ast}Q$), with initial conditions along the image of $\sigma$, reduces to
find the integral curves of the vector field $X_{H}^{\sigma}:=(\pi_{Q})_{\ast
}\circ X_{H}\circ\sigma$ on $Q$, where $\pi_{Q}:T^{\ast}Q\rightarrow Q$ is the
cotangent fibration. If one wants to find all the integral curves of $X_{H}$,
i.e. the solutions of the equations of motions of the system for every initial
condition, it is necessary to introduce the notion of a \emph{complete
solution}: a surjective local diffeomorphism $\Sigma:Q\times\Lambda\rightarrow
T^{\ast}Q$, with $\Lambda$ a manifold, such that $\sigma_{\lambda}%
:=\Sigma(\lambda,\cdot)$ is a solution of the classical HJE for each
$\lambda\in\Lambda$. In terms of $\Sigma$, around every point of $T^{\ast}Q$,
a local coordinate chart can be constructed such that the equations of motion
are easily solved. We can see that through the well-known connection between
the classical Hamilton-Jacobi Theory and the notion of \emph{commutative
integrability}. Let us describe such a connection (see, for example,
\cite{pepin-holo}). From every complete solution $\Sigma:Q\times
\Lambda\rightarrow T^{\ast}Q$, we can construct (unless locally) a Lagrangian
fibration $F:T^{\ast}Q\rightarrow\Lambda$ (transverse to $\pi_{Q}$) such that
the image of $X_{H}$ lives inside $\operatorname{Ker}F_{\ast}$. In other
terms, if $\Lambda$ is an open subset of $\mathbb{R}^{n}$ (with $n=\dim Q$),
we can build up $n$ first integrals $f_{i}:T^{\ast}Q\rightarrow\mathbb{R}$
(the components of $F$) which are independent and in involution with respect
to the canonical Poisson bracket. Conversely, from every Lagrangian fibration
$F:T^{\ast}Q\rightarrow\Lambda$ (transverse to $\pi_{Q}$) such that
$\operatorname{Im}X_{H}\subset\operatorname{Ker}F_{\ast}$, a complete solution
of the HJE can be constructed (unless locally). We recall that, if such a
fibration is given for a Hamiltonian system, the Arnold-Liouville Theorem
\cite{ar} establishes that the system in question is integrable by quadratures.

In the last few years, several generalized versions of the classical HJE have
been developed for Hamiltonian systems on different contexts: on general
symplectic, Poisson and almost-Poisson manifolds, and also on Lie algebroids
over vector bundles. The resulting Hamilton-Jacobi theories were applied to
nonholonomic systems, time-dependent Hamiltonian systems, reduced systems by
symmetries and systems with external forces \cite{bmmp, pepin-holo,
pepin-noholo, lmm, hjp}. In all of these contexts, a fibration $\Pi
:M\rightarrow N$ (i.e. a surjective submersion) is defined on the phase space
$M$ of each system, the solutions of the generalized HJE are sections
$\sigma:N\rightarrow M$ of such a fibration, and the complete solutions are
local diffeomorphisms $\Sigma:N\times\Lambda\rightarrow M$ such that
$\sigma_{\lambda}:=\Sigma(\lambda,\cdot)$ is a solution of the HJE for each
$\lambda\in\Lambda$. This clearly extends the classical situation, where the
involved fibration is the cotangent projection $\pi_{Q}:T^{\ast}Q\rightarrow
Q$ of a manifold $Q$. Unfortunately, no connection between complete solutions
and some kind of \emph{exact solvability} (as the integrability by
quadratures) has been given for any of those generalized versions of the HJE
(see for instance the Ref. \cite{bfs}).

\bigskip

This paper is the first of a series of papers in which we shall further extend
the previously mentioned Hamilton-Jacobi theories to general dynamical systems
on a fibered phase space. We are mainly interested in the connection between
complete solutions and exact solvability. In the present paper, as a first
step, we shall focus on (time-independent) Hamiltonian systems on general
Poisson manifolds. In forthcoming papers, we shall address the case of
Hamiltonian systems with external forces (including Hamiltonian systems with
constraints) and time-dependent Hamiltonian systems.

One of the contribution of this paper is to show, in the context of general
dynamical systems, that there exists a \emph{duality }between complete
solutions and first integrals, extending similar results that appear in the
literature. This enable us to establish, in the particular context of
Hamiltonian systems on Poisson manifolds, a deep connection between
\emph{(non)commutative integrability }and certain subclasses of complete solutions.

Recall that a Hamiltonian system on a Poisson manifold $M$, with Hamiltonian
function $H$, is a \emph{noncommutative integrable} system (see \cite{D,
fasso} for symplectic manifolds and \cite{rui,mf} for Poisson ones) if a
fibration $F:M\rightarrow\Lambda$ such that:

\begin{enumerate}
\item $\operatorname{Im}X_{H}\subset\operatorname{Ker}F_{\ast}$ ($F$ defines
first integrals for the system),

\item $\operatorname{Ker}F_{\ast}\subset\left(  \operatorname{Ker}F_{\ast
}\right)  ^{\perp}$, i.e. $F$ is isotropic,

\item $\left(  \operatorname{Ker}F_{\ast}\right)  ^{\perp}$ is integrable ($F$
has a polar),
\end{enumerate}
can be exhibited. In particular, the system is commutative integrable if in
addition $\operatorname{Ker}F_{\ast}=\left(  \operatorname{Ker}F_{\ast
}\right)  ^{\perp}$, i.e. $F$ is Lagrangian, as we have said above. By $\perp$
we are denoting the Poisson orthogonal. To be more precise, if $\Xi$ is the
Poisson bi-vector on $M$ and $\Xi^{\sharp}:T^{\ast}M\rightarrow TM$ is its
related linear bundle map, then $\left(  \operatorname{Ker}F_{\ast}\right)
^{\perp}:=\Xi^{\sharp}\left(  \left(  \operatorname{Ker}F_{\ast}\right)
^{0}\right)  $. All of these systems, among other things, are integrable by quadratures.

Another contribution of the paper (which can be seen as a first application of
our extended Hamilton-Jacobi Theory) is to show that conditions $1$ and $2$
listed above are enough in order to ensure integrability by quadratures of
Hamiltonian systems on Poisson manifolds. That is to say, the integrability of
$\left(  \operatorname{Ker}F_{\ast}\right)  ^{\perp}$ is not needed for that
purpose. Moreover, we show that condition $2$ can be replaced by a weaker
one:\emph{ } $\operatorname{Ker}F_{\ast}\cap\operatorname{Im}\Xi^{\sharp
}\subset\left(  \operatorname{Ker}F_{\ast}\right)  ^{\perp}$, which we call
\emph{weak isotropy} (together some regularity assumptions about the
symplectic leaves of $\Xi$). It is worth mentioning that the proof of this
result was done mainly in terms of complete solutions (instead of first
integrals). We think that, in spite of the duality between complete solutions
and first integrals, it would have been very hard to anticipate the mentioned
result by working with first integrals only.

\bigskip

The paper is organized as follows. In Section \ref{S2}, given a dynamical
system $\left(  M,X\right)  $ equipped with a fibration $\Pi:M\rightarrow N$,
being $M$ and $N$ smooth manifolds and $X$ a vector field on $M$, we introduce
the notion of $\Pi$-Hamilton-Jacobi equation ($\Pi$-HJE) for $\left(
M,X\right)  $. The unknown of such an equation is a section $\sigma
:N\rightarrow M$ of $\Pi$. The $\Pi$-HJE is defined in such a way that, if
$\sigma$ is a solution, then the vector field $X$ restricts to the closed
submanifold $\operatorname{Im}\sigma\subset M$ determined by the image of
$\sigma$. For the case in which $M$ is a Leibniz manifold\footnote{These
manifolds were introduced in \cite{gu} and they are used in the context of
generalized nonholonomic systems, gradient dynamical systems, in the study of
the interaction between non linear oscillators and the energy exchange between
them, and in the modeling of certain dissipative phenomena \cite{ortega}.} and
$X$ is a Hamiltonian vector field, we give characterizations of the
corresponding $\Pi$-HJEs in terms of (co)isotropic fibrations. In particular,
when $M$ is a symplectic manifold, we recover the well-known results for the
classical HJE (some of them contained in Ref. \cite{pepin-holo}). Moreover,
for Poisson and almost-Poisson manifolds, we recover the results of \cite{hjp}
and related works. In Section \ref{S3}, we define the complete solutions
$\Sigma:N\times\Lambda\rightarrow M$ of the $\Pi$-HJE for $\left(  M,X\right)
$, following the same ideas as in the previous works on the subject (see the
above discussion). In particular, for each $\lambda\in\Lambda$, the map
$\sigma_{\lambda}:=\Sigma(\lambda,\cdot):N\rightarrow M$ must be a solution of
the $\Pi$-HJE. This tool allows us to obtain any integral curve of $X$ from
the integral curves of the vector fields $X^{\sigma_{\lambda}}:=\Pi_{\ast
}\circ X\circ\sigma_{\lambda}$ on $N$. In the Leibniz scenario, several
characterizations and results about the complete solutions are given. In
particular, in the symplectic setting, we show that, if the fibration $\Pi$
and each submanifold $\operatorname{Im}\sigma_{\lambda}$ are Lagrangian, then
Darboux coordinates around every point of $M$, or equivalently local canonical
transformations, can be constructed in such a way that the solutions of the
Hamilton equations can be obtained by quadratures. This extends the
fundamental essence of the classical theory (valid for cotangent symplectic
manifolds and the cotangent fibration) to general symplectic manifolds and
general fibrations. Moreover, by choosing different Lagrangian fibrations
$\Pi$, we show that the mentioned canonical transformations can be taken of
type $1$, $2$, $3$ or $4$ (instead of type $2$ only, as in the classical
situation). We illustrate that with simple examples. Section \ref{S4} is
devoted to analyze the \textit{duality} between complete solutions and first
integrals. More precisely, we show that given a complete solution
$\Sigma:N\times\Lambda\rightarrow M$ of the $\Pi$-HJE for a dynamical system
$\left(  M,X\right)  $, for each $m\in M$ there exists an open neighborhood
$U$ of $m$ and a fibration $F:U\rightarrow F\left(  U\right)  \subset\Lambda$
(transverse to $\Pi$) such that $\operatorname{Im}\left.  X\right\vert _{U}%
\in\operatorname{Ker}F_{\ast}$. Reciprocally, given a fibration
$F:M\rightarrow\Lambda$ (transverse to $\Pi$) such that $\operatorname{Im}%
X\in\operatorname{Ker}F_{\ast},$ for each $m\in M$ there exists an open
neighborhood $U$ of $m$ such that $\left.  \left(  F,\Pi\right)  \right\vert
_{U}$ is invertible and its inverse defines a complete solution of the
$\left.  \Pi\right\vert _{U}$-HJE for $\left(  U,\left.  X\right\vert
_{U}\right)  $. Then, based on that duality, we prove an existence theorem for
complete solutions and, on the other hand, we establish a deep connection
between complete solutions of the $\Pi$-HJE and noncommutative integrability
on Poisson manifolds \cite{rui,mf}, i.e. a connection between complete
solutions and integrability by quadratures. Finally, in the last section of
the paper, we show for Hamiltonian systems on Poisson manifolds that, in order
to ensures integrability by quadratures, the above listed conditions $1$-$3$
that define a noncommutative system can be drastically weakened, as explained
in the previous paragraph.

\bigskip

\noindent We assume that the reader is familiar with the basic concepts of
Differential Geometry (see \cite{boot,kn,mrgm}), and with the basic ideas
related to the Lagrangian systems, Hamiltonian systems, Symplectic Geometry
and Poisson Geometry in the context of Geometric Mechanics (see
\cite{am,ar,mr}). We shall work in the smooth (i.e. $C^{\infty}$) category,
focusing exclusively on finite-dimensional smooth manifolds.

\section{The Hamilton-Jacobi equation for a fibration}

\label{S2}

In this section we introduce the notion of \emph{Hamilton-Jacobi equation} in
the context of general dynamical systems on a fibered phase space. We study
the consequences of having a solution of such an equation, focusing on how it
can help to find the trajectories of the system under consideration. We show
that our theory extends those developed for Hamiltonian systems on symplectic,
Poisson and almost-Poisson manifolds (including the nonholonomic systems). In
particular, in the framework of Hamiltonian systems on a cotangent bundle, we
show that our equation reduces to the classical Hamilton-Jacobi equation when
the cotangent fibration is considered.

\subsection{Definition and basic properties}

\label{hjemain}

Let us consider a dynamical system $\left(  M,X\right)  $, where $M$ is a
smooth manifold (the \emph{phase space}) and $X\in\mathfrak{X}\left(
M\right)  $ (i.e. $X$ is a vector field on $M$). Let $N$ be another smooth
manifold and $\Pi:M\rightarrow N$ a fibration, i.e. a surjective submersion
(\emph{ipso facto} an open map). Given an open subset $U\subset M$, we denote
by $\left.  \Pi\right\vert _{U}:U\rightarrow\Pi(U)$ the restricted fibration
and $\left.  X\right\vert _{U}:U\rightarrow TU$ the corresponding vector field
obtained by restricting $X$ to $U$.

\begin{definition}
\label{maindef}We shall call the equation
\begin{equation}
\sigma_{\ast}\circ\Pi_{\ast}\circ X\circ\sigma=X\circ\sigma\label{hje}%
\end{equation}
for a section $\sigma:N\rightarrow M$ of $\Pi$ (\emph{ipso facto} a closed
map), the $\Pi$-\textbf{Hamilton-Jacobi equation }($\Pi$-\textbf{HJE})\textbf{
for} $\left(  M,X\right)  $. If $\sigma$ satisfies such an equation, we shall
say that $\sigma$ is a \textbf{(global) solution of the }$\Pi$-\textbf{HJE
for} $\left(  M,X\right)  $. On the other hand, given an open subset $U\subset
M$, we shall say that a map $\sigma:\Pi\left(  U\right)  \rightarrow U$ is a
\textbf{local solution of the }$\Pi$-\textbf{HJE} \textbf{for} $\left(
M,X\right)  $ \textbf{on} $U$ if $\sigma$ is a solution of the $\left.
\Pi\right\vert _{U}$-HJE for $\left(  U,\left.  X\right\vert _{U}\right)  $.
\end{definition}

\begin{remark}
Given an open subset $V\subset N$, if a local section $\sigma:V\rightarrow M$
of $\Pi$ satisfies $\left(  \ref{hje}\right)  $ along $V$, then $\sigma$
defines (by co-restriction) a local solution of the $\Pi$-HJE for $\left(
M,X\right)  $ on every $U$ such that $\sigma\left(  V\right)  \subset
U\subset\Pi^{-1}\left(  V\right)  $ (\emph{ipso facto }$\Pi\left(  U\right)
=V$). Reciprocally, given an open subset $U\subset M$, if $\sigma:\Pi\left(
U\right)  \rightarrow U$ is a local solution of the $\Pi$-HJE for $\left(
M,X\right)  $ on $U$, then $\sigma$ defines a local section with domain
$V:=\Pi\left(  U\right)  $ which satisfies Eq. $\left(  \ref{hje}\right)  $
along $V$.
\end{remark}

\begin{example}
\label{fex}Suppose that $M=\mathbb{R}^{d}$, $N=\mathbb{R}^{k}$, with $k<d$,
and $\Pi:\mathbb{R}^{d}\rightarrow\mathbb{R}^{k}$ is the projection onto the
first $k$ components of $\mathbb{R}^{d}$. A section $\sigma:\mathbb{R}%
^{k}\rightarrow\mathbb{R}^{d}$ of $\Pi$ is a map of the form
\[
\sigma:n\in\mathbb{R}^{k}\longmapsto\left(  n,\hat{\sigma}\left(  n\right)
\right)  \in\mathbb{R}^{d}=\mathbb{R}^{k}\times\mathbb{R}^{l},
\]
with $l=d-k$. On the other hand, identifying $T\mathbb{R}^{d}$ with
$\mathbb{R}^{d}\times\mathbb{R}^{d}$, a vector field $X\in\mathfrak{X}\left(
\mathbb{R}^{d}\right)  $ can be described as a map%
\[
X:m\in\mathbb{R}^{d}\longmapsto\left(  m,\hat{X}\left(  m\right)  \right)
\in\mathbb{R}^{d}\times\mathbb{R}^{d}.
\]
Under this notation, it is easy to see that the $\Pi$-HJE for $\left(
M,X\right)  $ is given by%
\begin{equation}%
{\displaystyle\sum\nolimits_{i=1}^{k}}
\frac{\partial\hat{\sigma}_{j}\left(  n\right)  }{\partial n_{i}}\,\hat{X}%
_{i}\left(  n,\hat{\sigma}\left(  n\right)  \right)  =\hat{X}_{j+k}\left(
n,\hat{\sigma}\left(  n\right)  \right)  ,\ \ \ j=1,...,l, \label{achje}%
\end{equation}
where $n_{i}$, $\hat{\sigma}_{i}$ and $\hat{X}_{i}$ are the $i$-th components
of $n$, $\hat{\sigma}$ and $\hat{X}$, respectively.
\end{example}

It is clear that, for general manifolds $M$ and $N$, if we fix local
coordinates of $M$ and $N$ adapted\footnote{By \emph{adapted to }$\Pi$, we
mean local coordinates on $M$ and $N$ such that the local expression of $\Pi$
is the projection on the $k$-first coordinates, being $k:=\dim N$.} to $\Pi$,
the local expression of $\left(  \ref{hje}\right)  $ is exactly given by the
Eq. $\left(  \ref{achje}\right)  $. In particular, the local expression of the
$\Pi$-HJE for $\left(  M,X\right)  $ is a system of $l:=\dim M-\dim N$ first
order partial differential equations (PDE) with $l$ unknowns. So, a sufficient
condition for the existence of local solutions can be obtained from the
characteristic curve method. More precisely, in the notation of the above
example, if the vector field $X$ satisfies
\begin{equation}
\left(  \hat{X}_{1}\left(  m\right)  ,...,\hat{X}_{k}\left(  m\right)
\right)  \neq0,\ \ \text{i.e.\ \ \ }X\left(  m\right)  \notin
\operatorname{Ker}\Pi_{\ast,m}, \label{nonker0}%
\end{equation}
for some $m\in M$, it can be shown that there exists an open neighborhood
$U\subset M$ of $m$ and a local solution of the $\Pi$-HJE for $\left(
M,X\right)  $ on $U$. In Section \ref{exist}, we shall prove this fact as a
particular case of a more general result.

It is worth mentioning that the aim of this paper is not to study the
existence problem for global and local solutions of the $\Pi$-HJE (because
they use to exist around almost any point -see Section \ref{exist}-), but to
investigate what we can say about \emph{integrability} when a solution of such
an equation is given to us. By integrability we mean some kind of \emph{exact
solvability}, like the so-called \emph{solvability by quadratures}, i.e. the
possibility of finding an explicit expression of each trajectory of the system
by solving linear equations, using the inverse function theorem and
integrating known functions: \emph{the quadratures}.

\bigskip

Since the local solutions for $\left(  M,X\right)  $ are exactly the global
solutions for each $\left(  U,\left.  X\right\vert _{U}\right)  $, from now
on, we shall mainly focus, just for simplicity, on the global ones.

\bigskip

Let $\sigma:N\rightarrow M$ be a solution of the $\Pi$-HJE for $\left(
M,X\right)  $. Since $\sigma$ is a section of $\Pi$, then, by definition,
\begin{equation}
\Pi\circ\sigma=id_{N}\text{\ \ \ and\ \ \ }\sigma\circ\left.  \Pi\right\vert
_{\operatorname{Im}\sigma}=id_{\operatorname{Im}\sigma}. \label{secinv}%
\end{equation}
Consequently, $\sigma$ is an injective immersion and a closed map, so it is an
embedding. Then $\operatorname{Im}\sigma\subset M$ is a closed regular
submanifold and $N$ is diffeomorphic to $\operatorname{Im}\sigma$. Also,%
\begin{equation}
\operatorname{Ker}\left.  \Pi_{\ast}\right\vert _{\operatorname{Im}\sigma
}\oplus\operatorname{Im}\sigma_{\ast}=T_{\operatorname{Im}\sigma}M.
\label{deco}%
\end{equation}
These observations give rise to the following characterization of the
solutions of the $\Pi$-HJE.

\begin{proposition}
\label{Xims}A section $\sigma:N\rightarrow M$ of $\Pi$ is a solution of the
$\Pi$-HJE for $\left(  M,X\right)  $ if and only if .
\begin{equation}
\operatorname{Im}\left(  \left.  X\right\vert _{\operatorname{Im}\sigma
}\right)  \subset T\operatorname{Im}\sigma, \label{xins}%
\end{equation}
i.e. $X$ restricts to $\operatorname{Im}\sigma$.
\end{proposition}

\emph{Proof.} Given $m\in\operatorname{Im}\sigma$, we can write $m=\sigma
\left(  n\right)  $ for a unique $n\in N$. It follows from Eq. $\left(
\ref{hje}\right)  $ that%
\[
X\left(  m\right)  =X\left(  \sigma\left(  n\right)  \right)  =\sigma_{\ast
,n}\circ\Pi_{\ast,\sigma\left(  n\right)  }\left(  X\left(  \sigma\left(
n\right)  \right)  \right)  \in T_{\sigma\left(  n\right)  }\operatorname{Im}%
\sigma=T_{m}\operatorname{Im}\sigma,
\]
so $X$ restricts to $\operatorname{Im}\sigma$. Reciprocally, suppose that
$\sigma:N\rightarrow M$ is a section of $\Pi$ such that Eq. $\left(
\ref{xins}\right)  $ holds. Combining this last equation with $\left(
\ref{secinv}\right)  $, we have that%
\[
\sigma_{\ast,n}\left[  \Pi_{\ast,\sigma\left(  n\right)  }\left(  X\left(
\sigma\left(  n\right)  \right)  \right)  \right]  =\sigma_{\ast,n}\left[
\left(  \left.  \Pi\right\vert _{\operatorname{Im}\sigma}\right)
_{\ast,\sigma\left(  n\right)  }\left(  X\left(  \sigma\left(  n\right)
\right)  \right)  \right]  =X\left(  \sigma\left(  n\right)  \right)
,\ \ \ \ \forall n\in N,
\]
i.e. $\sigma\ $satisfies the Eq. $\left(  \ref{hje}\right)  $.$\ \ \ \square$

\bigskip

Since $\operatorname{Im}\sigma$ is closed, above result says, as it is
well-known, that $\operatorname{Im}\sigma$ is an $X$-\textbf{invariant
submanifold} of $M$. In other words, an integral curve $\Gamma:I\rightarrow M$
of $X$, with $I$ an open interval, intersects $\operatorname{Im}\sigma$ if and
only if $\operatorname{Im}\Gamma\subset\operatorname{Im}\sigma$. Also, Eq.
$\left(  \ref{hje}\right)  $ for $\sigma$ says exactly that the vector fields
$X\in\mathfrak{X}\left(  M\right)  $ and
\begin{equation}
X^{\sigma}:=\Pi_{\ast}\circ X\circ\sigma\in\mathfrak{X}\left(  N\right)
\label{Xs}%
\end{equation}
are $\sigma$-related, i.e.%
\begin{equation}
X\circ\sigma=\sigma_{\ast}\circ X^{\sigma}. \label{xsrel}%
\end{equation}
Then, if $\gamma$ is an integral curve of $X^{\sigma}$, the curve
$\Gamma:=\sigma\circ\gamma$ is an integral curve of $X$ (intersecting
$\operatorname{Im}\sigma$). Moreover, every integral curve of $X$ intersecting
$\operatorname{Im}\sigma$ can be obtained in this way, as we prove below.

\begin{theorem}
\label{parcial}Let $\sigma:N\rightarrow M$ be a solution of the $\Pi$-HJE for
$\left(  M,X\right)  $. If $\Gamma:I\rightarrow M$ is an integral curve of $X$
intersecting $\operatorname{Im}\sigma$, then we can write $\Gamma=\sigma
\circ\gamma$ for a unique integral curve $\gamma$ of $X^{\sigma}$.
\end{theorem}

\emph{Proof.} The $X$-invariance of $\operatorname{Im}\sigma$ ensures that
$\operatorname{Im}\Gamma\subset\operatorname{Im}\sigma$, and consequently
$\Pi\left(  \Gamma\left(  t\right)  \right)  =\left.  \Pi\right\vert
_{\operatorname{Im}\sigma}\left(  \Gamma\left(  t\right)  \right)  $, for all
$t\in I$. So, the second part of Eq. $\left(  \ref{secinv}\right)  $ implies
the identity $\sigma\circ\Pi\circ\Gamma=\Gamma$. Defining $\gamma:=\Pi
\circ\Gamma$, we have that $\sigma\circ\gamma=\Gamma$ and from $\left(
\ref{hje}\right)  $ that
\[
\gamma^{\prime}\left(  t\right)  =\Pi_{\ast}\left(  \Gamma^{\prime}\left(
t\right)  \right)  =\Pi_{\ast}\circ X\left(  \Gamma\left(  t\right)  \right)
=\Pi_{\ast}\circ X\circ\sigma\circ\gamma\left(  t\right)  =X^{\sigma}\left(
\gamma\left(  t\right)  \right)  ,
\]
$\forall t\in I$, i.e. $\gamma$ is an integral curve of $X^{\sigma}$. The
injectivity of $\sigma$ ensures that such a curve $\gamma$ is
unique.\ \ $\square$

\bigskip

\begin{remark}
\label{dkl}It is easy to see in Example \ref{fex} that the vector field
$X^{\sigma}\in\mathfrak{X}\left(  N\right)  $ is given by
\[
X^{\sigma}\left(  n\right)  =\left(  n,\hat{X}_{1}\left(  \sigma\left(
n\right)  \right)  ,...,\hat{X}_{k}\left(  \sigma\left(  n\right)  \right)
\right)  .
\]
So, in order to find a trajectory of $\left(  M,X\right)  $ starting in
$\operatorname{Im}\sigma$, instead of solving the system of $d$ first order
ordinary differential equations (ODE)
\[
\dot{m}_{i}\left(  t\right)  =\hat{X}_{i}\left(  m\left(  t\right)  \right)
,\ \ \ i=1,...,d,
\]
for a curve $m:I\rightarrow M$, being $I\subset\mathbb{R}$ an open interval,
it is enough to solve the system of $k$ first order ODEs
\[
\dot{n}_{i}\left(  t\right)  =\hat{X}_{i}\left(  \sigma\left(  n\left(
t\right)  \right)  \right)  ,\ \ \ i=1,...,k,
\]
for $n:I\rightarrow N$, and then define $m\left(  t\right)  =\sigma\left(
n\left(  t\right)  \right)  $.
\end{remark}

\subsection{The symplectic scenario}

\label{syscen}

In the case of a Hamiltonian system defined on a symplectic manifold $\left(
M,\omega\right)  $, we have other characterizations of the extended HJE. As
usual, given $m\in M$ and a subspace $V\subset T_{m}M$, we denote by
$V^{\omega}\subset T_{m}M$ the symplectic orthogonal of $V$, and we say that
$V$ is \textbf{isotropic} (resp. \textbf{co-isotropic}) if $V\subset
V^{\omega}$ (resp. $V^{\omega}\subset V$), and \textbf{Lagrangian}\emph{ }if
$V=V^{\omega}$. Also, by $\omega^{\flat}:TM\rightarrow T^{\ast}M$ we denote
the vector bundle isomorphism given by $\left\langle \omega^{\flat}\left(
X\right)  ,Y\right\rangle =\omega\left(  X,Y\right)  $, and by $\omega
^{\sharp}:T^{\ast}M\rightarrow TM$ we denote its inverse.

\begin{definition}
\label{paral}We shall say that a submersion $\Pi:M\rightarrow N$ is
\textbf{isotropic} (resp. \textbf{co-isotropic}) if each fiber of
$\operatorname{Ker}\Pi_{\ast}$ is isotropic (resp. co-isotropic), i.e.
\begin{equation}
\operatorname{Ker}\Pi_{\ast}\subset\left(  \operatorname{Ker}\Pi_{\ast
}\right)  ^{\omega}, \label{isot}%
\end{equation}
(resp. $\left(  \operatorname{Ker}\Pi_{\ast}\right)  ^{\omega}\subset
\operatorname{Ker}\Pi_{\ast}$). $\Pi$ is \textbf{Lagrangian} if it is both
isotropic and co-isotropic. We shall say that a section $\sigma:N\rightarrow
M$ of $\Pi$ is \textbf{isotropic} (resp. \textbf{co-isotropic}) if each fiber
of $\operatorname{Im}\sigma_{\ast}$ is isotropic (resp. co-isotropic), i.e.%
\[
\operatorname{Im}\sigma_{\ast}\subset\left(  \operatorname{Im}\sigma_{\ast
}\right)  ^{\omega}%
\]
(resp. $\left(  \operatorname{Im}\sigma_{\ast}\right)  ^{\omega}%
\subset\operatorname{Im}\sigma_{\ast}$), and that $\sigma$ is
\textbf{Lagrangian} if it is both isotropic and co-isotropic.
\end{definition}

It is easy to show that $\sigma$ is isotropic if and only if $\sigma^{\ast
}\omega=0$. Note that isotropy (resp. co-isotropy) condition on $\sigma$
implies that $2\dim N\leq\dim M$ (resp. $\dim M\leq2\dim N$).

\begin{theorem}
\label{phjsim}Consider a symplectic manifold $\left(  M,\omega\right)  $, a
function $H:M\rightarrow\mathbb{R}$, its Hamiltonian vector field $X_{H}$
w.r.t. $\omega$, a fibration $\Pi:M\rightarrow N$ and a section $\sigma
:N\rightarrow M$ of $\Pi$.

\begin{enumerate}
\item If $\sigma$ is a solution of the $\Pi$-HJE for $\left(  M,X_{H}\right)
$, then
\begin{equation}
d\left(  H\circ\sigma\right)  =i_{X_{H}^{\sigma}}\sigma^{\ast}\omega,
\label{hjes2}%
\end{equation}
being $X_{H}^{\sigma}$ the vector field on $N$ defined as in (\ref{Xs}), for
$X=X_{H}$.

\item If in addition $\sigma$ is isotropic, then
\begin{equation}
d\left(  H\circ\sigma\right)  =0. \label{hjes5}%
\end{equation}

\item On the other hand, if $\sigma$ satisfies $\left(  \ref{hjes2}\right)  $,
then
\[
\operatorname{Im}\left(  X_{H}\circ\sigma-\sigma_{\ast}\circ X_{H}^{\sigma
}\right)  \subset\left(  \operatorname{Ker}\Pi_{\ast}\right)  \cap\left(
\operatorname{Im}\sigma_{\ast}\right)  ^{\omega}.
\]

\end{enumerate}
\end{theorem}

\emph{Proof.} $\left(  1\right)  $ If $\sigma$ is a solution of the $\Pi$-HJE
for $\left(  M,X_{H}\right)  $, then [see Eqs. $\left(  \ref{Xs}\right)  $ and
$\left(  \ref{xsrel}\right)  $]
\begin{equation}
X_{H}\circ\sigma=\sigma_{\ast}\circ X_{H}^{\sigma}. \label{o1}%
\end{equation}
Using that $X_{H}=\omega^{\sharp}\circ dH$, for all $n\in N$ and $y\in T_{n}%
N$, we have that
\begin{align}
\omega\left(  X_{H}\left(  \sigma\left(  n\right)  \right)  ,\sigma_{\ast
,n}\left(  y\right)  \right)   &  =\left\langle dH\left(  \sigma\left(
n\right)  \right)  ,\sigma_{\ast,n}\left(  y\right)  \right\rangle \nonumber\\
&  =\left\langle \sigma^{\ast}dH\left(  n\right)  ,y\right\rangle
=\left\langle d\left(  H\circ\sigma\right)  \left(  n\right)  ,y\right\rangle
. \label{lhjm}%
\end{align}
On the other hand,
\begin{equation}
\omega\left(  \sigma_{\ast,n}\left(  X_{H}^{\sigma}\left(  n\right)  \right)
,\sigma_{\ast,n}\left(  y\right)  \right)  =\sigma^{\ast}\omega\left(
X_{H}^{\sigma}\left(  n\right)  ,y\right)  . \label{o2}%
\end{equation}
Using $\left(  \ref{o1}\right)  $ on the first members of $\left(
\ref{lhjm}\right)  $ and $\left(  \ref{o2}\right)  $, we deduce $\left(
\ref{hjes2}\right)  $.

\smallskip

$\left(  2\right)  $ The isotropy condition for $\sigma$ means that
$\sigma^{\ast}\omega=0$. Then, Eq. $\left(  \ref{hjes2}\right)  $ translates
to $\left(  \ref{hjes5}\right)  $.

\smallskip

$\left(  3\right)  $ Given a section $\sigma$ of $\Pi$, it is clear that
[recall Eqs. $\left(  \ref{secinv}\right)  $ and $\left(  \ref{Xs}\right)  $]
$\Pi_{\ast}\circ\left(  X_{H}\circ\sigma-\sigma_{\ast}\circ X_{H}^{\sigma
}\right)  =0$, i.e.
\[
\operatorname{Im}\left(  X_{H}\circ\sigma-\sigma_{\ast}\circ X_{H}^{\sigma
}\right)  \subset\operatorname{Ker}\Pi_{\ast}.
\]
On the other hand, if $\sigma$ satisfies $\left(  \ref{hjes2}\right)  $, from
$\left(  \ref{lhjm}\right)  $ and $\left(  \ref{o2}\right)  $ we have that
\[
\omega\left(  X_{H}\left(  \sigma\left(  n\right)  \right)  -\sigma_{\ast
,n}\left(  X_{H}^{\sigma}\left(  n\right)  \right)  ,\sigma_{\ast,n}\left(
y\right)  \right)  =0
\]
for all $n\in N$ and $y\in T_{n}N$, i.e. $\operatorname{Im}\left(  X_{H}%
\circ\sigma-\sigma_{\ast}\circ X_{H}^{\sigma}\right)  \subset\left(
\operatorname{Im}\sigma_{\ast}\right)  ^{\omega}$, and the Theorem is
proved.\ \ \ \ $\ \square$

\bigskip

Let us emphasize that, in general, we can not ensure that the $\Pi$-HJE for a
Hamiltonian system on a symplectic manifold is equivalent to Eq. $\left(
\ref{hjes2}\right)  $. Some conditions that enable us to do that are given below.

\begin{corollary}
\label{coro1}Under the conditions of Theorem \ref{phjsim}, if in addition
the fibration $\Pi$ is isotropic (or in particular if $\Pi$ is Lagrangian), we
have that $\sigma$ is a solution of the $\Pi$-HJE for $\left(  M,X_{H}\right)
$ if and only if $\left(  \ref{hjes2}\right)  $ holds.
\end{corollary}

\emph{Proof.} The first implication was shown in the point $\left(  1\right)
$ of Theorem \ref{phjsim}. We have to show the converse. Since $\Pi$ is
isotropic, then [see Eqs. $\left(  \ref{deco}\right)  $ and $\left(
\ref{isot}\right)  $], for all $m\in\operatorname{Im}\sigma$,
\begin{align*}
\operatorname{Ker}\Pi_{\ast,m}\cap\left(  \operatorname{Im}\sigma_{\ast
,\Pi\left(  m\right)  }\right)  ^{\omega}  &  \subset\left(
\operatorname{Ker}\Pi_{\ast,m}\right)  ^{\omega}\cap\left(  \operatorname{Im}%
\sigma_{\ast,\Pi\left(  m\right)  }\right)  ^{\omega}\\
&  =\left(  \operatorname{Ker}\Pi_{\ast,m}+\operatorname{Im}\sigma_{\ast
,\Pi\left(  m\right)  }\right)  ^{\omega}=\left(  T_{m}M\right)  ^{\omega}=0,
\end{align*}
and accordingly [see the item $\left(  3\right)  $ of Theorem \ref{phjsim}]
$X_{H}\circ\sigma-\sigma_{\ast}\circ X_{H}^{\sigma}=0$, i.e. $\sigma$ is a
solution of the $\Pi$-HJE for $\left(  M,X_{H}\right)  $ [see Eq. $\left(
\ref{xsrel}\right)  $ for $X=X_{H}$].$\ \ \ \ \square$

\bigskip

\begin{example}
\label{eglhje}Consider a Lagrangian system on a manifold $Q$ with regular
Lagrangian function $L:TQ\rightarrow\mathbb{R}$, and let $E_{L}:TQ\rightarrow
\mathbb{R}$ be the \emph{energy function }associated to $L$ (see \cite{am,ar}
and references therein). Regularity of $L$ implies that its \emph{Cartan }%
$2$\emph{-form} $\omega_{L}$ is a symplectic form on $TQ$, and consequently
the trajectories of the system are given by the integral curves of the vector
field
\begin{equation}
X_{L}=\omega_{L}^{\sharp}\left(  dE_{L}\right)  . \label{xL}%
\end{equation}
So, we have the dynamical system $\left(  TQ,X_{L}\right)  $. Consider the
canonical projection $\tau_{Q}:TQ\rightarrow Q$. The $\tau_{Q}$-HJE for
$\left(  TQ,X_{L}\right)  $ is given by%
\[
\sigma_{\ast}\circ\left(  \tau_{Q}\right)  _{\ast}\circ X_{L}\circ\sigma
=X_{L}\circ\sigma,
\]
being $\sigma:Q\rightarrow TQ$ a section of $\tau_{Q}$ [note that $\sigma
\in\mathfrak{X}\left(  Q\right)  $]. As it is well-known, $X_{L}$ is a section
of $\left(  \tau_{Q}\right)  _{\ast}$, i.e. $X_{L}$ is a \emph{second order
vector field}, then%
\begin{equation}
X_{L}^{\sigma}:=\left(  \tau_{Q}\right)  _{\ast}\circ X_{L}\circ\sigma=\sigma,
\label{xsL}%
\end{equation}
and accordingly the $\tau_{Q}$-HJE for $\left(  TQ,X_{L}\right)  $ reduces to%
\begin{equation}
\sigma_{\ast}\circ\sigma=X_{L}\circ\sigma. \label{glhje}%
\end{equation}
Last equation was introduced in Ref. \cite{pepin-holo} under the name of
\emph{generalized Lagrangian Hamilton-Jacobi problem}. Since $\tau_{Q}$ is a
Lagrangian fibration, Corollary \ref{coro1} says that Eq. $\left(
\ref{glhje}\right)  $ is equivalent to $d\left(  E_{L}\circ\sigma\right)
=i_{\sigma}\sigma^{\ast}\omega_{L}$ [recall Eq. $\left(  \ref{xsL}\right)  $],
as it was also shown in \cite{pepin-holo}.
\end{example}

The following result tells us what happen in the case when the section
$\sigma$ is co-isotropic or Lagrangian.

\begin{corollary}
\label{coro2}Under the conditions of Theorem \ref{phjsim}, if in addition

\begin{enumerate}
\item $\sigma$ is co-isotropic, then $\sigma$ is a solution of the $\Pi$-HJE
for $\left(  M,X_{H}\right)  $ if and only if $d\left(  H\circ\sigma\right)
=i_{X_{H}^{\sigma}}\sigma^{\ast}\omega$.

\item $\sigma$ is Lagrangian, then $\sigma$ is a solution of the $\Pi$-HJE for
$\left(  M,X_{H}\right)  $ if and only if $d\left(  H\circ\sigma\right)  =0$.
\end{enumerate}
\end{corollary}

\emph{Proof.} $\left(  1\right)  $ As in the previous corollary, we must show
that, for all $m\in\operatorname{Im}\sigma$,
\[
\operatorname{Ker}\Pi_{\ast,m}\cap\left(  \operatorname{Im}\sigma_{\ast
,\Pi\left(  m\right)  }\right)  ^{\omega}=0.
\]
From $\left(  \ref{deco}\right)  $, it follows that%
\[
\left(  \operatorname{Ker}\left.  \Pi_{\ast}\right\vert _{\operatorname{Im}%
\sigma}\right)  ^{\omega}\oplus\left(  \operatorname{Im}\sigma_{\ast}\right)
^{\omega}=T_{\operatorname{Im}\sigma}M.
\]
So, using the co-isotropy of $\sigma$, for all $m\in\operatorname{Im}\sigma$,%
\begin{align*}
\operatorname{Ker}\Pi_{\ast,m}\cap\left(  \operatorname{Im}\sigma_{\ast
,\Pi\left(  m\right)  }\right)  ^{\omega}  &  \subset\operatorname{Ker}%
\Pi_{\ast,m}\cap\operatorname{Im}\sigma_{\ast,\Pi\left(  m\right)  }\\
&  =\left(  \left(  \operatorname{Ker}\Pi_{\ast,m}\right)  ^{\omega}+\left(
\operatorname{Im}\sigma_{\ast,\Pi\left(  m\right)  }\right)  ^{\omega}\right)
^{\omega}=\left(  T_{m}M\right)  ^{\omega}=0,
\end{align*}
as we wanted to show.\ \ \ 

$\smallskip$

$\left(  2\right)  $ Combining the point $\left(  1\right)  $ of this
corollary and the point $\left(  2\right)  $ of Theorem \ref{phjsim}, the
result easily follows.$\ \ \ \square$

\smallskip

\begin{example}
Coming back to Example \ref{eglhje}, if we ask $\sigma$ to be isotropic w.r.t.
$\omega_{L}$, then $\sigma$ is Lagrangian, for dimensional reasons, and
consequently it must satisfy $d\left(  E_{L}\circ\sigma\right)  =0$ and
$\sigma^{\ast}\omega_{L}=0$. These equations were called the \emph{Lagrangian
Hamilton-Jacobi problem }in Ref. \cite{pepin-holo}.
\end{example}

\subsubsection{The standard Hamilton-Jacobi Theory}

\label{canej}

In this section we shall see how to recover the classical HJE. Suppose that
$M=T^{\ast}Q$ for some manifold $Q$, $\omega=\omega_{Q}=d\theta$, where
$\theta$ is the Liouville $1$-form (i.e. $\omega$ is the canonical $2$-form on
$T^{\ast}Q$), $H:T^{\ast}Q\rightarrow\mathbb{R}$ is a function and $\Pi
=\pi_{Q}:T^{\ast}Q\rightarrow Q$ is the canonical cotangent projection. Note
that the sections of $\pi_{Q}$ are the $1$-forms on $Q$ [i.e. they belong to
$\Omega^{1}\left(  Q\right)  $], $\pi_{Q}$ is a Lagrangian fibration,
and\footnote{Given a function $f:T^{\ast}Q\rightarrow\mathbb{R}$, by
$\mathbb{F}f:T^{\ast}Q\rightarrow TQ$ we are denoting the \emph{fiber
derivative of }$f$, given by
\begin{equation}
\left\langle \alpha,\mathbb{F}f\left(  \beta\right)  \right\rangle =\left.
\frac{d}{ds}\right\vert _{0}f\left(  \alpha+s\,\beta\right)  ,\ \ \ \forall
q\in Q,\ \ \forall\alpha,\,\beta\in T_{q}^{\ast}Q. \label{fd}%
\end{equation}
}%
\begin{equation}
\pi_{Q\ast}\circ X_{H}=\mathbb{F}H. \label{pifh}%
\end{equation}
Then, the $\pi_{Q}$-HJE for $\left(  T^{\ast}Q,X_{H}\right)  $ is an equation
for a $1$-form $\sigma\in\Omega^{1}\left(  Q\right)  $ and it reads%
\begin{equation}
\sigma_{\ast}\circ\mathbb{F}H\circ\sigma=X_{H}\circ\sigma. \label{qshje}%
\end{equation}
Moreover, $X_{H}^{\sigma}=\mathbb{F}H\circ\sigma$ and, since $\pi_{Q}$ is
Lagrangian, Corollary \ref{coro1} says that $\left(  \ref{qshje}\right)  $ is
equivalent to [see Eq. $\left(  \ref{hjes2}\right)  $]%
\[
d\left(  H\circ\sigma\right)  =i_{\mathbb{F}H\circ\sigma}\sigma^{\ast}\omega.
\]
Recall that $\sigma^{\ast}\theta=\sigma$ for every $\sigma\in\Omega^{1}\left(
Q\right)  $, and consequently%
\begin{equation}
\sigma^{\ast}\omega=\sigma^{\ast}d\theta=d\sigma^{\ast}\theta=d\sigma.
\label{ests}%
\end{equation}
Hence, the $\pi_{Q}$-HJE for $\left(  T^{\ast}Q,X_{H}\right)  $ can be written%
\begin{equation}
d\left(  H\circ\sigma\right)  =i_{\mathbb{F}H\circ\sigma}d\sigma. \label{ghje}%
\end{equation}
In Ref. \cite{pepin-holo}, this last equation was called \emph{generalized
Hamiltonian Hamilton-Jacobi problem}.

\begin{remark}
\label{kozl}Another extension of the Hamilton-Jacobi Theory was presented in
Reference \cite{koz}. Its related (time-independent) Hamilton-Jacobi equation
is given precisely by the Eq. $\left(  \ref{ghje}\right)  $ (see the Equation
$1.6$ of that reference). It was called \emph{Lamb's equation} in \cite{koz}.
\end{remark}

The \emph{classical} case is obtained when we look for a solution $\sigma$ of
$\left(  \ref{ghje}\right)  $ such that $d\sigma=0$. In this case, Eq.
$\left(  \ref{ghje}\right)  $ translates to\footnote{This also can be seen in
the following way. The condition $d\sigma=0$ says exactly that $\sigma$ is a
Lagrangian section. Thus, according to Corollary \ref{coro2}, the $\pi_{Q}%
$-HJE $\left(  \ref{qshje}\right)  $ reduces precisely to $\left(
\ref{ehje}\right)  $.}%
\begin{equation}
d\left(  H\circ\sigma\right)  =0\text{\ \ \ \ and\ \ \ \ }d\sigma=0,
\label{ehje}%
\end{equation}
which defines the so-called \emph{Hamiltonian Hamilton-Jacobi problem} of Ref.
\cite{pepin-holo}.

\begin{definition}
\label{standd}From now on, we shall mean by \textbf{standard} (resp.
\textbf{standard Lagrangian or classical})\textbf{ situation} the case in
which $M=T^{\ast}Q$ for some manifold $Q$, $\omega$ is its canonical
symplectic form, $X=X_{H}$ for some function $H:T^{\ast}Q\rightarrow
\mathbb{R}$, $N=Q$ and $\Pi=\pi_{Q}$ (resp. and $d\sigma=0$). And we shall
refer to $\left(  \ref{ghje}\right)  $ [resp. $\left(  \ref{ehje}\right)  $]
as the \textbf{standard} (resp. \textbf{standard Lagrangian or classical}%
)\textbf{ HJE}.
\end{definition}

Every solution of $\left(  \ref{ehje}\right)  $ is locally given by a
primitive function $W:V\rightarrow\mathbb{R}$, i.e. $\left.  \sigma\right\vert
_{V}=dW$, where $V$ is an open submanifold of $Q$. The function $W$ is usually
called \textbf{characteristic Hamilton function}.

\subsubsection{Some examples: standard and nonstandard equations}

\label{stanon}

Let us take $M:=\mathbb{R}^{k}\times\mathbb{R}^{k}$ and the symplectic
structure $\omega:=dq^{i}\wedge dp_{i}$ (summation over repeated indices is
assumed from now on), being $\left(  q,p\right)  $ the natural global
coordinates of $\mathbb{R}^{k}\times\mathbb{R}^{k}$. Note that $\theta
=p_{i}\,dq^{i}$. Fix $N=\mathbb{R}^{k}$ and consider the submersions
$\mathfrak{p}_{r}:\mathbb{R}^{k}\times\mathbb{R}^{k}\rightarrow\mathbb{R}^{k}%
$, $r=1,2$, given by the projectors on the first and second $\mathbb{R}^{k}$
factors, respectively. It is clear that $\operatorname{Ker}\left(
\mathfrak{p}_{r}\right)  _{\ast}$ is Lagrangian for $r=1,2$. As a consequence,
given any function $H:\mathbb{R}^{k}\times\mathbb{R}^{k}\rightarrow\mathbb{R}%
$, the $\mathfrak{p}_{r}$-HJE for $\left(  \mathbb{R}^{k}\times\mathbb{R}%
^{k},X_{H}\right)  $ is given by
\begin{equation}
d\left(  H\circ\sigma\right)  =i_{X_{H}^{\sigma,r}}\sigma^{\ast}\omega,
\label{hjei}%
\end{equation}
where $\sigma:\mathbb{R}^{k}\rightarrow\mathbb{R}^{k}\times\mathbb{R}^{k}$ is
a section of $\mathfrak{p}_{r}$ and $X_{H}^{\sigma,r}=\left(  \mathfrak{p}%
_{r}\right)  _{\ast}\circ X_{H}\circ\sigma$. It is also clear that the $r=1$
case corresponds to the standard HJE [see Eq. $\left(  \ref{ghje}\right)  $].
In fact, under the identifications mentioned above, $X_{H}^{\sigma
,1}=\mathbb{F}H\circ\sigma$ and $\sigma^{\ast}\omega=d\sigma$, so Eq. $\left(
\ref{hjei}\right)  $ reduces to Eq. $\left(  \ref{ghje}\right)  $. Writing
$\sigma\left(  q\right)  =\left(  q,\hat{\sigma}\left(  q\right)  \right)  $,
it is easy to see that Eq. $\left(  \ref{hjei}\right)  $ results

\smallskip%

\begin{equation}
\underset{d\left(  H\circ\sigma\right)  }{\underbrace{\frac{\partial H\left(
q,\hat{\sigma}\left(  q\right)  \right)  }{\partial q^{i}}}}=\underset
{\sigma^{\ast}\omega=d\sigma}{\underbrace{\left(  \frac{\partial\hat{\sigma
}_{j}\left(  q\right)  }{\partial q^{i}}-\frac{\partial\hat{\sigma}_{i}\left(
q\right)  }{\partial q^{j}}\right)  }}\underset{X_{H}^{\sigma,1}%
=\mathbb{F}H\circ\sigma}{\,\underbrace{\frac{\partial H}{\partial p_{j}%
}\left(  q,\hat{\sigma}\left(  q\right)  \right)  }},\ \ \ \ \ i=1,...,k,
\label{lshje}%
\end{equation}
where $\hat{\sigma}_{i}$ is the $i$-th component of $\hat{\sigma}$.

\begin{remark}
Since we can write $X_{H}\left(  q,p\right)  =\left(  q,p,\hat{X}_{H}\left(
q,p\right)  \right)  $ with%
\[
\hat{X}_{H}\left(  q,p\right)  =\left(  \frac{\partial H}{\partial p_{1}%
}\left(  q,p\right)  ,...,\frac{\partial H}{\partial p_{k}}\left(  q,p\right)
,-\frac{\partial H}{\partial q^{1}}\left(  q,p\right)  ,...,-\frac{\partial
H}{\partial q^{k}}\left(  q,p\right)  \right)  ,
\]
Eq. $\left(  \ref{lshje}\right)  $ can be derived from the PDEs $\left(
\ref{achje}\right)  $ for $d=2l=2k$ and
\[
\hat{X}_{j}\left(  q,p\right)  =\left\{
\begin{array}
[c]{cc}%
\displaystyle\frac{\partial H}{\partial p_{j}}\left(  q,p\right)  , & 1\leq
j\leq k,\\
& \\
\displaystyle-\frac{\partial H}{\partial q^{j-k}}\left(  q,p\right)  , &
k+1\leq j\leq2k.
\end{array}
\right.
\]

\end{remark}

For the $r=2$ case, writing $\sigma\left(  p\right)  =\left(  \hat{\sigma
}\left(  p\right)  ,p\right)  $, Eq. $\left(  \ref{hjei}\right)  $ is given by

\smallskip\
\begin{equation}
\underset{d\left(  H\circ\sigma\right)  }{\underbrace{\frac{\partial H\left(
\hat{\sigma}\left(  p\right)  ,p\right)  }{\partial p_{i}}}}=\underset
{\sigma^{\ast}\omega}{\underbrace{\left(  \frac{\partial\hat{\sigma}%
^{j}\left(  p\right)  }{\partial p_{i}}-\frac{\partial\hat{\sigma}^{i}\left(
p\right)  }{\partial p_{j}}\right)  }}\underset{X_{H}^{\sigma,2}%
}{\,\underbrace{\frac{\partial H}{\partial q^{j}}\left(  \hat{\sigma}\left(
p\right)  ,p\right)  }}\,,\ \ \ \ \ i=1,...,k. \label{lnshje}%
\end{equation}
Note that isotropy condition $\sigma^{\ast}\omega=0$ reads%
\begin{equation}
\partial_{i}\hat{\sigma}^{j}-\partial_{j}\hat{\sigma}^{i}=0 \label{isocon}%
\end{equation}
in both cases. Eq. $\left(  \ref{lnshje}\right)  $ is an example of a
nonstandard HJE. Let us focus on the $k=1$ case. The Equations $\left(
\ref{lshje}\right)  $ and $\left(  \ref{lnshje}\right)  $ translate to%
\[
\frac{\partial H\left(  q,\hat{\sigma}\left(  q\right)  \right)  }{\partial
q}=0\ \ \ \text{and\ \ \ }\frac{\partial H\left(  \hat{\sigma}\left(
p\right)  ,p\right)  }{\partial p}=0,
\]
respectively, since the second member is identically zero in both of the
equations. Suppose that $H:\mathbb{R}\times\mathbb{R}\rightarrow\mathbb{R}$ is
defined by%
\[
H\left(  q,p\right)  =\frac{1}{2}\,\left(  p^{2}+f\left(  q\right)  \right)
,
\]
for some function $f:\mathbb{R}\rightarrow\mathbb{R}$. Then, $\left(
\ref{lshje}\right)  $ and $\left(  \ref{lnshje}\right)  $ reduce to the ODEs%
\[
\frac{1}{2}\left[  \left(  \hat{\sigma}^{2}\left(  q\right)  \right)
^{\prime}+f^{\prime}\left(  q\right)  \right]  =0\ \ \ \text{and\ \ \ }%
p+\frac{1}{2}\,\left[  f\left(  \hat{\sigma}\left(  p\right)  \right)
\right]  ^{\prime}=0.
\]
For the first one, the solutions are given by the formulae%
\[
\hat{\sigma}_{\lambda}^{+}\left(  q\right)  =\,\sqrt{\lambda-f\left(
q\right)  }\ \ \ \ \text{and }\ \ \ \hat{\sigma}_{\lambda}^{-}\left(
q\right)  =\,-\sqrt{\lambda-f\left(  q\right)  },\ \ \ \lambda\in\mathbb{R};
\]
while for the second one, if $f$ is monotonic in some open interval, the
solutions are given by%
\[
\hat{\sigma}_{\lambda}\left(  p\right)  =f^{-1}\left(  \lambda-p^{2}\right)
,\ \ \ \lambda\in\mathbb{R}.
\]

\begin{remark}
Since we are working on a $1$-dimensional manifold, functions $\hat{\sigma
}_{\lambda}$ above always have a primitive $W_{\lambda}$.
\end{remark}

Of course, the Equations $\left(  \ref{lshje}\right)  $ and $\left(
\ref{lnshje}\right)  $ (and consequently their solutions) coincide when
$f\left(  q\right)  =\frac{1}{2}\,q^{2}$: the $1$\emph{-dimensional harmonic
oscillator}.

\subsection{The Poisson, almost-Poisson and Leibniz scenarios}

\label{almost}

Given a manifold $M$, a \emph{Leibniz structure }on $M$ is a $\left(
0,2\right)  $-tensor $\Xi:T^{\ast}M\times T^{\ast}M\rightarrow\mathbb{R}$ (see
Ref. \cite{ortega}). If $\Xi$ is anti-symmetric, then $\Xi$ is an
\emph{almost-Poisson structure} \cite{bates,can,van}, and if, in addition,
\[
\Xi\left(  d\left[  \Xi\left(  df,dg\right)  \right]  ,dh\right)  +\Xi\left(
d\left[  \Xi\left(  dh,df\right)  \right]  ,dg\right)  +\Xi\left(  d\left[
\Xi\left(  dg,dh\right)  \right]  ,df\right)  =0,
\]
for all $f,g,h\in C^{\infty}\left(  M\right)  $, i.e. the\emph{ Jacobi
identity} holds, then $\Xi$ is a Poisson structure. The pair $\left(
M,\Xi\right)  $ is called a Leibniz, an almost-Poisson and a Poisson manifold,
respectively. In any case, we shall define the morphisms of vector bundles
$\Xi^{\sharp},\Xi^{\sharp op}:T^{\ast}M\rightarrow TM$ by the formulae%
\[
\left\langle \alpha,\Xi^{\sharp}\left(  \beta\right)  \right\rangle
=\Xi\left(  \alpha,\beta\right)  \ \ \ \text{and\ \ \ }\left\langle \alpha
,\Xi^{\sharp op}\left(  \beta\right)  \right\rangle =-\Xi\left(  \beta
,\alpha\right)  .
\]
Note that $\Xi^{\sharp\ast}=-\Xi^{\sharp op}$, and consequently
\begin{equation}
\operatorname{Ker}\Xi^{\sharp}=\left(  \operatorname{Im}\Xi^{\sharp
op}\right)  ^{0}. \label{top}%
\end{equation}
Of course, $\Xi^{\sharp}=\Xi^{\sharp op}$ for almost-Poisson and Poisson
structures. Also, for a Poisson manifold, $\operatorname{Im}\Xi^{\sharp}=TM$
(i.e. the Poisson manifold is \emph{transitive}) if and only if there exists a
symplectic form $\omega$ on $M$ such that $\Xi^{\sharp}=\omega^{\sharp}$.

\bigskip

Consider a Leibniz manifold $\left(  M,\Xi\right)  $, a function
$H:M\rightarrow\mathbb{R}$ and the related vector field $X_{H}=\Xi^{\sharp
}\circ dH$. The dynamical system $\left(  M,X_{H}\right)  $ is called a
\emph{Hamiltonian system in a Leibniz manifold} or simply a \emph{Leibniz
system}. Examples of Leibniz systems are given by the \emph{generalized
nonholonomic systems} with linear constraints (see \cite{cg} and references
therein; see also the examples appearing in Ref. \cite{ortega}). Standard
nonholonomic systems correspond to the case in which $\Xi$ is anti-symmetric,
i.e. an almost Poisson structure (see again Refs. \cite{bates,can,van}).

\begin{remark}
\label{nhs}The manifold $M$ corresponding to a generalized and a standard
nonholonomic system is given by a co-distribution\footnote{By a
co-distribution on $Q$ we mean the canonical dual notion of a regular
distribution of constant rank.} $D\subset T^{\ast}Q$ on some manifold $Q$: the
\emph{constraint co-distribution}. So, the restriction $\Pi:=\left.  \pi
_{Q}\right\vert _{D}:D\rightarrow Q$ is a natural fibration for this kind of systems.
\end{remark}

Given a fibration $\Pi:M\rightarrow N$, we shall briefly study the $\Pi$-HJE
for a Leibniz system $\left(  M,X_{H}\right)  $. In particular, we shall show
that, in our setting, the Hamilton-Jacobi Theory developed in \cite{hjp} for
almost-Poisson systems is recovered.

\bigskip

Fix a Leibniz manifold $\left(  M,\Xi\right)  $ and a point $m\in M$. We shall
say that a subspace $V\subset T_{m}M$ is \textbf{weakly isotropic} (resp.
\textbf{co-isotropic}) if
\begin{equation}
V\cap\operatorname{Im}\Xi^{\sharp}\subset\Xi^{\sharp}\left(  V^{0}\right)
\label{tiso}%
\end{equation}
(resp. $\Xi^{\sharp}\left(  V^{0}\right)  \subset V$), and \textbf{weakly
Lagrangian}\emph{ }if $V\cap\operatorname{Im}\Xi^{\sharp}=\Xi^{\sharp}\left(
V^{0}\right)  $. On the other hand, $V$ is \textbf{isotropic} if it is weakly
isotropic and
\begin{equation}
V\subset\operatorname{Im}\Xi^{\sharp op}, \label{stiso}%
\end{equation}
and \textbf{Lagrangian}\emph{ }if it is isotropic and
co-isotropic.\footnote{Our notion of weakly isotropic (resp. weakly
Lagrangian), corresponds to the notion of isotropic (resp. Lagrangian) of Ref.
\cite{hjp}.} In parallel with Definition \ref{paral}, we have the next one.

\begin{definition}
\label{all}Given a fibration $\Pi:M\rightarrow N$ and a section $\sigma
:N\rightarrow M$ of $\Pi$, we shall say that $\Pi$ (resp. $\sigma$) is
\textbf{(weakly) isotropic}, \textbf{co-isotropic} or \textbf{(weakly)
Lagrangian} if so is each fiber of $\operatorname{Ker}\Pi_{\ast}$ (resp.
$\operatorname{Im}\sigma_{\ast}$).
\end{definition}

Note that, when $\Xi$ is an almost-Poisson structure, $\sigma:N\rightarrow M$
is isotropic if and only if
\[
\operatorname{Im}\sigma_{\ast}\subset\Xi^{\sharp}\left[  \left(
\operatorname{Im}\sigma_{\ast}\right)  ^{0}\right]  ,
\]
since $\operatorname{Im}\Xi^{\sharp}=\operatorname{Im}\Xi^{\sharp op}$. Of
course, if $\left(  M,\Xi\right)  $ is a symplectic manifold, since
$\operatorname{Im}\Xi^{\sharp}=TM$, the above notions of weakly isotropic
(resp. weakly Lagrangian) and isotropic (resp. Lagrangian) coincide with that
given at the beginning of Section \ref{syscen}.

\begin{proposition}
\label{siso}Fix a Leibniz manifold $\left(  M,\Xi\right)  $, a function
$H:M\rightarrow\mathbb{R}$, a fibration $\Pi:M\rightarrow N$ and a section
$\sigma:N\rightarrow M$ of $\Pi$.

\begin{enumerate}
\item If $\sigma$ is a solution of the $\Pi$-HJE for $\left(  M,X_{H}\right)
$ and $\sigma$ is weakly isotropic, then%
\begin{equation}
\operatorname{Im}\left(  \left.  dH\right\vert _{\operatorname{Im}\sigma
}\right)  \subset\left(  \operatorname{Im}\sigma_{\ast}\cap\operatorname{Im}%
\Xi^{\sharp op}\right)  ^{0}. \label{aphje}%
\end{equation}

If in addition $\sigma$ is isotropic, then%
\begin{equation}
d\left(  H\circ\sigma\right)  =0. \label{saphje}%
\end{equation}

\item On the other hand, if $\sigma$ satisfies Eq. $\left(  \ref{aphje}%
\right)  $ [or Eq. $\left(  \ref{saphje}\right)  $] and $\sigma$ is
co-isotropic, then $\sigma$ is a solution of the $\Pi$-HJE for $\left(
M,X_{H}\right)  $.
\end{enumerate}
\end{proposition}

\emph{Proof.} $\left(  1\right)  $ Suppose that a solution $\sigma$ of the
$\Pi$-HJE for $\left(  M,X_{H}\right)  $ is given. According to Proposition
\ref{Xims}, this is equivalent to $\operatorname{Im}\left(  \left.
X_{H}\right\vert _{\operatorname{Im}\sigma}\right)  \subset T\operatorname{Im}%
\sigma=\operatorname{Im}\sigma_{\ast}$. Since in addition
\[
\operatorname{Im}\left(  \left.  X_{H}\right\vert _{\operatorname{Im}\sigma
}\right)  =\operatorname{Im}\left(  \Xi^{\sharp}\circ\left.  dH\right\vert
_{\operatorname{Im}\sigma}\right)  \subset\operatorname{Im}\Xi^{\sharp},
\]
then $\sigma$ is a solution of the $\Pi$-HJE if and only if%
\begin{equation}
\operatorname{Im}\left(  \left.  X_{H}\right\vert _{\operatorname{Im}\sigma
}\right)  \subset\operatorname{Im}\sigma_{\ast}\cap\operatorname{Im}%
\Xi^{\sharp}. \label{imxh}%
\end{equation}
On the other hand, since $\operatorname{Im}\sigma_{\ast}\cap\operatorname{Im}%
\Xi^{\sharp}\subset\Xi^{\sharp}\left[  \left(  \operatorname{Im}\sigma_{\ast
}\right)  ^{0}\right]  $ [see $\left(  \ref{tiso}\right)  $], Eq. $\left(
\ref{imxh}\right)  $ implies that
\[
\operatorname{Im}\left(  \left.  X_{H}\right\vert _{\operatorname{Im}\sigma
}\right)  \subset\Xi^{\sharp}\left[  \left(  \operatorname{Im}\sigma_{\ast
}\right)  ^{0}\right]  ,
\]
i.e. $\Xi^{\sharp}\left(  dH\left(  m\right)  \right)  =X_{H}\left(  m\right)
\in\Xi^{\sharp}\left[  \left(  \operatorname{Im}\sigma_{\ast,m}\right)
^{0}\right]  $ for all $m\in\operatorname{Im}\sigma$. Then,
\[
dH\left(  m\right)  \in\left(  \operatorname{Im}\sigma_{\ast,m}\right)
^{0}+\operatorname{Ker}\Xi_{m}^{\sharp}=\left(  \operatorname{Im}\sigma
_{\ast,m}\cap\left(  \operatorname{Ker}\Xi_{m}^{\sharp}\right)  ^{0}\right)
^{0},
\]
or using Eq. $\left(  \ref{top}\right)  $,
\[
dH\left(  m\right)  \in\left(  \operatorname{Im}\sigma_{\ast,m}\cap
\operatorname{Im}\Xi_{m}^{\sharp op}\right)  ^{0},\ \ \ \forall m\in
\operatorname{Im}\sigma.
\]
This is precisely Eq. $\left(  \ref{aphje}\right)  $. If in addition
$\operatorname{Im}\sigma_{\ast}\cap\operatorname{Im}\Xi^{\sharp op}%
=\operatorname{Im}\sigma_{\ast}$ [see Eq. $\left(  \ref{stiso}\right)  $], the
last equation implies that $\sigma_{m}^{\ast}\left(  dH\left(  m\right)
\right)  =0$ for all $m\in\operatorname{Im}\sigma$, i.e. $d\left(
H\circ\sigma\right)  =0$, as we wanted to show.

\bigskip

$\left(  2\right)  $ Applying $\Xi^{\sharp}$ to $\left(  \ref{aphje}\right)  $
we obtain that%
\begin{align*}
\operatorname{Im}\left(  \left.  X_{H}\right\vert _{\operatorname{Im}\sigma
}\right)   &  \subset\Xi^{\sharp}\left[  \left(  \operatorname{Im}\sigma
_{\ast}\cap\operatorname{Im}\Xi^{\sharp op}\right)  ^{0}\right]  =\Xi^{\sharp
}\left[  \left(  \operatorname{Im}\sigma_{\ast}\right)  ^{0}\right]
+\Xi^{\sharp}\left[  \left(  \operatorname{Im}\Xi^{\sharp op}\right)
^{0}\right]  =\Xi^{\sharp}\left[  \left(  \operatorname{Im}\sigma_{\ast
}\right)  ^{0}\right]  ,
\end{align*}
where we have used Eq. $\left(  \ref{top}\right)  $ in the last step. Thus,
\[
\operatorname{Im}\left(  \left.  X_{H}\right\vert _{\operatorname{Im}\sigma
}\right)  \subset\Xi^{\sharp}\left[  \left(  \operatorname{Im}\sigma_{\ast
}\right)  ^{0}\right]  .
\]
The same we would obtain from Eq. $\left(  \ref{saphje}\right)  $. Since in
addition $\Xi^{\sharp}\left[  \left(  \operatorname{Im}\sigma_{\ast}\right)
^{0}\right]  \subset\operatorname{Im}\sigma_{\ast}$ (i.e. $\sigma$ is
co-isotropic), the Equation $\left(  \ref{imxh}\right)  $ follows immediately
from the above one, i.e. $\sigma$ is a solution of the $\Pi$-HJE for $\left(
M,X_{H}\right)  $.\ \ \ $\square$

\bigskip

Eq. $\left(  \ref{aphje}\right)  $ is exactly the Hamilton-Jacobi equation
that appears in Ref. \cite{hjp} for almost-Poisson systems with a fibration
$\Pi$ (see Theorem 2.3 of Ref. \cite{hjp}). Eq. $\left(  \ref{saphje}\right)
$, on the other hand, is just the equation found in the same reference for the
systems defined on a symplectic manifold (see Theorem 2.4 of \cite{hjp}):
the classical HJE. The next result, which is an immediate consequence of the
above proposition, shows that the mentioned equations are equivalent to our
$\Pi$-HJE for weakly Lagrangian and Lagrangian sections, respectively.

\begin{corollary}
\label{slag}Under the conditions of the above proposition,

\begin{enumerate}
\item if $\sigma$ is weakly Lagrangian, then $\sigma$ is a solution of the
$\Pi$-HJE for $\left(  M,X_{H}\right)  $ if and only if Eq. $\left(
\ref{aphje}\right)  $ holds.

\item if $\sigma$ is Lagrangian, then $\sigma$ is a solution of the $\Pi$-HJE
for $\left(  M,X_{H}\right)  $ if and only if Eq. $\left(  \ref{saphje}%
\right)  $ holds.
\end{enumerate}
\end{corollary}

\bigskip

Consider a nonholonomic system with constraint co-distribution $D\subset
T^{\ast}Q$ (see Remark \ref{nhs}). The Hamilton-Jacobi equation obtained in
Ref. \cite{hjp} for this system (seen as an almost-Poisson system), and for
the natural fibration $\Pi:=\left.  \pi_{Q}\right\vert _{D}:D\rightarrow Q$,
coincides with the so-called \emph{nonholonomic Hamilton-Jacobi equation} of
Refs. \cite{pepin-noholo, lmm, bmmp}. Then, based on the last corollary, we
can conclude that the nonholonomic Hamilton-Jacobi equation is precisely our
$\left.  \pi_{Q}\right\vert _{D}$-HJE (restricted to weakly Lagrangian and
Lagrangian sections).

\section{Complete solutions}

\label{S3}

Fix again a dynamical system $\left(  M,X\right)  $ and a fibration
$\Pi:M\rightarrow N$. Following previous works on the subject (see for
instance Ref. \cite{pepin-holo}), in this section we shall introduce and study
the notion of a complete solution of the $\Pi$-HJE for $\left(  M,X\right)  $.
We shall see that this tool can be used to construct local coordinates in
which the equation of motions of the system can be substantially simplified.
Moreover, we shall see, in the context of Hamiltonian systems on symplectic
manifolds, that for certain complete solutions and certain fibrations $\Pi$,
the trajectories of $X$ can be obtained up to quadratures (as in the classical
situation). This result will be extended to Poisson manifolds (and to more
general fibrations $\Pi$) in the last section of the paper.

\subsection{Definition and basic properties}

\label{totsol}

Let $\Lambda$ be an $l$-manifold with $l:=\dim M-\dim N$.

\begin{definition}
We shall say that $\Sigma:N\times\Lambda\rightarrow M$ is a \textbf{complete
solution of the }$\Pi$-\textbf{HJE} for $\left(  M,X\right)  $ if

\begin{description}
\item[T1] $\Sigma$ is surjective,

\item[T2] $\Sigma$ is a local diffeomorphism (\emph{ipso facto} an open map),

\item[T3] for each $\lambda\in\Lambda$, the map%
\begin{equation}
\sigma_{\lambda}:=\Sigma\left(  \cdot,\lambda\right)  :n\in N\longmapsto
\Sigma\left(  n,\lambda\right)  \in M \label{ps}%
\end{equation}
is a solution of the $\Pi$-HJE for $\left(  M,X\right)  $, i.e. is a section
of $\Pi$ solving $\left(  \ref{hje}\right)  $. Each map $\sigma_{\lambda}$
will be called \textbf{partial solution}.
\end{description}

Given an open subset $U\subset M$, we shall say that $\Sigma:\Pi\left(
U\right)  \times\Lambda\rightarrow U$ is a \textbf{local complete solution on}
$U$\textbf{ of the }$\Pi$\textbf{-HJE} for $\left(  M,X\right)  $ if it is a
complete solution of the $\left.  \Pi\right\vert _{U}$-HJE for $\left(
U,\left.  X\right\vert _{U}\right)  $ (recall Definition \ref{maindef}).
\end{definition}

It is easy to show that, if the surjectivity condition (i.e. condition $T1$)
on $\Sigma:N\times\Lambda\rightarrow M$ is not fulfilled, then $\Sigma$ still
defines (by co-restriction) a global complete solution, provided we change the
phase space $M$ of the system by the open submanifold $\Sigma\left(
N\times\Lambda\right)  \subset M$.

\begin{example}
\label{sex}Going back to Example \ref{fex}, suppose that a family of solutions
$\hat{\sigma}_{\lambda}$ of the Eq. $\left(  \ref{achje}\right)  $, with
$\lambda\in\mathbb{R}^{l}$, is given. It is easy to show that
\[
\Sigma:\left(  n,\lambda\right)  \in\mathbb{R}^{k}\times\mathbb{R}%
^{l}\longmapsto\left(  n,\hat{\sigma}_{\lambda}\left(  n\right)  \right)
\in\mathbb{R}^{k}\times\mathbb{R}^{l}%
\]
is a local diffeomorphism if and only if the $l\times l$ matrix with
coefficients $\left.  \partial\left(  \hat{\sigma}_{\lambda}\right)
_{i}\right/  \partial\lambda_{j}$ is non-degenerate for all $\left(
n,\lambda\right)  $. In such a case, according to above discussion, $\Sigma$
defines a global complete solution with $\Lambda=\mathbb{R}^{l}$ (taking the
phase space equal to the image of $\Sigma$).
\end{example}

Existence conditions of complete solutions will be studied in Section
\ref{exist}. Now, let us study some consequences of having a complete solution
$\Sigma$. First, we shall focus on the properties $T1$ and $T3$ of $\Sigma$.
Let $\sigma_{\lambda}$ be a partial solution, and define%
\begin{equation}
M_{\lambda}:=\operatorname{Im}\sigma_{\lambda}\ \ \ \text{and\ \ \ }%
X^{\sigma_{\lambda}}:=\Pi_{\ast}\circ X\circ\sigma_{\lambda}\in\mathfrak{X}%
\left(  N\right)  . \label{Ml}%
\end{equation}
According to the discussions we made in Section \ref{hjemain}, each
$M_{\lambda}$ is a closed regular submanifold diffeomorphic to $N$, the vector
field $X\in\mathfrak{X}\left(  M\right)  $ restricts to $M_{\lambda}$ (see
Proposition \ref{Xims}), and $X$ and $X^{\sigma_{\lambda}}$ are $\sigma
_{\lambda}$-related. On the other hand, since $\Sigma$ is surjective,%
\[
M=\Sigma\left(  N\times\Lambda\right)  =%
{\displaystyle\bigcup\nolimits_{\lambda\in\Lambda}}
\sigma_{\lambda}\left(  N\right)  =%
{\displaystyle\bigcup\nolimits_{\lambda\in\Lambda}}
M_{\lambda}.
\]
From this point, using Theorem \ref{parcial}, it can be shown that all of
the trajectories of $X$ can be constructed from those of the vector fields
$X^{\sigma_{\lambda}}$'s. We shall see that below, in terms of the following
characterization of the complete solutions. Denote by $p_{N}:N\times
\Lambda\rightarrow N$ and $p_{\Lambda}:N\times\Lambda\rightarrow\Lambda$ the
canonical projections.

\begin{proposition}
\label{pisrel}A surjective local diffeomorphism $\Sigma:N\times\Lambda
\rightarrow M$ is a complete solution of the $\Pi$-HJE for $\left(
M,X\right)  $ if and only if%
\begin{equation}
\Pi\circ\Sigma=p_{N}\ \ \ \text{and}\ \ \ \Sigma_{\ast}\circ X^{\Sigma}%
=X\circ\Sigma, \label{Srel}%
\end{equation}
being $X^{\Sigma}\in\mathfrak{X}\left(  N\times\Lambda\right)  $ the unique
vector field on $N\times\Lambda$ satisfying
\begin{equation}
\left(  p_{N}\right)  _{\ast}\circ X^{\Sigma}=\Pi_{\ast}\circ X\circ
\Sigma\ \ \ \text{and}\ \ \ \left(  p_{\Lambda}\right)  _{\ast}\circ
X^{\Sigma}=0. \label{XS}%
\end{equation}
In particular, the fields $X$ and $X^{\Sigma}$ are $\Sigma$-related.
\end{proposition}

\emph{Proof.} Let $\Sigma:N\times\Lambda\rightarrow M$ be a complete solution
of the $\Pi$-HJE with partial solutions $\sigma_{\lambda}:N\rightarrow M$.
Since each $\sigma_{\lambda}$ is a section of $\Pi$, then%
\[
\Pi\circ\Sigma\left(  n,\lambda\right)  =\Pi\circ\sigma_{\lambda}\left(
n\right)  =n=p_{N}\left(  n,\lambda\right)  ,
\]
from which the first part of $\left(  \ref{Srel}\right)  $ follows. Let us
prove the second one. First, note that the vector field satisfying $\left(
\ref{XS}\right)  $ is given on $\left(  n,\lambda\right)  $ by%
\[
X^{\Sigma}\left(  n,\lambda\right)  =\left(  \Pi_{\ast}\circ X\circ
\Sigma\left(  n,\lambda\right)  ,0\right)  .
\]
Then, using Eq. $\left(  \ref{hje}\right)  $ for each $\sigma_{\lambda}$, we
have that%
\begin{align*}
\Sigma_{\ast}\circ X^{\Sigma}\left(  n,\lambda\right)   &  =\Sigma
_{\ast,\left(  n,\lambda\right)  }\left(  \Pi_{\ast}\circ X\circ\Sigma\left(
n,\lambda\right)  ,0\right)  =\left(  \sigma_{\lambda}\right)  _{\ast
,n}\left(  \Pi_{\ast}\circ X\circ\sigma_{\lambda}\left(  n\right)  \right) \\
&  =\left(  \sigma_{\lambda}\right)  _{\ast}\circ\Pi_{\ast}\circ X\circ
\sigma_{\lambda}\left(  n\right)  =X\circ\sigma_{\lambda}\left(  n\right)
=X\circ\Sigma\left(  n,\lambda\right)  ,
\end{align*}
and the second part of $\left(  \ref{Srel}\right)  $ is obtained.\ The
converse is left to the reader.\ \ \ \ $\square$

\bigskip

Using $\left(  \ref{Ml}\right)  $ and $\left(  \ref{XS}\right)  $, it is easy
to see that
\begin{equation}
X^{\Sigma}\left(  n,\lambda\right)  =\left(  \Pi_{\ast}\circ X\circ
\Sigma\left(  n,\lambda\right)  ,0\right)  =\left(  X^{\sigma_{\lambda}%
}\left(  n\right)  ,0\right)  ,\ \ \ \forall\left(  n,\lambda\right)  \in
N\times\Lambda. \label{XSXs}%
\end{equation}
Then, any integral curve of $X^{\Sigma}$ is of the form $t\mapsto\left(
\gamma\left(  t\right)  ,\lambda\right)  $, for some $\lambda\in\Lambda$ and
some integral curve $\gamma$ of $X^{\sigma_{\lambda}}$. On the other hand,
according to the last proposition, $X$ and $X^{\Sigma}$ are $\Sigma$-related.
Consequently, given an integral curve $\gamma$ of $X^{\sigma_{\lambda}}$,%
\begin{equation}
\Gamma\left(  t\right)  =\Sigma\left(  \gamma\left(  t\right)  ,\lambda
\right)  =\sigma_{\lambda}\left(  \gamma\left(  t\right)  \right)  \label{Gg}%
\end{equation}
is an integral curve of $X$. Moreover, using the surjectivity of $\Sigma$ and
Theorem \ref{parcial}, every trajectory of $\left(  M,X\right)  $ can be
obtained in that way. More precisely,

\begin{theorem}
\label{complete}Let $\Sigma:N\times\Lambda\rightarrow M$ be a complete
solution of the $\Pi$-HJE for $\left(  M,X\right)  $. For every integral curve
$\Gamma:I\rightarrow M$ of $X$ there exists $\lambda\in\Lambda$ such that we
can write $\Gamma=\sigma_{\lambda}\circ\gamma$ for a unique integral curve
$\gamma$ of $X^{\sigma_{\lambda}}$.
\end{theorem}

Let us now exploit the condition $T2$, i.e. the fact that $\Sigma$ is a local
diffeomorphism. It is clear that, for every $m\in M$ there exist $\left(
n,\lambda\right)  \in N\times\Lambda$ and open charts $\left(  U,\psi\right)
$, $\left(  V_{N},\psi_{N}\right)  $ and $\left(  V_{\Lambda},\psi_{\Lambda
}\right)  $ of $M$, $N$ and $\Lambda$, respectively, such that $\Sigma\left(
n,\lambda\right)  =m$, $n\in V_{N}$, $\lambda\in V_{\Lambda}$, $\Sigma\left(
V_{N}\times V_{\Lambda}\right)  \subset U$ and%
\[
\left.  \Sigma\right\vert _{V_{N}\times V_{\Lambda}}:V_{N}\times V_{\Lambda
}\rightarrow U
\]
is a diffeomorphism.

\begin{proposition}
\label{loccoor}If $\Sigma:N\times\Lambda\rightarrow M$ is a complete solution
of the $\Pi$-HJE for $\left(  M,X\right)  $, then for every $m\in M$ there
exist an open neighborhood $U$ of $m$ and coordinates $\left(  n_{1}%
,...,n_{d-l},\lambda_{1},...,\lambda_{l}\right)  :U\rightarrow\mathbb{R}^{d}$,
with $d:=\dim M$, such that
\[
\left.  X\right\vert _{U}=%
{\displaystyle\sum\nolimits_{i=1}^{d-l}}
f_{i}\,\frac{\partial}{\partial n_{i}}%
\]
for some functions $f_{i}:U\rightarrow\mathbb{R}$.
\end{proposition}

\emph{Proof.} Given $m\in M$, consider the local charts $\left(  V_{N}%
,\psi_{N}\right)  $ and $\left(  V_{\Lambda},\psi_{\Lambda}\right)  $ and the
open subset $U$ described previously. Then, $U$ and the map
\[
\Phi:=\left(  \psi_{N}\times\psi_{\Lambda}\right)  \circ\left(  \left.
\Sigma\right\vert _{V_{N}\times V_{\Lambda}}\right)  ^{-1}:U\rightarrow
\mathbb{R}^{d-l}\times\mathbb{R}^{l}=\mathbb{R}^{d}%
\]
define a local chart of $M$. Fix $u\in U$ and let $\left(  n,\lambda\right)
\in V_{N}\times V_{\Lambda}$ such that $u=\Sigma\left(  n,\lambda\right)  \in
U$. Then, using $\left(  \ref{XSXs}\right)  $ and the fact that $X$ and
$X^{\Sigma}$ are $\Sigma$-related, we have that
\begin{align*}
\Phi_{\ast}\left(  X\left(  u\right)  \right)   &  =\left(  \psi_{N}\times
\psi_{\Lambda}\right)  _{\ast}\circ\left(  \left.  \Sigma\right\vert
_{V_{N}\times V_{\Lambda}}\right)  _{\ast}^{-1}\left(  X\left(  u\right)
\right)  =\left(  \psi_{N}\times\psi_{\Lambda}\right)  _{\ast}\left[
X^{\Sigma}\left(  \left.  \Sigma\right\vert _{V_{N}\times V_{\Lambda}}\right)
_{\ast}^{-1}\left(  u\right)  \right] \\
&  =\left(  \psi_{N}\times\psi_{\Lambda}\right)  _{\ast}\left[  X^{\Sigma
}\left(  n,\lambda\right)  \right]  =\left(  \psi_{N}\times\psi_{\Lambda
}\right)  _{\ast}\left(  X^{\sigma_{\lambda}}\left(  n\right)  ,0\right)
=\left(  \left(  \psi_{N}\right)  _{\ast}\left(  X^{\sigma_{\lambda}}\left(
n\right)  \right)  ,0\right)  ,
\end{align*}
from which the proposition immediately follows.\ \ \ $\square$

\bigskip

Given a local chart as in proposition above, it is clear that the equations of
motion for $X$ along $U$ read%
\begin{equation}
\dot{n}_{i}\left(  t\right)  =f_{i}\left(  n\left(  t\right)  ,\lambda\left(
t\right)  \right)  \ \ \ \text{and}\ \ \ \dot{\lambda}_{j}\left(  t\right)
=0, \label{gs}%
\end{equation}
for $i=1,...,d-l$ and $j=1,...,l$. So, as we have seen in Remark \ref{dkl}, in
order to find the trajectories of $X$, we only need to solve $d-l$ equations,
instead of $d$. In particular, if $d-l=1$, it is easy to see that the
Equations $\left(  \ref{gs}\right)  $ can be solved up to the quadrature
\begin{equation}
t=%
{\displaystyle\int\nolimits_{n_{0}}^{n\left(  t\right)  }}
\frac{ds}{f_{1}\left(  s,\lambda_{0}\right)  },\ \ \ \text{with\ \ \ }%
n_{0}:=n\left(  0\right)  \ \ \ \text{and\ \ \ }\lambda_{0}:=\left(
\lambda_{1}\left(  0\right)  ,...,\lambda_{d-1}\left(  0\right)  \right)  .
\label{ibq1}%
\end{equation}
Here, we have supposed that $f_{1}(s,\lambda_{0})$ is not zero in a certain
neighborhood of $(n_{0},\lambda_{0})$. If $d-l>1$, we shall see in Section
\ref{isoquad} that, in the symplectic and Poisson contexts, certain complete
solutions enable us to find coordinates $\left(  n_{1},...,n_{d-l}\right)  $
for $N$ such that above equations can be solved by quadratures too. A special
situation is described in Section \ref{totcan}, where the mentioned
coordinates are Darboux coordinates.

\subsection{The Leibniz and symplectic scenarios}

In the context of Leibniz systems (see Section \ref{almost}), we have a very
simple characterization of the complete solutions. Let us first note that, if
$\left(  M,\Xi\right)  $ is a Leibniz (resp. almost-Poisson and Poisson)
manifold and $\Phi:P\rightarrow M$ is a local diffeomorphism, then the
assignment%
\begin{equation}
\left(  \alpha,\beta\right)  \mapsto\left\langle \alpha,\Phi_{\ast,p}%
^{-1}\circ\Xi_{\Phi\left(  p\right)  }^{\sharp}\circ\left(  \Phi_{p}^{\ast
}\right)  ^{-1}\left(  \beta\right)  \right\rangle ,\ \ \ \forall p\in
P\ \ \text{and}\ \ \forall\alpha,\beta\in T_{p}^{\ast}P, \label{deflei}%
\end{equation}
defines a Leibniz (resp. almost-Poisson and Poisson) structure on $P$. We
shall denote it $\Phi^{\ast}\Xi$. On the other hand, in the particular case of
a symplectic manifold $\left(  M,\omega\right)  $, its clear that the
pull-back $\Phi^{\ast}\omega\in\Omega^{2}\left(  P\right)  $ is closed and
non-degenerated, so $\Phi^{\ast}\omega$ is a symplectic form on $P$.

\begin{theorem}
\label{char}Consider a Leibniz (resp. almost-Poisson, Poisson and symplectic)
manifold $\left(  M,\Xi\right)  $, a function $H:M\rightarrow\mathbb{R}$ and a
fibration $\Pi:M\rightarrow N$. Then, a surjective local diffeomorphism
$\Sigma:N\times\Lambda\rightarrow M$ is a complete solution of the $\Pi$-HJE
for $\left(  M,X_{H}\right)  $ if and only if $\Pi\circ\Sigma=p_{N}$ and
$X_{H}^{\Sigma}\in\mathfrak{X}\left(  N\times\Lambda\right)  $ [see $\left(
\ref{XS}\right)  $] is the Hamiltonian vector field of $H\circ\Sigma$ w.r.t.
the Leibniz (resp. almost-Poisson, Poisson and symplectic) structure
$\Sigma^{\ast}\Xi$, i.e.%
\begin{equation}
X_{H}^{\Sigma}=\left(  \Sigma^{\ast}\Xi\right)  ^{\sharp}\left(  \Sigma^{\ast
}dH\right)  . \label{leihje}%
\end{equation}

\end{theorem}

\emph{Proof.} We have to show that the second part of $\left(  \ref{Srel}%
\right)  $ for $X=X_{H}$ is equivalent to $\left(  \ref{leihje}\right)  $.
Since $X_{H}=\Xi^{\sharp}\circ dH$, Eq. $\left(  \ref{Srel}\right)  $ says
that
\[
\Sigma_{\ast,\left(  n,\lambda\right)  }\left(  X_{H}^{\Sigma}\left(
n,\lambda\right)  \right)  =\Xi_{\Sigma\left(  n,\lambda\right)  }^{\sharp
}\left(  dH\left(  \Sigma\left(  n,\lambda\right)  \right)  \right)  ,
\]
and since $dH\left(  \Sigma\left(  n,\lambda\right)  \right)  =\left(
\Sigma_{\left(  n,\lambda\right)  }^{\ast}\right)  ^{-1}\left(  d\left(
H\circ\Sigma\right)  \left(  n,\lambda\right)  \right)  $, the Theorem
immediately follows from the definition of $\Sigma^{\ast}\Xi$ [see Eq.
$\left(  \ref{deflei}\right)  $].\ \ \ $\square$

\bigskip

So, every complete solution $\Sigma$ for a Leibniz system $\left(
M,X_{H}\right)  $, with Leibniz structure $\Xi$, is a Leibniz map between
$\left(  M,\Xi\right)  $ and $\left(  N\times\Lambda,\Sigma^{\ast}\Xi\right)
$, and the related dynamical system $\left(  N\times\Lambda,X_{H}^{\Sigma
}\right)  $ is also a Leibniz system, with Hamiltonian function $H\circ\Sigma
$. The same is true if we replace \textquotedblleft Leibniz\textquotedblright%
\ by \textquotedblleft almost-Poisson,\textquotedblright\ \textquotedblleft
Poisson\textquotedblright\ and \textquotedblleft symplectic.\textquotedblright

\begin{definition}
\label{isodef}We shall say that a map $\Sigma:N\times\Lambda\rightarrow M$ is
\textbf{(weakly) isotropic}, \textbf{co-isotropic} or \textbf{(weakly)
Lagrangian} if so is each map $\sigma_{\lambda}:=\Sigma\left(  \cdot
,\lambda\right)  $ (see Definition \ref{all}), for all $\lambda\in\Lambda$.
\end{definition}

The next result is immediate, and will be useful later.

\begin{proposition}
\label{equio}Statements below are equivalent:

\begin{enumerate}
\item $\Sigma$ is (weakly) (co)isotropic;

\item the subspaces $T_{m}M_{\lambda}\subset T_{m}M$ are (weakly)
(co)isotropic, for all $m\in M$ and $\lambda\in\Lambda$;

\item the subspaces $T_{n}N\times0_{\lambda}$ are (weakly) (co)isotropic
w.r.t. the Leibniz structure $\Sigma^{\ast}\Xi$, for all $n\in N$ and
$\lambda\in\Lambda$.
\end{enumerate}
\end{proposition}

And for the fibration $\Pi$, it follows similarly that,

\begin{proposition}
\label{equios}$\Pi$ is (weakly) (co)isotropic if and only if the subspaces
$0_{n}\times T_{\lambda}\Lambda$ are (weakly) (co)isotropic w.r.t. the Leibniz
structure $\Sigma^{\ast}\Xi$, for all $n\in N$ and $\lambda\in\Lambda$.
\end{proposition}

For isotropic and Lagrangian complete solutions $\Sigma$, Proposition
\ref{siso} and Corollary \ref{slag} imply the following results.

\begin{proposition}
\label{HShp}Consider a Leibniz manifold $\left(  M,\Xi\right)  $, a function
$H:M\rightarrow\mathbb{R}$ and a fibration $\Pi:M\rightarrow N$. Assume that
$N$ is connected.

\begin{enumerate}
\item If $\Sigma:N\times\Lambda\rightarrow M$ is an isotropic solution of the
$\Pi$-HJE for $\left(  M,X_{H}\right)  $, then there exists a unique function
$h:\Lambda\rightarrow\mathbb{R}$ such that
\begin{equation}
H\circ\Sigma=h\circ p_{\Lambda}. \label{hshp}%
\end{equation}

\item A Lagrangian surjective local diffeomorphism $\Sigma:N\times
\Lambda\rightarrow M$ is a complete solution of the $\Pi$-HJE for $\left(
M,X_{H}\right)  $ if and only if $\Pi\circ\Sigma=p_{N}$ and Eq. $\left(
\ref{hshp}\right)  $ holds for a unique function $h:\Lambda\rightarrow
\mathbb{R}$.
\end{enumerate}
\end{proposition}

\bigskip

In the case of Hamiltonian systems on symplectic manifolds, the
characterization given by Theorem \ref{char} can be slightly
re-formulated. (Compare to Theorem \ref{phjsim}).

\begin{theorem}
\label{hjsim}Consider a symplectic manifold $\left(  M,\omega\right)  $, a
function $H:M\rightarrow\mathbb{R}$ and a fibration $\Pi:M\rightarrow N$. Then
a surjective local diffeomorphism $\Sigma:N\times\Lambda\rightarrow M$ is a
complete solution of the $\Pi$-HJE for $\left(  M,X_{H}\right)  $ if and only
if $\Pi\circ\Sigma=p_{N}$ and [see Eq. $\left(  \ref{XS}\right)  $]
\begin{equation}
d\left(  H\circ\Sigma\right)  =i_{X_{H}^{\Sigma}}\Sigma^{\ast}\omega.
\label{hjes3}%
\end{equation}
In particular, for all $\lambda\in\Lambda$ [see Eq. $\left(  \ref{Ml}\right)
$],
\begin{equation}
d\left(  H\circ\sigma_{\lambda}\right)  =i_{X_{H}^{\sigma_{\lambda}}}%
\sigma_{\lambda}^{\ast}\omega. \label{hjes4}%
\end{equation}

\end{theorem}

\emph{Proof.} Since $\Sigma^{\ast}\omega$ is a symplectic form on
$N\times\Lambda$, we can write Eq. $\left(  \ref{hjes3}\right)  $ as
\begin{equation}
\Sigma^{\ast}dH=\left(  \Sigma^{\ast}\omega\right)  ^{\flat}\left(
X_{H}^{\Sigma}\right)  , \label{1f}%
\end{equation}
or equivalently $\left(  \Sigma^{\ast}\omega\right)  ^{\sharp}\left(
\Sigma^{\ast}dH\right)  =X_{H}^{\Sigma}$. So, the first affirmation follows
from the previous theorem. To show the second affirmation, note that Eq.
$\left(  \ref{1f}\right)  $ implies that%
\[
\left\langle \Sigma^{\ast}dH,Y\right\rangle =\left(  \Sigma^{\ast}%
\omega\right)  \left(  X_{H}^{\Sigma},Y\right)  ,
\]
for all $Y\in\mathfrak{X}\left(  N\times\Lambda\right)  $. So, if we take
$Y=\left(  y,0\right)  $ with $y\in\mathfrak{X}\left(  N\right)  $, Eq.
$\left(  \ref{hjes4}\right)  $ follows form the above one.\ \ \ \ $\square$

\bigskip

\begin{remark}
According to the Corollaries \ref{coro1} and \ref{coro2}, if $\Pi$ is
isotropic or if $\Sigma$ is co-isotropic, then $\left(  \ref{hjes3}\right)  $
is equivalent to $\left(  \ref{hjes4}\right)  $, for all $\lambda\in\Lambda$.
\end{remark}

Concluding, given a symplectic manifold $\left(  M,\omega\right)  $, every
complete solution $\Sigma$ for a Hamiltonian system $\left(  M,X_{H}\right)  $
is a local symplectomorphism (i.e. a local diffeomorphism and a symplectic
map) between $\left(  M,\omega\right)  $ and $\left(  N\times\Lambda
,\Sigma^{\ast}\omega\right)  $, and the related dynamical system $\left(
N\times\Lambda,X_{H}^{\Sigma}\right)  $ is also a Hamiltonian system (with
Hamiltonian function $H\circ\Sigma$).

\subsection{Lagrangian fibrations and related canonical transformations}

\label{totcan}

Let us continue working with a symplectic manifold $\left(  M,\omega\right)
$, a function $H:M\rightarrow\mathbb{R}$, a fibration $\Pi:M\rightarrow N$ and
a complete solution $\Sigma:N\times\Lambda\rightarrow M$ of the $\Pi$-HJE for
$\left(  M,X_{H}\right)  $. We shall see, under certain hypothesis, that
$\Sigma$ defines Darboux coordinates on which the equations of motion can be
easily solved, as it is well-know in the standard Lagrangian situation (see
Definition \ref{standd}).

\begin{remark}
\label{darsym}Recall that giving Darboux coordinates $\psi:U\subset
M\rightarrow\mathbb{R}^{k}\times\mathbb{R}^{k}$, with $2k=\dim M$, it is the same
as giving a canonical transformation, i.e. a symplectomorphism from $U$ to
$T^{\ast}\mathbb{R}^{k}$ (with its canonical symplectic structure). We just
need to identify $T^{\ast}\mathbb{R}^{k}$ with $\mathbb{R}^{k}\times
\mathbb{R}^{k}$ in the usual way.
\end{remark}

\subsubsection{A related local symplectomorphism}

\label{S3.3.1}

\label{spr}

Assume that $\omega=d\theta$, for some $1$-form $\theta$, and that $N$ is
simply-connected. Let us also assume that $\Sigma$ is isotropic, i.e.
$\sigma_{\lambda}^{\ast}\omega=0$ for all $\lambda$. Then
\[
0=\sigma_{\lambda}^{\ast}\omega=\sigma_{\lambda}^{\ast}d\theta=d\sigma
_{\lambda}^{\ast}\theta,
\]
and accordingly, since $N$ is simply-connected, $\sigma_{\lambda}^{\ast}%
\theta$ is exact for each $\lambda$. This implies that there exists a function
$W_{\lambda}:N\rightarrow\mathbb{R}$ such that $dW_{\lambda}=\sigma_{\lambda
}^{\ast}\theta$.

\begin{remark}
Each function $W_{\lambda}$ can be seen as a generalization of the idea of a
characteristic Hamilton function, presented in the classical setting.
\end{remark}

In turn, the family of functions $W_{\lambda}$'s gives rise to a function
$W:N\times\Lambda\rightarrow\mathbb{R}$ satisfying%
\[
\left\langle \left(  \Sigma^{\ast}\theta-dW\right)  \left(  n,\lambda\right)
,\left(  y,0\right)  \right\rangle =\left\langle \sigma_{\lambda}^{\ast}%
\theta\left(  n\right)  -dW_{\lambda}\left(  n\right)  ,y\right\rangle =0
\]
for all $\left(  n,\lambda\right)  \in N\times\Lambda$ and $y\in T_{n}N$. Now,
let us define $\varphi:N\times\Lambda\rightarrow T^{\ast}\Lambda$ by the
formula%
\begin{equation}
\left\langle \varphi\left(  n,\lambda\right)  ,z\right\rangle =\left\langle
\left(  \Sigma^{\ast}\theta-dW\right)  \left(  n,\lambda\right)  ,\left(
0,z\right)  \right\rangle ,\ \ \ \ \ z\in T_{\lambda}\Lambda. \label{fitot}%
\end{equation}

\begin{proposition}
\label{pr1}Under the previous conditions, the map $\varphi$ is an immersion,
and if $\Sigma$ is Lagrangian, then $\varphi$ is a local diffeomorphism.
\end{proposition}

\emph{Proof.} In order to prove the first affirmation, because of the form of
$\varphi$, it is enough to prove that each map
\[
\varphi_{\lambda}:=\varphi\left(  \cdot,\lambda\right)  :N\rightarrow
T_{\lambda}^{\ast}\Lambda
\]
is an immersion. Given $n\in N$ and $x\in T_{n}N$, and identifying (as usual)
the tangent of the linear space $T_{\lambda}^{\ast}\Lambda$ with itself, it
can be shown that
\[
\left\langle \left(  \varphi_{\lambda}\right)  _{\ast,n}\left(  x\right)
,z\right\rangle =d\left(  \Sigma^{\ast}\theta-dW\right)  \left(  \left(
x,0\right)  ,\left(  0,z\right)  \right)  ,
\]
for all $z\in T_{\lambda}\Lambda$. But $d\left(  \Sigma^{\ast}\theta
-dW\right)  =\Sigma^{\ast}\omega$, so
\[
\left\langle \left(  \varphi_{\lambda}\right)  _{\ast,n}\left(  x\right)
,z\right\rangle =\Sigma^{\ast}\omega\left(  \left(  x,0\right)  ,\left(
0,z\right)  \right)  .
\]
Accordingly, if $\left(  \varphi_{\lambda}\right)  _{\ast,n}\left(  x\right)
=0$, then
\begin{equation}
\Sigma^{\ast}\omega\left(  \left(  x,0\right)  ,\left(  0,z\right)  \right)
=0,\ \ \ \forall\left(  0,z\right)  \in0_{n}\times T_{\lambda}\Lambda.
\label{cis1}%
\end{equation}
On the other hand, since $\Sigma$ is isotropic, each subspace $T_{n}%
N\times0_{\lambda}$ is isotropic w.r.t. $\Sigma^{\ast}\omega$ (see Proposition
\ref{equio}), what implies that%
\begin{equation}
\Sigma^{\ast}\omega\left(  \left(  x,0\right)  ,\left(  y,0\right)  \right)
=0,\ \ \ \forall\left(  y,0\right)  \in T_{n}N\times0_{\lambda}. \label{cis2}%
\end{equation}
Combining $\left(  \ref{cis1}\right)  $ and $\left(  \ref{cis2}\right)  $, it
follows that
\[
\Sigma^{\ast}\omega\left(  \left(  x,0\right)  ,\left(  y,z\right)  \right)
=0,\ \ \ \forall\left(  y,z\right)  \in T_{n}N\times T_{\lambda}\Lambda.
\]
Finally, since $\Sigma^{\ast}\omega$ is non-degenerated, then $x$ must
vanishes. This proves that $\varphi$ is an immersion. It is clear that, if
$\Sigma$ is Lagrangian, then, for dimensional reasons, $\varphi$ is a local
diffeomorphism. \ \ \ $\square$

\bigskip

\begin{proposition}
\label{pr2}Under the previous conditions, if in addition $\Pi$ is isotropic
(and consequently $\Pi$ and $\Sigma$ are Lagrangian), then $\varphi
:N\times\Lambda\rightarrow T^{\ast}\Lambda$ is a local symplectomorphism, i.e.
$\varphi^{\ast}\omega_{\Lambda}=\Sigma^{\ast}\omega$, being $\omega_{\Lambda}$
the canonical symplectic structure on $T^{\ast}\Lambda$.
\end{proposition}

\emph{Proof.} Let us fix a coordinate chart $\psi:V\subset\Lambda
\rightarrow\psi\left(  V\right)  \subset\mathbb{R}^{k}$ for $\Lambda$. Note
that $\left(  \pi_{\Lambda}^{-1}\left(  V\right)  ,\left(  \psi^{\ast}\right)
^{-1}\right)  $ is a Darboux coordinate chart for $\left(  T^{\ast}%
\Lambda,\omega_{\Lambda}\right)  $. Let us write%
\[
\psi\left(  \lambda\right)  =\left(  \lambda^{1},...,\lambda^{k}\right)
\ \ \ \text{and\ \ \ }\left(  \psi^{\ast}\right)  ^{-1}\left(  z\right)
=\left(  \lambda^{1},...,\lambda^{k},\alpha_{1},...,\alpha_{k}\right)
\]
for all $\lambda\in V$ and $z\in T_{\lambda}^{\ast}V$. In this notation, we
have that%
\[
\omega_{\Lambda}\left(  \frac{\partial}{\partial\lambda^{i}},\frac{\partial
}{\partial\lambda^{j}}\right)  =\omega_{\Lambda}\left(  \frac{\partial
}{\partial\alpha_{i}},\frac{\partial}{\partial\alpha_{j}}\right)
=0\ \ \ \text{and\ \ \ }\omega_{\Lambda}\left(  \frac{\partial}{\partial
\lambda^{i}},\frac{\partial}{\partial\alpha_{j}}\right)  =\delta_{i}^{j},
\]
for all $i,j=1,...,k$. Also, since $d\left(  \Sigma^{\ast}\theta-dW\right)
=\Sigma^{\ast}\omega$, it can be shown that
\[
\left(  \psi\circ\varphi\right)  _{\ast}\left(  0,\frac{\partial}%
{\partial\lambda^{i}}\right)  =\left(  \mathbf{e}_{i},\mathbf{a}_{i}\right)
,
\]
where $\mathbf{e}_{i}\in\mathbb{R}^{k}$ is the $i$-th canonical vector and
each $\mathbf{a}_{i}\in\mathbb{R}^{k}$ has $j$-th component%
\[
\left(  \mathbf{a}_{i}\right)  _{j}=\Sigma^{\ast}\omega\left(  \left(
0,\frac{\partial}{\partial\lambda^{i}}\right)  ,\left(  0,\frac{\partial
}{\partial\lambda^{j}}\right)  \right)  ,
\]
and, for all $x\in TN$,
\[
\left(  \psi\circ\varphi\right)  _{\ast}\left(  x,0\right)  =\left(
0,\mathbf{b}\left(  x\right)  \right)  ,
\]
where each $\mathbf{b}\left(  x\right)  \in\mathbb{R}^{k}$ has $j$-th
component%
\[
\mathbf{b}_{j}\left(  x\right)  =\Sigma^{\ast}\omega\left(  \left(
x,0\right)  ,\left(  0,\frac{\partial}{\partial\lambda^{j}}\right)  \right)
.
\]
Using that $\Pi$ is Lagrangian, we have that $\mathbf{a}_{i}=0$ for all $i$
(see Proposition \ref{equios}), and consequently,%
\begin{align*}
\varphi^{\ast}\omega_{\Lambda}\left(  \left(  0,\frac{\partial}{\partial
\lambda^{i}}\right)  ,\left(  0,\frac{\partial}{\partial\lambda^{j}}\right)
\right)   &  =\omega_{\Lambda}\left(  \varphi_{\ast}\left(  0,\frac{\partial
}{\partial\lambda^{i}}\right)  ,\varphi_{\ast}\left(  0,\frac{\partial
}{\partial\lambda^{j}}\right)  \right)  =\left(  \mathbf{e}_{i},0\right)
^{t}\cdot\mathbb{J}\cdot\left(  \mathbf{e}_{j},0\right) \\
& \\
&  =0=\Sigma^{\ast}\omega\left(  \left(  0,\frac{\partial}{\partial\lambda
^{i}}\right)  ,\left(  0,\frac{\partial}{\partial\lambda^{j}}\right)  \right)
,
\end{align*}
where \textquotedblleft$t$\textquotedblright\ indicates transposition and
$\mathbb{J}$ is the real $(2k\times2k)$-matrix
\begin{equation}
\mathbb{J}:=\left[
\begin{array}
[c]{cc}%
0 & -\mathbb{I}\\
\mathbb{I} & 0
\end{array}
\right]  , \label{jota}%
\end{equation}
where $\mathbb{I}$ is the identity $(k\times k)$-matrix. On the other hand,
\begin{align*}
\varphi^{\ast}\omega_{\Lambda}\left(  \left(  x,0\right)  ,\left(  y,0\right)
\right)   &  =\omega_{\Lambda}\left(  \varphi_{\ast}\left(  x,0\right)
,\varphi_{\ast}\left(  y,0\right)  \right)  =\left(  0,\mathbf{b}\left(
x\right)  \right)  ^{t}\cdot\mathbb{J}\cdot\left(  0,\mathbf{b}\left(
y\right)  \right) \\
&  =0=\Sigma^{\ast}\omega\left(  \left(  x,0\right)  ,\left(  y,0\right)
\right)  ,
\end{align*}
for all $x,y\in TN$, where in the last equality we have used that $\Sigma$ is
isotropic. Finally,%
\[
\varphi^{\ast}\omega_{\Lambda}\left(  \left(  x,0\right)  ,\left(
0,\frac{\partial}{\partial\lambda^{j}}\right)  \right)  =\left(
0,\mathbf{b}\left(  x\right)  \right)  ^{t}\cdot\mathbb{J}\cdot\left(
\mathbf{e}_{j},0\right)  =\mathbf{b}_{j}\left(  x\right)  =\Sigma^{\ast}%
\omega\left(  \left(  x,0\right)  ,\left(  0,\frac{\partial}{\partial
\lambda^{j}}\right)  \right)  ,
\]
and the proposition is proved.\ \ \ $\square$

\bigskip

Above results are also valid if we do not assume that $\omega=d\theta$ and $N$
is connected, but at a local level (i.e. around every $m\in M$). This is
justified by the following reasoning:

Let $\left(  M,\omega\right)  $ be a symplectic manifold, $\Pi:M\rightarrow N$
a fibration and $\Sigma:N\times\Lambda\rightarrow M$ a surjective local
diffeomorphism such that $\Pi\circ\Sigma=p_{N}$. Then, for each $m\in M$ there
exist an open neighborhood $U$ of $m$ and a submanifold $\Lambda_{1}%
\subset\Lambda$ such that $\left(  1\right)  $ $i^{\ast}\omega=d\theta$ for
some $\theta\in\Omega^{1}\left(  U\right)  $, with $i:U\rightarrow M$ the
canonical inclusion, $\left(  2\right)  $ $\Pi\left(  U\right)  $ is simply
connected, and $\left(  3\right)  $ $\Sigma$ restricted to $\Pi\left(
U\right)  \times\Lambda_{1}$ is a surjection onto $U$. In particular, if
$\Sigma$ is a (resp. (co)isotropic) complete solution of the $\Pi$-HJE for
some dynamical system $\left(  M,X\right)  $, then $\Sigma$ restricted to
$\Pi\left(  U\right)  \times\Lambda_{1}$ is a (resp. (co)isotropic) complete
solution of the $\left.  \Pi\right\vert _{U}$-HJE for $\left(  U,\left.
X\right\vert _{U}\right)  $.

\subsubsection{Special Darboux coordinates}

Since a complete solution $\Sigma$ is a local symplectomorphism and, according
to the previous results, so is $\varphi$, for every Darboux coordinate chart
$\Psi$ on $T^{\ast}\Lambda$ (that can be composed with $\varphi$) we have a
Darboux coordinate chart $\Phi$ on $M$ giving, roughly speaking, by the
formula%
\begin{equation}
\Phi:=\Psi\circ\varphi\circ\Sigma^{-1}. \label{rsct}%
\end{equation}
What is even more interesting is that, in such charts, the expression of the
equations of motion for $\left(  M,X_{H}\right)  $ are particularly simple.
Note that
\[
\dim M=2\dim N=2\dim\Lambda=2k.
\]

\begin{theorem}
Consider a symplectic manifold $\left(  M,\omega\right)  $, a function
$H:M\rightarrow\mathbb{R}$, a fibration $\Pi:M\rightarrow N$ and a complete
solution $\Sigma:N\times\Lambda\rightarrow M$ of the $\Pi$-HJE for $\left(
M,X_{H}\right)  $. If $\Pi$ and $\Sigma$ are Lagrangian, then for each
$m_{0}\in M$ there exist a Darboux coordinate chart
\[
\Phi:m\in U\mapsto\left(  \lambda,\alpha\right)  \in\mathbb{R}^{k}%
\times\mathbb{R}^{k}%
\]
of $M$ and a coordinate chart $\psi:V_{\Lambda}\rightarrow\mathbb{R}^{k}$ of
$\Lambda$ such that $m_{0}\in U$ and
\[
\left[  \Phi_{\ast}\circ\left(  \left.  X_{H}\right\vert _{U}\right)  \right]
\left(  \lambda,\alpha\right)  =\left(  0,-\nabla\left(  h\circ\psi
^{-1}\right)  \left(  \lambda\right)  \right)  ,\ \ \ \forall\left(
\lambda,\alpha\right)  \in\Phi\left(  U\right)  ,
\]
for some function $h:V_{\Lambda}\rightarrow\mathbb{R}$.
\end{theorem}

\emph{Proof.} Again, for simplicity, we shall assume that $\omega=d\theta$ and
that $N$ is simply-connected (\emph{ipso facto} connected). If this is not
true, we must restrict ourselves to appropriate open subsets of $M$, $N$ and
$\Lambda$, as we explained at the end of Section \ref{S3.3.1}. This enable us
to use Propositions \ref{pr1} and \ref{pr2}.

Given $m_{0}\in M$, consider an open neighborhood $U\subset M$ of $m_{0}$, an
open subset $V_{N}\subset N$, and a local chart $\left(  V_{\Lambda}%
,\psi\right)  $ of $\Lambda$, with $\psi:V_{\Lambda}\rightarrow\psi\left(
V_{\Lambda}\right)  \subset\mathbb{R}^{k}$, such that $\Sigma_{|V_{N}\times
V_{\Lambda}}:V_{N}\times V_{\Lambda}\rightarrow U$ is a diffeomorphism. Now,
shrink $V_{N}$ and $V_{\Lambda}$, if necessary, in such a way that%
\[
\left.  \varphi\right\vert _{V_{N}\times V_{\Lambda}}:V_{N}\times V_{\Lambda
}\rightarrow\varphi\left(  V_{N}\times V_{\Lambda}\right)  \subset T^{\ast
}\Lambda
\]
is a diffeomorphism. Note that $\varphi\left(  V_{N}\times V_{\Lambda}\right)
\subset T^{\ast}V_{\Lambda}$. Consider the map%
\[
\Phi:=\left(  \psi^{\ast}\right)  ^{-1}\circ\left.  \varphi\right\vert
_{V_{N}\times V_{\Lambda}}\circ\left(  \left.  \Sigma\right\vert _{V_{N}\times
V_{\Lambda}}\right)  ^{-1}:U\rightarrow\Phi\left(  U\right)  \subset
\mathbb{R}^{k}\times\mathbb{R}^{k}.
\]
Since $\left(  \psi^{\ast}\right)  ^{-1}$ defines Darboux coordinates on
$T^{\ast}\Lambda$ and the maps $\left.  \varphi\right\vert _{V_{N}\times
V_{\Lambda}}$ and $\left(  \left.  \Sigma\right\vert _{V_{N}\times V_{\Lambda
}}\right)  ^{-1}$ are symplectomorphisms, then $\Phi$ defines Darboux
coordinates on $M$, or equivalently, $\Phi$ defines a symplectomorphism
between $U$ and $T^{\ast}\mathbb{R}^{k}$ (see Remark \ref{darsym}). Writing
$\left(  \psi^{\ast}\right)  ^{-1}\left(  z\right)  =\left(  \lambda
,\alpha\right)  $ for all $\lambda\in V_{\Lambda}$ and $z\in T_{\lambda}%
^{\ast}V_{\Lambda}$, this implies that%
\[
\Phi_{\ast}\circ\left(  \left.  X_{H}\right\vert _{U}\right)  =X_{H\circ
\Phi^{-1}}=\mathbb{J}^{-1}\cdot\left(
\begin{array}
[c]{c}%
\partial\left.  \left(  H\circ\Phi^{-1}\right)  \right/  \partial\lambda\\
\partial\left.  \left(  H\circ\Phi^{-1}\right)  \right/  \partial\alpha
\end{array}
\right)  ,
\]
where $\mathbb{J}$ is given by Eq. $\left(  \ref{jota}\right)  $. On the other
hand, since $N$ is connected and $\Sigma$ is a Lagrangian complete solution,
it follows from Proposition \ref{HShp} that $H\circ\Sigma=h\circ p_{\Lambda}$
for a unique function $h:\Lambda\rightarrow\mathbb{R}$. As a consequence,
\begin{align*}
H\circ\Phi^{-1}  &  =H\circ\left.  \Sigma\right\vert _{V_{N}\times V_{\Lambda
}}\circ\left(  \left.  \varphi\right\vert _{V_{N}\times V_{\Lambda}}\right)
^{-1}\circ\left.  \psi^{\ast}\right\vert _{\Phi\left(  U\right)  }\\
&  =h\circ\left.  p_{\Lambda}\right\vert _{V_{N}\times V_{\Lambda}}%
\circ\left(  \left.  \varphi\right\vert _{V_{N}\times V_{\Lambda}}\right)
^{-1}\circ\left.  \psi^{\ast}\right\vert _{\Phi\left(  U\right)  }=h\circ
\pi_{\Lambda}\circ\left.  \psi^{\ast}\right\vert _{\Phi\left(  U\right)  }.
\end{align*}
So, for every pair $\left(  \lambda,\alpha\right)  \in\Phi\left(  U\right)  $,
we have that%
\[
H\circ\Phi^{-1}\left(  \lambda,\alpha\right)  =h\circ\pi_{\Lambda}\circ
\psi^{\ast}\left(  \lambda,\alpha\right)  =h\circ\psi^{-1}\left(
\lambda\right)  ,
\]
and accordingly
\[
\Phi_{\ast}\circ\left(  \left.  X_{H}\right\vert _{U}\right)  =\mathbb{J}%
^{-1}\cdot\left(
\begin{array}
[c]{c}%
\partial\left.  \left(  h\circ\psi^{-1}\right)  \right/  \partial\lambda\\
0
\end{array}
\right)  =\left(  0,-\frac{\partial\left(  h\circ\psi^{-1}\right)  }%
{\partial\lambda}\right)  .\ \ \ \ \ \square
\]

\bigskip

The equations of motion on the coordinates $\left(  \lambda,\alpha\right)
=\left(  \lambda^{1},...,\lambda^{k},\alpha_{1},...,\alpha_{k}\right)  $,
defined in the last theorem, are%
\begin{equation}
\dot{\alpha}_{i}\left(  t\right)  =-\frac{\partial\left(  h\circ\psi
^{-1}\right)  }{\partial\lambda^{i}}\left(  \lambda\left(  t\right)  \right)
\ \ \ \text{and\ \ \ \ }\dot{\lambda}^{i}\left(  t\right)
=0,\ \ \ \ \ i=1,...,k, \label{sem}%
\end{equation}
which can be easily integrated [compare with the general situation, given by
Eq. $\left(  \ref{gs}\right)  $]. Note that, to arrive at the last equations,
we need to construct (at least locally) the function $W$ and the $1$-form
$\theta$, which give rise to the map $\varphi$ [see Eq. $\left(
\ref{fitot}\right)  $]. According to the Poincar\'{e} Lemma, they can be found
by quadratures. So, we can conclude that, given a Lagrangian fibration $\Pi$,
if a Lagrangian solution of the $\Pi$-HJE for $\left(  M,X_{H}\right)  $ is
known, then $\left(  M,X_{H}\right)  $ can be integrated by quadratures.
Weaker conditions for integrability by quadratures will be given in Section
\ref{ghjm}.

\subsection{Some examples: type $1$ to $4$ generating functions}

\label{canchan}

Let us illustrate previous results with simple examples. Let us go back to
Section \ref{stanon}, and concentrate first on the $r=1$ case, i.e. in the
standard situation. Suppose that a family of (local) solutions $\sigma
_{\lambda}:q\mapsto\sigma_{\lambda}\left(  q\right)  =\left(  q,\hat{\sigma
}_{\lambda}\left(  q\right)  \right)  $ is given for the $\mathfrak{p}_{1}%
$-HJE, with $\lambda\in\mathbb{R}^{k}$, such that the square matrix with
coefficients $\left.  \partial\left(  \hat{\sigma}_{\lambda}\right)
_{i}\right/  \partial\lambda_{j}$ is non-degenerate for all $\left(
q,\lambda\right)  $. Then, $\Sigma\left(  q,\lambda\right)  :=\left(
q,\hat{\sigma}_{\lambda}\left(  q\right)  \right)  $ is a local
diffeomorphism, and consequently $\Sigma$ defines a (local) complete solution
of the $\mathfrak{p}_{1}$-HJE.

\begin{example}
Take $k=1$ and
\[
H\left(  q,p\right)  =\frac{1}{2}\,\left(  p^{2}+q\right)  .
\]
According to the results of Section \ref{stanon} for $f\left(  q\right)  =q$,
the solutions of the $\mathfrak{p}_{1}$-HJE are
\[
\hat{\sigma}_{\lambda}^{+}\left(  q\right)  =\,\sqrt{\lambda-q}%
\ \ \ \ \text{and }\ \ \ \hat{\sigma}_{\lambda}^{-}\left(  q\right)
=\,-\sqrt{\lambda-q},\ \ \ \lambda\in\mathbb{R}.
\]
As a consequence, for every $a\in\mathbb{R}$, we have the local complete
solutions
\[
\Sigma^{\pm}:\left(  -\infty,a\right)  \times\left(  a,\infty\right)
\rightarrow\mathbb{R}^{2}%
\]
given by%
\[
\Sigma^{+}\left(  q,\lambda\right)  =\left(  q,\sqrt{\lambda-q}\right)
\ \ \ \text{and\ \ \ }\Sigma^{-}\left(  q,\lambda\right)  =\left(
q,-\sqrt{\lambda-q}\right)  .
\]
They are local diffeomorphisms since%
\[
\frac{\partial}{\partial\lambda}\left(  \pm\sqrt{\lambda-q}\right)  =\mp
\frac{1}{2\sqrt{\lambda-q}}\neq0.
\]

\end{example}

In addition, suppose that each solution $\sigma_{\lambda}$ is isotropic (and
consequently Lagrangian), i.e. $d\sigma_{\lambda}=0$ for each $\lambda$ [see
Eq. $\left(  \ref{ests}\right)  $]. Then, there exists a family of functions
$W_{\lambda}$ such that%
\begin{equation}
\left(  \hat{\sigma}_{\lambda}\right)  _{i}\left(  q\right)  =\frac{\partial
W_{\lambda}\left(  q\right)  }{\partial q^{i}},\ \ \ \frac{\partial H\left(
q,\left.  \partial W_{\lambda}\left(  q\right)  \right/  \partial q\right)
}{\partial q^{i}}=0,\ \ \ \ i=1,...k, \label{can12}%
\end{equation}
[see Eq. $\left(  \ref{lshje}\right)  $] and the matrix with coefficients
$\displaystyle\frac{\partial^{2}W_{\lambda}\left(  q\right)  }{\partial\lambda_{j}\partial
q^{i}}$ is non-degenerate. In this situation, as it is well-known, the
formulae%
\[
p_{i}=\frac{\partial W_{\lambda}\left(  q\right)  }{\partial q^{i}%
}\ \ \ \text{and\ \ \ }\alpha^{i}=-\frac{\partial W_{\lambda}\left(  q\right)
}{\partial\lambda_{i}}%
\]
define two canonical transformations (i.e. a change of coordinates between
Darboux coordinates)%
\begin{equation}
\left(  q,p\right)  \longmapsto\left(  Q,P\right)  :=\left(  \lambda
,\alpha\right)  \label{ct1}%
\end{equation}
and
\begin{equation}
\left(  q,p\right)  \longmapsto\left(  Q,P\right)  :=\left(  -\alpha
,\lambda\right)  , \label{ct2}%
\end{equation}
and the function $W\left(  q,\lambda\right)  =W_{\lambda}\left(  q\right)  $
plays the role of a \emph{generating function} of type $1$ and type $2$ (see
Ref. \cite{gold}), respectively. Moreover, Eq. $\left(  \ref{can12}\right)  $
implies that the related Hamiltonians $K_{1}\left(  Q,P\right)  $ and
$K_{2}\left(  Q,P\right)  $ only depend on $\lambda$, i.e. $K_{1}=K_{1}\left(
Q\right)  $ and $K_{2}=K_{2}\left(  P\right)  $. This makes the integration of
the equation of motion, in the new Darboux coordinates $\left(  Q,P\right)  $,
a trivial task. For instance, for the type $2$ case, the equations of motion
translate to [compare to Eq. $\left(  \ref{sem}\right)  $]%
\[
\dot{Q}^{i}\left(  t\right)  =\frac{\partial K_{2}}{\partial P_{i}}\left(
P\left(  t\right)  \right)  ,\ \ \ \dot{P}_{i}\left(  t\right)  =0,
\]
so, the solutions are%
\[
\left(  Q^{i}\left(  t\right)  ,P_{i}\left(  t\right)  \right)  =\left(
\frac{\partial K_{2}}{\partial P_{i}}\left(  P\left(  0\right)  \right)
\,t+Q^{i}\left(  0\right)  ,P_{i}\left(  0\right)  \right)  .
\]

\bigskip

Let us analyze the previous construction in terms of the map $\varphi$ and the
canonical transformation $\Phi$ [see Eqs. $\left(  \ref{fitot}\right)  $ and
$\left(  \ref{rsct}\right)  $] defined in the previous section. Since%
\[
\theta\left(  \Sigma\left(  q,\lambda\right)  \right)  =\theta\left(
q,\hat{\sigma}_{\lambda}\left(  q\right)  \right)  =\left(  \hat{\sigma
}_{\lambda}\right)  _{i}\left(  q\right)  \,dq^{i},
\]
then (omitting the dependence on $q$ and $\lambda$)
\[
\left\langle \sigma_{\lambda}^{\ast}\theta,\frac{\partial}{\partial q^{i}%
}\right\rangle =\left\langle \theta,\left(  \mathbf{e}_{i},\frac{\partial
\hat{\sigma}_{\lambda}}{\partial q^{i}}\right)  \right\rangle =\left(
\hat{\sigma}_{\lambda}\right)  _{i},
\]
where $\mathbf{e}_{i}$ is the $i$-th vector of the canonical basis of
$\mathbb{R}^{k}$, and%
\[
\left\langle \sigma_{\lambda}^{\ast}\theta,\frac{\partial}{\partial\lambda
_{i}}\right\rangle =\left\langle \theta,\left(  0,\frac{\partial\hat{\sigma
}_{\lambda}}{\partial\lambda_{i}}\right)  \right\rangle =0.
\]
Accordingly,
\[
\varphi\left(  q,\lambda\right)  =\left\langle \sigma_{\lambda}^{\ast}%
\theta-dW_{\lambda},\frac{\partial}{\partial\lambda_{i}}\right\rangle
\,d\lambda_{i}=-\frac{\partial W_{\lambda}}{\partial\lambda_{i}}\,d\lambda
_{i}=\left(  \lambda,-\frac{\partial W_{\lambda}}{\partial\lambda}\right)  ,
\]
with $W_{\lambda}$ such that $\left(  \hat{\sigma}_{\lambda}\right)
_{i}=\left.  \partial W_{\lambda}\right/  \partial q^{i}$. Thus, for $\left(
q,p\right)  =\Sigma\left(  q,\lambda\right)  $, the map $\Phi:=\Psi
\circ\varphi\circ\Sigma^{-1}$ [recall Eq. $\left(  \ref{rsct}\right)  $] is
given by%
\[
\Phi\left(  q,p\right)  =\Psi\circ\varphi\left(  q,\lambda\right)
=\Psi\left(  \lambda,-\frac{\partial W_{\lambda}}{\partial\lambda}\right)
=:\left(  Q,P\right)  ,
\]
So, if we take $\Psi=\Psi_{1}$ or $\Psi=\Psi_{2}$, with
\begin{equation}
\Psi_{1}\left(  \lambda,\alpha\right)  :=\left(  \lambda,\alpha\right)
\ \ \ \ \text{and \ \ \ }\Psi_{2}\left(  \lambda,\alpha\right)  :=\left(
-\alpha,\lambda\right)  , \label{fi12}%
\end{equation}
we have the canonical transformations $\left(  \ref{ct1}\right)  $ or $\left(
\ref{ct2}\right)  $, respectively.

\bigskip

Consider now the $\mathfrak{p}_{2}$-HJE. If we have a family of (local)
solutions $\sigma_{\lambda}:p\mapsto\sigma_{\lambda}\left(  p\right)  =\left(
\hat{\sigma}_{\lambda}\left(  p\right)  ,p\right)  $, with $\lambda
\in\mathbb{R}^{k}$, such that the matrix with coefficients $\left.
\partial\left(  \hat{\sigma}_{\lambda}\right)  ^{i}\right/  \partial
\lambda_{j}$ is non-degenerate for all $\left(  p,\lambda\right)  $, then
$\Sigma\left(  p,\lambda\right)  :=\left(  \hat{\sigma}_{\lambda}\left(
p\right)  ,p\right)  $ is a local complete solution of the $\mathfrak{p}_{2}%
$-HJE. Again, suppose in addition each $\sigma_{\lambda}$ is Lagrangian.
According to $\left(  \ref{isocon}\right)  $, this also means that there exist
functions $\widetilde{W}_{\lambda}$ such that%
\begin{equation}
\left(  \hat{\sigma}_{\lambda}\right)  ^{i}\left(  p\right)  =\frac
{\partial\widetilde{W}_{\lambda}\left(  p\right)  }{\partial p_{i}%
},\ \ \ \ \ \ \ \ \frac{\partial H\left(  \left.  \partial\widetilde
{W}_{\lambda}\left(  p\right)  \right/  \partial p,p\right)  }{\partial p_{i}%
}=0,\ \ \ \ i=1,...k, \label{can34}%
\end{equation}
and the numbers $\displaystyle\frac{\partial^{2}\widetilde{W}_{\lambda}\left(  p\right)
}{\partial\lambda_{j}\partial p_{i}}$ define a non-degenerate matrix. In this
situation, the formulae%
\[
q^{i}=-\frac{\partial\widetilde{W}_{\lambda}\left(  p\right)  }{\partial
p_{i}}\ \ \ \text{and\ \ \ }\alpha_{i}=-\frac{\partial\widetilde{W}_{\lambda
}\left(  p\right)  }{\partial\lambda^{i}}%
\]
define the canonical transformations
\begin{equation}
\left(  q,p\right)  \longmapsto\left(  Q,P\right)  :=\left(  \lambda
,\alpha\right)  \ \ \ \text{and\ \ \ }\left(  q,p\right)  \longmapsto\left(
Q,P\right)  :=\left(  \alpha,-\lambda\right)  , \label{ct34}%
\end{equation}
and the function $\widetilde{W}\left(  p,\lambda\right)  =\widetilde
{W}_{\lambda}\left(  p\right)  $ plays the role of a type $3$ and type $4$
generating function, respectively. Again, because of Eq. $\left(
\ref{can34}\right)  $, the integration of the equations of motion in the new
coordinates is trivial.

\begin{example}
Recalling again the results of Section \ref{stanon}, we have for $f\left(
q\right)  =q$ that the solutions of the $\mathfrak{p}_{2}$-HJE are
$\hat{\sigma}_{\lambda}\left(  p\right)  =\lambda-p^{2}$, $\lambda
\in\mathbb{R}$. Since $\frac{\partial}{\partial\lambda}\left(  \lambda
-p^{2}\right)  =1\neq0$, such solutions give rise to a global complete
solution $\Sigma:\mathbb{R}^{2}\rightarrow\mathbb{R}^{2}$ of the
$\mathfrak{p}_{2}$-HJE given by $\Sigma\left(  p,\lambda\right)  =\left(
\lambda-p^{2},p\right)  $. On the other hand, since
\[
\hat{\sigma}_{\lambda}\left(  p\right)  =\frac{\partial\widetilde{W}_{\lambda
}\left(  p\right)  }{\partial p}\ \ \ \ \text{with\ \ \ }\widetilde
{W}_{\lambda}\left(  p\right)  =\lambda p-p^{3}/3,
\]
the function $\widetilde{W}:\left(  p,\lambda\right)  \in\mathbb{R}^{2}%
\mapsto\lambda p-p^{3}/3\in\mathbb{R}$ defines a type $3$ and a type $4$
generating function.
\end{example}

Let us construct the maps $\varphi$ and $\Phi$ for this case. Since%
\[
\left\langle \sigma_{\lambda}^{\ast}\theta,\frac{\partial}{\partial p_{i}%
}\right\rangle =\left\langle \theta,\left(  \frac{\partial\hat{\sigma
}_{\lambda}}{\partial p_{i}},\mathbf{e}_{i}\right)  \right\rangle
=p_{j}\,\frac{\partial\left(  \hat{\sigma}_{\lambda}\right)  ^{j}}{\partial
p_{i}}%
\]
and%
\[
\left\langle \sigma_{\lambda}^{\ast}\theta,\frac{\partial}{\partial\lambda
_{i}}\right\rangle =\left\langle \theta,\left(  \frac{\partial\hat{\sigma
}_{\lambda}}{\partial\lambda_{i}},0\right)  \right\rangle =p_{j}%
\,\frac{\partial\left(  \hat{\sigma}_{\lambda}\right)  ^{j}}{\partial
\lambda_{i}},
\]
then
\begin{align*}
\varphi\left(  p,\lambda\right)   &  =\left\langle \sigma_{\lambda}^{\ast
}\theta-dW_{\lambda},\frac{\partial}{\partial\lambda_{i}}\right\rangle
\,d\lambda_{i}=\left(  p_{j}\,\frac{\partial\left(  \hat{\sigma}_{\lambda
}\right)  ^{j}}{\partial\lambda_{i}}-\frac{\partial W_{\lambda}}%
{\partial\lambda_{i}}\right)  \,d\lambda_{i}=\left(  \lambda,p_{j}%
\,\frac{\partial\left(  \hat{\sigma}_{\lambda}\right)  ^{j}}{\partial\lambda
}-\frac{\partial W_{\lambda}}{\partial\lambda}\right) \\
& \\
&  =\left(  \lambda,\partial\left.  \left(  p_{j}\,\left(  \hat{\sigma
}_{\lambda}\right)  ^{j}-W_{\lambda}\right)  \right/  \partial\lambda\right)
,
\end{align*}
with $W_{\lambda}$ such that
\[
p_{j}\,\frac{\partial\left(  \hat{\sigma}_{\lambda}\right)  ^{j}}{\partial
p_{i}}=\frac{\partial W_{\lambda}}{\partial p_{i}},
\]
or equivalently $\left(  \hat{\sigma}_{\lambda}\right)  ^{i}=$ $\partial
\left.  \left(  p_{j}\,\left(  \hat{\sigma}_{\lambda}\right)  ^{j}-W_{\lambda
}\right)  \right/  \partial p_{i}$. And for $\left(  q,p\right)
=\Sigma\left(  p,\lambda\right)  $,
\[
\Phi\left(  q,p\right)  =\Psi\circ\varphi\left(  p,\lambda\right)
=\Psi\left(  \lambda,\partial\left.  \left(  p_{j}\,\left(  \hat{\sigma
}_{\lambda}\right)  ^{j}-W_{\lambda}\right)  \right/  \partial\lambda\right)
=:\left(  Q,P\right)  .
\]
Then, defining $\widetilde{W}_{\lambda}:=p_{j}\,\left(  \hat{\sigma}_{\lambda
}\right)  ^{j}-W_{\lambda}$ [see the first part of $\left(  \ref{can34}%
\right)  $], the map $\Phi$ corresponding to $\Psi=\Psi_{1}$ and $\Psi
=\Psi_{2}$ [see $\left(  \ref{fi12}\right)  $] give rise just to the first and
the second canonical transformation of $\left(  \ref{ct34}\right)  $,
respectively. In the table below, a dictionary relating the classical
terminology and the canonical transformations constructed in the previous
section is presented.%
\[%
\begin{array}
[c]{ll}%
\mathbf{Type\ 1:} & \ \ r=1,\ \ \Psi=\Psi_{1}\\
\mathbf{Type\ 2:} & \ \ r=1,\ \ \Psi=\Psi_{2}\\
\mathbf{Type\ 3:} & \ \ r=2,\ \ \Psi=\Psi_{1}\\
\mathbf{Type\ 4:} & \ \ r=2,\ \ \Psi=\Psi_{2}%
\end{array}
\]

\section{The complete solution - first integral's \emph{duality}}

\label{S4}

\label{fihj}

It is well-known that, in the context of the standard Hamilton-Jacobi Theory,
complete solutions and first integrals are closely related (see for example
\cite{pepin-holo}). We show in this section that the same is true in our
extended setting. More precisely, we show that, given a dynamical system
$\left(  M,X\right)  $ and a fibration $\Pi$, every complete solution of the
$\Pi$-HJE for $\left(  M,X\right)  $ gives rise, around each $m\in M$, to a
set of $l:=\dim M-\dim N$ local first integrals of $X$, and reciprocally,
every set of $l$ first integrals of $X$ (\emph{transversal} to $\Pi$) gives
rise, around each $m\in M$, to a local complete solution of the $\Pi$-HJE for
$\left(  M,X\right)  $. Using these results, we give a sufficient condition
for existence of complete solutions for any dynamical system and any
fibration. At the end of the section we shall study the connection between
complete solutions and the notions of \emph{commutative} and
\emph{non-commutative integrability}.

\subsection{From complete solutions to first integrals}

\label{tsfi}

\subsubsection{The related local momentum maps}

Consider three manifolds $M$, $N$ and $\Lambda$, and a surjective local
diffeomorphism $\Sigma:N\times\Lambda\rightarrow M$. For every $m\in M$,
consider a point $\left(  n,\lambda\right)  \in N\times\Lambda$ and open
neighborhoods $U\subset M$, $V_{N}\subset N$ and $V_{\Lambda}\subset\Lambda$
of $m$, $n$ and $\lambda$, respectively, such that $\Sigma\left(
n,\lambda\right)  =m$ and $\Sigma_{|V_{N}\times V_{\Lambda}}:V_{N}\times
V_{\Lambda}\rightarrow U$ is a diffeomorphism. Define
\begin{equation}
\pi:=\left.  p_{N}\right\vert _{V_{N}\times V_{\Lambda}}\circ\left(  \left.
\Sigma\right\vert _{V_{N}\times V_{\Lambda}}\right)  ^{-1}\ \ \ \text{and}%
\ \ \ F:=\left.  p_{\Lambda}\right\vert _{V_{N}\times V_{\Lambda}}\circ\left(
\left.  \Sigma\right\vert _{V_{N}\times V_{\Lambda}}\right)  ^{-1}.
\label{pif}%
\end{equation}
It is clear that they are submersions and
\begin{equation}
\pi\left(  U\right)  =V_{N}\ \ \ \ \text{and}\ \ \ F\left(  U\right)
=V_{\Lambda}. \label{piu}%
\end{equation}
Then, given $m\in M$, there exists an open neighborhood $U\subset M$ of $m$
and fibrations $\pi:U\rightarrow\pi\left(  U\right)  \subset N$ and
$F:U\rightarrow F\left(  U\right)  \subset\Lambda$ such that $\left.
\Sigma\right\vert _{\pi\left(  U\right)  \times F\left(  U\right)  }%
:\pi\left(  U\right)  \times F\left(  U\right)  \rightarrow U$ is a
diffeomorphism and%
\[
\left(  \pi,F\right)  =\left(  \left.  \Sigma\right\vert _{\pi\left(
U\right)  \times F\left(  U\right)  }\right)  ^{-1},
\]
according to $\left(  \ref{pif}\right)  $ and $\left(  \ref{piu}\right)  $.
The fact that $\left(  \pi,F\right)  $ is a diffeomorphism implies that
\[
TU=\operatorname{Ker}\pi_{\ast}\oplus\operatorname{Ker}F_{\ast},
\]
as stated in the lemma below.

\begin{lemma}
\label{fgsum}Consider three finite dimensional manifolds $M$, $N$ and
$\Lambda$, and two submersions $h:M\rightarrow N$ and $g:M\rightarrow\Lambda$.
Then, $\left(  h,g\right)  $ is a local diffeomorphism if and only if
$TM=\operatorname{Ker}h_{\ast}\oplus\operatorname{Ker}g_{\ast}$. We shall say
in this case that $h$ and $g$ are \textbf{transverse}.
\end{lemma}

\emph{Proof. }Given $m\in M$ and $x\in T_{m}M$, since
\[
\left(  h,g\right)  _{\ast,m}\left(  x\right)  =\left(  h_{\ast,m},g_{\ast
,m}\right)  \left(  x\right)  =\left(  h_{\ast,m}\left(  x\right)  ,g_{\ast
,m}\left(  x\right)  \right)  ,
\]
it is clear that $\left(  h,g\right)  _{\ast,m}$ is injective if and only if
$\operatorname{Ker}h_{\ast,m}\cap\operatorname{Ker}g_{\ast,m}=\left\{
0\right\}  $. On the other hand, since $h$ and $g$ are submersions, it is easy
to prove, just by counting dimensions, that $\left(  h,g\right)  _{\ast,m}$ is
surjective if and only if
\[
\dim M=\dim\left(  \operatorname{Ker}h_{\ast,m}\right)  +\dim\left(
\operatorname{Ker}g_{\ast,m}\right)  .
\]
This completes the proof.$\ \ \ \square$

\bigskip

Suppose now that $\Pi\circ\Sigma=p_{N}$ for some fibration $\Pi:M\rightarrow
N$. It is easy to show that, for any triple $(U,\pi,F)$ as above, we have that
$\pi=\left.  \Pi\right\vert _{U}$. This means that, for every $m\in M$, there
exists an open neighborhood $U\subset M$ of $m$ and a fibration
$F:U\rightarrow F\left(  U\right)  $ such that
\[
\left.  \Sigma\right\vert _{\Pi\left(  U\right)  \times F\left(  U\right)
}:\Pi\left(  U\right)  \times F\left(  U\right)  \rightarrow U
\]
is a diffeomorphism and $\left(  \left.  \Pi\right\vert _{U},F\right)
=\left(  \left.  \Sigma\right\vert _{\Pi\left(  U\right)  \times F\left(
U\right)  }\right)  ^{-1}$. In particular,
\begin{equation}
F=\left.  p_{\Lambda}\right\vert _{\Pi\left(  U\right)  \times F\left(
U\right)  }\circ\left(  \left.  \Sigma\right\vert _{\Pi\left(  U\right)
\times F\left(  U\right)  }\right)  ^{-1}. \label{fpl}%
\end{equation}
Note that, if $\Sigma$ is a global diffeomorphism, we can take $U:=M$ and
$F:=p_{\Lambda}\circ\Sigma^{-1}$.

\begin{remark}
\label{comp}If $N$ is a compact manifold, using standard techniques, $U$ can
be taken such that $\Pi\left(  U\right)  =N$. In other words, there exist $U$
and $F$ such that $\left.  \Sigma\right\vert _{N\times F\left(  U\right)
}:N\times F\left(  U\right)  \rightarrow U$ is a diffeomorphism and $\left(
\left.  \Pi\right\vert _{U},F\right)  =\left(  \left.  \Sigma\right\vert
_{N\times F\left(  U\right)  }\right)  ^{-1}$.
\end{remark}

Finally, assume that $\Sigma$ is a complete solution of the $\Pi$-HJE for some
dynamical system $\left(  M,X\right)  $. Let us show that, in this case, for
any pair $(U,F)$ as given above, we have that $\left.  X\right\vert _{U}%
\in\operatorname{Ker}F_{\ast}$. Because of the definition of $X^{\Sigma}$ [see
Eq. $\left(  \ref{XS}\right)  $],
\[
\left(  p_{\Lambda}\right)  _{\ast}\left(  X^{\Sigma}\left(  n,\lambda\right)
\right)  =0,\ \ \ \forall\left(  n,\lambda\right)  \in N\times\Lambda.
\]
On the other hand, it follows from Eqs. $\left(  \ref{Srel}\right)  $ and
$\left(  \ref{fpl}\right)  $,
\[
F_{\ast}\left(  X\left(  \Sigma\left(  n,\lambda\right)  \right)  \right)
=F_{\ast}\left(  \Sigma_{\ast}\left(  X^{\Sigma}\left(  n,\lambda\right)
\right)  \right)  =\left(  p_{\Lambda}\right)  _{\ast}\left(  X^{\Sigma
}\left(  n,\lambda\right)  \right)  =0
\]
for all $\left(  n,\lambda\right)  \in\Pi\left(  U\right)  \times F\left(
U\right)  $. So, the wanted result follows from the fact that $U=\Sigma\left(
\Pi\left(  U\right)  \times F\left(  U\right)  \right)  $. Summing up, we have
proved the next theorem.

\begin{theorem}
\label{pc}Let $\left(  M,X\right)  $ be a dynamical system, $\Pi:M\rightarrow
N$ a fibration and $\Sigma:N\times\Lambda\rightarrow M$ a complete solution of
the $\Pi$-HJE for $\left(  M,X\right)  $. Then, for every $m\in M$ there exist
an open neighborhood $U$ of $m$ and a fibration $F:U\rightarrow F\left(
U\right)  \subset\Lambda$ such that%
\[
\operatorname{Im}\left.  X\right\vert _{U}\subset\operatorname{Ker}F_{\ast}%
\]
and $F$ is transverse to $\left.  \Pi\right\vert _{U}$, i.e.%
\[
TU=\operatorname{Ker}\left(  \left.  \Pi\right\vert _{U}\right)  _{\ast}%
\oplus\operatorname{Ker}F_{\ast}.
\]
Moreover, $U$ and $F$ can be taken such that

\begin{enumerate}
\item $\Sigma\left(  \Pi\left(  U\right)  \times F\left(  U\right)  \right)
\subset U$,

\item $\left.  \Sigma\right\vert _{\Pi\left(  U\right)  \times F\left(
U\right)  }$ is a global diffeomorphism onto $U$,

\item and
\begin{equation}
F\circ\left.  \Sigma\right\vert _{\Pi\left(  U\right)  \times F\left(
U\right)  }=\left.  p_{\Lambda}\right\vert _{\Pi\left(  U\right)  \times
F\left(  U\right)  }. \label{fsi}%
\end{equation}

\end{enumerate}
\end{theorem}

In the context of the previous Theorem, suppose that $\Lambda$ is an open
subset of $\mathbb{R}^{l}$. If we call $f_{i}$ to each component of $F$, the
fact that $\operatorname{Im}\left.  X\right\vert _{U}\subset\operatorname{Ker}%
F_{\ast}$ is equivalent to $\left\langle df_{i},\left.  X\right\vert
_{U}\right\rangle =0$ for all $i=1,...,l$, i.e. each $f_{i}$ is a first
integral for $\left.  X\right\vert _{U}$. We can see that from another point
of view. We know from Theorem \ref{complete} that every integral curve
$\Gamma:I\rightarrow M$ of $X$ is given by the formula $\Gamma\left(
t\right)  =\Sigma\left(  \gamma\left(  t\right)  ,\lambda\right)  $ for some
$\lambda\in\Lambda$ and some integral curve $\gamma$ of $X^{\sigma_{\lambda}}$
[see $\left(  \ref{Ml}\right)  $]. Taking an appropriate pair $U$ and $F$ as
in the previous theorem, Eq. $\left(  \ref{fsi}\right)  $ implies that
\[
F\left(  \Gamma\left(  t\right)  \right)  =F\left(  \Sigma\left(
\gamma\left(  t\right)  ,\lambda\right)  \right)  =p_{\Lambda}\left(
\gamma\left(  t\right)  ,\lambda\right)  =\lambda,
\]
for all $t$ such that $\Gamma\left(  t\right)  \in U$. Concluding, related to
any complete solution $\Sigma:N\times\Lambda\rightarrow M$ of the $\Pi$-HJE
for the dynamical system $\left(  M,X\right)  $, we have $\dim\Lambda$ (local)
constant of motion for $\left(  M,X\right)  $. This motives us to the next definition.

\begin{definition}
\label{setlfi}Let $\left(  M,X\right)  $ be a dynamical system, $\Pi
:M\rightarrow N$ a fibration and $\Sigma:N\times\Lambda\rightarrow M$ a
complete solution of the $\Pi$-HJE for $\left(  M,X\right)  $. A pair $\left(
U,F\right)  $, with $U$ an open submanifold of $M$ and $F:U\rightarrow
F\left(  U\right)  \subset\Lambda$ a fibration satisfying the items $1$, $2$
and $3$ of the last theorem, will be called a \textbf{local momentum map} of
$\left(  M,X\right)  $ \textbf{related} to $\Sigma$. If $\Sigma$ is a global
diffeomorphism, we shall call $F:=p_{\Lambda}\circ\Sigma^{-1}:M\rightarrow
\Lambda$ the \textbf{global momentum map} of $\left(  M,X\right)  $\textbf{
related }to $\Sigma$.
\end{definition}

It is clear that a local momentum map $\left(  U,F\right)  $ of $\left(
M,X\right)  $ related to $\Sigma$ is the global momentum map $F:U\rightarrow
F\left(  U\right)  $ of $\left(  U,\left.  X\right\vert _{U}\right)  $ related
to $\left.  \Sigma\right\vert _{\Pi\left(  U\right)  \times F\left(  U\right)
}$.

\begin{remark}
\label{setlfi2}Given a pair $\left(  U,F\right)  $ as in definition above, the
level surfaces of $F$ are
\[
F^{-1}\left(  \lambda\right)  =M_{\lambda}\cap U,\ \ \ \forall\lambda\in
F\left(  U\right)  ,
\]
where $M_{\lambda}=\Sigma\left(  N\times\left\{  \lambda\right\}  \right)
=\operatorname{Im}\sigma_{\lambda}$ [see $\left(  \ref{Ml}\right)  $]. Of
course, if $\Sigma$ is a global diffeomorphism, we have for the global
momentum map that $F^{-1}\left(  \lambda\right)  =M_{\lambda}$, $\forall
\lambda\in\Lambda$. In any case,
\begin{equation}
\operatorname{Ker}F_{\ast,m}=T_{m}\left(  M_{F\left(  m\right)  }\right)
\label{kerf}%
\end{equation}
for all $m$ in the domain of $F$. On the other hand, if $N$ is compact (see
Remark \ref{comp}), the pairs $\left(  U,F\right)  $ can be chosen such that
\[
F^{-1}\left(  \lambda\right)  =M_{\lambda},\ \ \ \forall\lambda\in F\left(
U\right)  .
\]

\end{remark}

It is well-known that, when the number of first integrals is $l=d-1$, the
system can be integrated by quadratures. This can be deduced from the last
theorem and a comment we made at the end of the Section \ref{totsol} [see Eq.
$\left(  \ref{ibq1}\right)  $].

\subsubsection{Some commutation relations}

\label{scr}

Consider a Leibniz manifold $\left(  M,\Xi\right)  $ (see Section
\ref{almost}), a function $H:M\rightarrow\mathbb{R}$, a fibration
$\Pi:M\rightarrow N$ and a complete solution $\Sigma:N\times\Lambda\rightarrow
M$ of the $\Pi$-HJE for $\left(  M,X_{H}\right)  $. Consider also a local
momentum map $\left(  U,F\right)  $ related to $\Sigma$, and assume that
$\Lambda$, and consequently\emph{ }$F\left(  U\right)  $, is an open subset of
$\mathbb{R}^{l}$. Denote by $f_{i}$ the $i$-th component of $F$. We shall
study conditions on $\Sigma$ under which the functions $f_{i}$'s are in
involution. We are mainly interested on the symplectic and the Poisson cases,
because of the connection between involutivity and integrability that exists
in those cases. The next proposition is an immediate consequence of
Proposition \ref{equio} and the Eq. $\left(  \ref{kerf}\right)  $.

\begin{proposition}
\label{equi}Consider a Leibniz manifold $\left(  M,\Xi\right)  $, a function
$H:M\rightarrow\mathbb{R}$, a fibration $\Pi:M\rightarrow N$ and a complete
solution $\Sigma:N\times\Lambda\rightarrow M$ of the $\Pi$-HJE for $\left(
M,X_{H}\right)  $. $\Sigma$ is (weakly) (co)isotropic if and only if, for each
local momentum map $\left(  U,F\right)  $ related to $\Sigma$, the fibration
$F:U\rightarrow F\left(  U\right)  $ is (weakly) (co)isotropic. And if
$\Sigma$ is a global diffeomorphism, $\Sigma$ is (weakly) (co)isotropic if and
only if the global momentum map related to $\Sigma$ is (weakly) (co)isotropic.
\end{proposition}

Since $\operatorname{Ker}F_{\ast}=\left\langle df_{1},...,df_{l}\right\rangle
^{0}$ and $\Xi^{\sharp}\left[  \left(  \operatorname{Ker}F_{\ast}\right)
^{0}\right]  =\left\langle X_{f_{1}},...,X_{f_{l}}\right\rangle $, from the
identities
\[
\left\{  f_{i},f_{j}\right\}  =\Xi\left(  df_{i},df_{j}\right)  =\left\langle
df_{i},\Xi^{\sharp}\left(  df_{j}\right)  \right\rangle =\left\langle
df_{i},X_{f_{j}}\right\rangle ,
\]
we have that the functions $f_{i}$'s are in involution if and only if
$\Xi^{\sharp}\left[  \left(  \operatorname{Ker}F_{\ast}\right)  ^{0}\right]
\subset\operatorname{Ker}F_{\ast}$. Using Proposition \ref{equi}, this is true
if $\Sigma$ is co-isotropic, and in particular if it is (weakly) Lagrangian.
(This result includes the analogous ones obtained in Refs. \cite{pepin-noholo}
and \cite{hjp} for almost-Poisson manifolds.) In the case of symplectic
manifolds, unless in some special cases, we can deduce that from the way in
which the brackets $\left\{  f_{i},f_{j}\right\}  $ depend on the partial
solutions $\sigma_{\lambda}$'s. From now on, we shall assume that $\Xi
^{\sharp}=\omega^{\sharp}$ for some symplectic form $\omega$ on $M$.

\begin{proposition}
Under above notation, if $\Pi$ is isotropic w.r.t. $\omega$ [recall Eq.
$\left(  \ref{isot}\right)  $], then
\begin{equation}
\left\{  f_{i},f_{j}\right\}  \left(  m\right)  =\left(  \sigma_{F\left(
m\right)  }\right)  ^{\ast}\omega\left(  \Pi_{\ast,m}\left(  X_{f_{i}}\left(
m\right)  \right)  ,\Pi_{\ast,m}\left(  X_{f_{j}}\left(  m\right)  \right)
\right)  ,\ \ \ \forall m\in U, \label{bf}%
\end{equation}
for all $i,j=1,...,l$. In particular, in the standard situation (see
Definition \ref{standd}),%
\[
\left\{  f_{i},f_{j}\right\}  \left(  m\right)  =\left(  \sigma_{F\left(
m\right)  }\right)  ^{\ast}\omega\left(  \mathbb{F}f_{i}\left(  m\right)
,\mathbb{F}f_{j}\left(  m\right)  \right)  ,\ \ \ \forall m\in U.
\]

\end{proposition}

\emph{Proof.} To simplify the notation, assume for simplicity that $U=M$.
Given $m\in M$,%
\[
X_{f_{i}}\left(  m\right)  -\left(  \sigma_{F\left(  m\right)  }\right)
_{\ast,\Pi\left(  m\right)  }\circ\Pi_{\ast,m}\left(  X_{f_{i}}\left(
m\right)  \right)  \in\operatorname{Ker}\Pi_{\ast,m},
\]
for all $i=1,...,l$, so%
\[
\omega\left(  X_{f_{i}}\left(  m\right)  -\left(  \sigma_{F\left(  m\right)
}\right)  _{\ast,\Pi\left(  m\right)  }\circ\Pi_{\ast,m}\left(  X_{f_{i}%
}\left(  m\right)  \right)  ,X_{f_{j}}\left(  m\right)  -\left(
\sigma_{F\left(  m\right)  }\right)  _{\ast,\Pi\left(  m\right)  }\circ
\Pi_{\ast,m}\left(  X_{f_{j}}\left(  m\right)  \right)  \right)  =0,
\]
for all $i,j=1,...,l$. Note that, since $F\circ\Sigma\left(  n,\lambda\right)
=\lambda$, i.e. $F\circ\sigma_{\lambda}:N\rightarrow\mathbb{R}$ is a constant
function, the same is true for each component $f_{i}\circ\sigma_{\lambda}$.
This ensures that%
\[%
\begin{array}
[c]{l}%
\omega\left(  X_{f_{i}}\left(  m\right)  ,\left(  \sigma_{F\left(  m\right)
}\right)  _{\ast,\Pi\left(  m\right)  }\circ\Pi_{\ast,m}\left(  X_{f_{j}%
}\left(  m\right)  \right)  \right)  =\left\langle df_{i}\left(  m\right)
,\left(  \sigma_{F\left(  m\right)  }\right)  _{\ast,\Pi\left(  m\right)
}\circ\Pi_{\ast,m}\left(  X_{f_{j}}\left(  m\right)  \right)  \right\rangle
=\\
\\
=\left\langle \left(  \sigma_{F\left(  m\right)  }\right)  _{\Pi\left(
m\right)  }^{\ast}\left(  df_{i}\left(  m\right)  \right)  ,\Pi_{\ast
,m}\left(  X_{f_{j}}\left(  m\right)  \right)  \right\rangle =\left\langle
d\left(  f_{i}\circ\sigma_{F\left(  m\right)  }\right)  \left(  \Pi\left(
m\right)  \right)  ,\Pi_{\ast,m}\left(  X_{f_{j}}\left(  m\right)  \right)
\right\rangle =0,
\end{array}
\]
from which the first part of the proposition easily follows. For the standard
situation, recall Eq. $\left(  \ref{pifh}\right)  $.\ \ \ $\square$

\bigskip

Let us characterize the case in which the functions $f_{i}$'s are in
involution (unless when $\Pi$ is an isotropic fibration).

\begin{proposition}
Under the conditions of the last proposition, the functions $f_{i}$'s are in
involution if and only if the partial solutions $\sigma_{\lambda}$'s are
Lagrangian for all $\lambda\in F\left(  U\right)  $. In the standard
situation, the functions $f_{i}$'s are in involution if and only if each
$1$-form $\sigma_{\lambda}$ is closed. In any case, $\Pi$ and $F$ must also be Lagrangian.
\end{proposition}

\emph{Proof.} If the functions $f_{i}$'s are in involution, then $F$ is
co-isotropic, i.e.%
\[
\omega^{\sharp}\left[  \left(  \operatorname{Ker}F_{\ast}\right)  ^{0}\right]
=\left\langle X_{f_{1}},...,X_{f_{l}}\right\rangle \subset\operatorname{Ker}%
F_{\ast}.
\]
So, taking into account that $F$ is transverse to $\Pi$, it is easy to show
that the subset of vectors
\[
\left\{  \Pi_{\ast,m}\left(  X_{f_{1}}\left(  m\right)  \right)
,...,\Pi_{\ast,m}\left(  X_{f_{l}}\left(  m\right)  \right)  \right\}  \subset
T_{\Pi\left(  m\right)  }N
\]
is l.i. (linearly independent). Using Eq. $\left(  \ref{bf}\right)  $, the
last observation implies that $\sigma_{\lambda}^{\ast}\omega=0$ (i.e.
$\sigma_{\lambda}$ is isotropic) for all $\lambda\in F\left(  U\right)  $. But
if $\sigma_{\lambda}$ is isotropic, since the same is assumed for $\Pi$, then
both $\Pi$ and $\sigma_{\lambda}$ must be Lagrangian (because the dimensions
of $\operatorname{Ker}\Pi_{\ast}$ and $\operatorname{Im}\left(  \sigma
_{\lambda}\right)  _{\ast}$ are complementary). This proves the first
affirmation. The converse follows immediately from Eq. $\left(  \ref{bf}%
\right)  $. For the standard situation, recall Eq. $\left(  \ref{ests}\right)
$ of Section \ref{canej}. \ \ \ $\square$

\bigskip

It is clear that the last proposition establishes a deep connection between
Lagrangian complete solutions and \emph{commutative} (or
\emph{Liouville-Arnold})\emph{ integrability} \cite{ar}. In the standard
situation, it recovers the following well-know result: \emph{closedness of
partial solutions is equivalent to the involutivity of the related first
integrals}. In Section \ref{INC}, we shall extend the above connection to
general Poisson manifolds, and even to the notion of \emph{non-commutative
}(or\emph{ Mischenko-Fomenko})\emph{ integrability} \cite{mf}.

\subsection{From first integrals to complete solutions}

\label{fits}

In this section, we are going to show a result which can be seen as a converse
of Theorem \ref{pc}. Let us fix a dynamical system $\left(  M,X\right)  $ and
a fibration $\Pi:M\rightarrow N$.

\begin{theorem}
\label{inthj}Let $l:=\dim M-\dim N$, $\Lambda$ an $l$-manifold and
$F:M\rightarrow\Lambda$ a fibration such that%
\[
\operatorname{Im}X\subset\operatorname{Ker}F_{\ast}%
\]
and $F$ and $\Pi$ are transverse, i.e.
\begin{equation}
TM=\operatorname{Ker}\Pi_{\ast}\oplus\operatorname{Ker}F_{\ast}. \label{ker}%
\end{equation}
Then, for every $m\in M$, there exists an open neighborhood $U$ of $m$ such
that $\left.  \left(  \Pi,F\right)  \right\vert _{U}:U\rightarrow\Pi\left(
U\right)  \times F\left(  U\right)  $ is a diffeormorphism and
\[
\Sigma:=\left(  \left.  \left(  \Pi,F\right)  \right\vert _{U}\right)  ^{-1}%
\]
is a complete solution of the $\left.  \Pi\right\vert _{U}$-HJE for $\left(
U,\left.  X\right\vert _{U}\right)  $. Moreover, $\left.  F\right\vert _{U}$
is the global momentum map related to $\Sigma$.
\end{theorem}

\emph{Proof.} Since $\Pi$ and $F$ are submersions and condition $\left(
\ref{ker}\right)  $ holds, we know from Lemma \ref{fgsum} that $\left(
\Pi,F\right)  $ is a local diffeomorphism. Let $U\subset M$ such that%
\[
\left.  \left(  \Pi,F\right)  \right\vert _{U}:U\rightarrow\Pi\left(
U\right)  \times F\left(  U\right)
\]
is a diffeomorphism, and define $\Sigma:=\left(  \left.  \left(  \Pi,F\right)
\right\vert _{U}\right)  ^{-1}:\Pi\left(  U\right)  \times F\left(  U\right)
\rightarrow U$. We have to prove that $\Sigma$ is a complete solution of the
$\left.  \Pi\right\vert _{U}$-HJE for $\left(  U,\left.  X\right\vert
_{U}\right)  $. To do that, we shall show that $\left(  \ref{Srel}\right)  $
holds for $\Sigma$ along $\Pi\left(  U\right)  \times F\left(  U\right)  $.
Let us fix $\left(  n,\lambda\right)  \in\Pi\left(  U\right)  \times F\left(
U\right)  $. Then $\left(  \Pi,F\right)  \left(  \Sigma\left(  n,\lambda
\right)  \right)  =\left(  n,\lambda\right)  $ and, applying $p_{N}$ on both
sides, $\Pi\left(  \Sigma\left(  n,\lambda\right)  \right)  =n$, which proves
the first part of $\left(  \ref{Srel}\right)  $. On the other hand, the
condition $\operatorname{Im}X\subset\operatorname{Ker}F_{\ast}$ implies that
\[
\left(  \Pi_{\ast},F_{\ast}\right)  \left(  X\left(  \Sigma\left(
n,\lambda\right)  \right)  \right)  =\left(  \Pi_{\ast}\left(  X\left(
\Sigma\left(  n,\lambda\right)  \right)  \right)  ,0\right)  =X^{\Sigma
}\left(  n,\lambda\right)  .
\]
Thus, $\Sigma_{\ast}\left(  X^{\Sigma}\left(  n,\lambda\right)  \right)
=\Sigma_{\ast}\left[  \left(  \Pi_{\ast},F_{\ast}\right)  \left(  X\left(
\Sigma\left(  n,\lambda\right)  \right)  \right)  \right]  =X\left(
\Sigma\left(  n,\lambda\right)  \right)  $, which implies the second part of
$\left(  \ref{Srel}\right)  $.\ \ \ $\square$

\bigskip

This theorem says that, given a set of first integrals for a dynamical system
$\left(  M,X\right)  $, defined by a fibration $F:M\rightarrow\Lambda$
transverse to $\Pi$, we can construct local complete solutions of the $\Pi
$-HJE for $\left(  M,X\right)  $ around every point $m\in M$.

\begin{example}
Suppose that $M=T^{\ast}Q$, $N=Q$, $\Pi=\pi_{Q}$ and $\Lambda\subset
\mathbb{R}^{l}$ is an open subset. Let us call $f_{i}$ the components of the
fibration $F:T^{\ast}Q\rightarrow\Lambda$. Then, transversality condition
$\left(  \ref{ker}\right)  $ in this case means that, for all $m\in T^{\ast}%
Q$, the set of fiber derivatives [recall Eq. $\left(  \ref{fd}\right)  $]
$\left\{  \mathbb{F}f_{1}\left(  m\right)  ,...,\mathbb{F}f_{l}\left(
m\right)  \right\}  \subset T_{\pi_{Q}\left(  m\right)  }Q$ is linearly independent.
\end{example}

If no fibration $\Pi$ is considered \emph{a priori}, we have the next result.

\begin{theorem}
\label{fint}Let $\Lambda$ be an $l$-manifold, with $l\leq d$, and
$F:M\rightarrow\Lambda$ a fibration such that $\operatorname{Im}%
X\subset\operatorname{Ker}F_{\ast}$. Then, for every $m\in M$, there exist an
open neighborhood $U$ of $m$ and a fibration $\pi:U\rightarrow\pi\left(
U\right)  \subset\mathbb{R}^{d-l}$ such that
\[
\left(  \pi,\left.  F\right\vert _{U}\right)  :U\rightarrow\pi\left(
U\right)  \times F\left(  U\right)
\]
is diffeormorphism and $\Sigma:=\left(  \pi,\left.  F\right\vert _{U}\right)
^{-1}$ is a complete solution of the $\pi$-HJE for $\left(  U,\left.
X\right\vert _{U}\right)  $. Moreover, $\left.  F\right\vert _{U}$ is the
global momentum map related to $\Sigma$.
\end{theorem}

\emph{Proof.} Consider, around a point $m\in M$, coordinates charts adapted of
$M$ and $\Lambda$ adapted to the fibration $F$, i.e. charts $\varphi:U\subset
M\rightarrow\varphi\left(  U\right)  \subset\mathbb{R}^{d}$ and $\psi
:\Lambda_{1}\subset\Lambda\rightarrow\psi\left(  \Lambda_{1}\right)
\subset\mathbb{R}^{l}$ such that $m\in U$, $F\left(  U\right)  \subset
\Lambda_{1}$ and, for all $\left(  x_{1},...,x_{d}\right)  \in\varphi\left(
U\right)  $,
\[
\psi\left[  F\left(  \varphi^{-1}\left(  x_{1},...,x_{d}\right)  \right)
\right]  =\left(  x_{d-l+1},...,x_{d}\right)  .
\]
Defining $\pi:U\rightarrow\mathbb{R}^{d-l}$ such that%
\[
\pi\left(  \varphi^{-1}\left(  x_{1},...,x_{d}\right)  \right)  =\left(
x_{1},...,x_{d-l}\right)  ,
\]
it is clear that
\[
TU=\operatorname{Ker}\pi_{\ast}\oplus\operatorname{Ker}\left(  \left.
F\right\vert _{U}\right)  _{\ast}.
\]
The rest is a consequence of the last theorem applied to $\left(  U,\left.
X\right\vert _{U}\right)  $, $\pi$ and $\left.  F\right\vert _{U}%
$.\ \ \ $\square$

\bigskip

In other words, related to a set of first integrals of a dynamical system
$\left(  M,X\right)  $, we always have, unless locally, a complete solution of
the $\pi$-HJE for $\left(  M,X\right)  $, for some fibration $\pi$.

\begin{remark}
\label{adap}It is worth mentioning that, given a fibration $F$, its adapted
charts can be found by quadratures from $F$ (actually, by linear calculations
only). Then, the fibration $\pi$ defined in the last proposition, can be found
by quadratures too.
\end{remark}

\subsection{Some applications of the duality}

Fix again a dynamical system $\left(  M,X\right)  $ and a fibration
$\Pi:M\rightarrow N$. We shall apply the previous results about the duality
between complete solutions and first integrals in different situations.

\subsubsection{Existence of (local) complete solutions}

\label{exist}

A sufficient condition for existence of complete solutions will be presented
in this section. As above, write $l=\dim M-\dim N$ and $d=\dim M$, and note
that $\dim\left(  \operatorname{Ker}\Pi_{\ast}\right)  =l$.

\begin{proposition}
\label{ex}Suppose that $X\left(  m\right)  \notin\operatorname{Ker}\Pi_{\ast}$
on some point $m\in M$. Then, there exists an open neighborhood $U$ of $m$ and
a fibration $F:U\rightarrow F\left(  U\right)  \subset\mathbb{R}^{l}$ such
that
\[
\operatorname{Im}\left.  X\right\vert _{U}\subset\operatorname{Ker}F_{\ast
}\ \ \ \text{and\ \ \ }TU=\operatorname{Ker}\left(  \left.  \Pi\right\vert
_{U}\right)  _{\ast}\oplus\operatorname{Ker}F_{\ast}.
\]

\end{proposition}

In order to prove that, we need the following Lemma.

\begin{lemma}
Let $V$ be an $n$-dimensional space, $S\subset V$ a $k$-dimensional subspace
and $B=\left\{  v_{1},...,v_{n}\right\}  $ a basis of $V$. Suppose that
$v_{1}\notin S$. Then, there exists $B_{1}\subset B$ containing $v_{1}$ such
that $B_{1}$ is a basis of a complement of $S$.
\end{lemma}

\emph{Proof.} Let us consider the quotient space $\left.  V\right/  S$ and the
classes $\left[  v_{i}\right]  $'s. Since $B$ is a basis of $V$, the elements
$\left[  v_{i}\right]  \in\left.  V\right/  S$ generate $\left.  V\right/  S$.
It is well-known that, since $\left[  v_{1}\right]  \neq0$ (because
$v_{1}\notin S$), we can extract from the set of generators $\left\{  \left[
v_{1}\right]  ,...,\left[  v_{n}\right]  \right\}  $ a basis of $\left.
V\right/  S$ containing $\left[  v_{1}\right]  $. Let $\left\{  \left[
v_{i_{1}}\right]  ,...,\left[  v_{i_{k}}\right]  \right\}  $, with $v_{i_{1}%
}=v_{1}$, such a basis. Then, $B_{1}:=\left\{  v_{i_{1}},...,v_{i_{k}%
}\right\}  $ is the subset we are looking for.\ \ \ $\square$

\bigskip

\emph{Proof of Proposition \ref{ex}.} Since $X\left(  m\right)  \notin
\operatorname{Ker}\Pi_{\ast}$, and consequently $X\left(  m\right)  \neq0$,
using the Rectification Theorem (see, for instance \cite{ar2}) for vector
fields we know that there exists an open neighborhood $U$ of $m$ and a chart
\[
\varphi=\left(  \varphi_{1},...,\varphi_{d}\right)  :U\rightarrow
\varphi\left(  U\right)  \subset\mathbb{R}^{d}%
\]
such that
\[
\left.  X\right\vert _{U}=\frac{\partial}{\partial\varphi_{1}}.
\]
Applying the previous Lemma to $V=T_{m}M$, $S=\operatorname{Ker}\Pi_{\ast,m}$
and $v_{i}=\left.  \partial\right/  \partial\varphi_{i}\left(  m\right)  $, we
can construct a complement of $\operatorname{Ker}\Pi_{\ast,m}$, generated by
some of the vectors $\left.  \partial\right/  \partial\varphi_{i}\left(
m\right)  $, and containing $\left.  \partial\right/  \partial\varphi
_{1}\left(  m\right)  =\left.  X\right\vert _{U}\left(  m\right)  $. Since
$\dim\left(  \operatorname{Ker}\Pi_{\ast,m}\right)  =l$, we shall need $d-l$
of such vectors, including $\left.  \partial\right/  \partial\varphi
_{1}\left(  m\right)  $. Let us re-order the components of $\varphi$ in such a
way that those vectors are%
\[
\frac{\partial}{\partial\varphi_{1}}\left(  m\right)  ,...,\frac{\partial
}{\partial\varphi_{d-l}}\left(  m\right)  .
\]
Then,%
\[
\left\langle \frac{\partial}{\partial\varphi_{1}}\left(  m\right)
,...,\frac{\partial}{\partial\varphi_{d-l}}\left(  m\right)  \right\rangle
\cap\operatorname{Ker}\Pi_{\ast,m}=\left\{  0\right\}  .
\]
By continuity, shrinking $U$ if it is necessary, previous condition is true
not only for $m$, but for all the points inside $U$. Now, define
$F:U\rightarrow\mathbb{R}^{l}$ such that%
\[
F\left(  \varphi^{-1}\left(  x_{1},...,x_{d}\right)  \right)  =\left(
x_{d-l+1},...,x_{d}\right)  .
\]
It is easy to see that the vector fields $\left.  \partial\right/
\partial\varphi_{1},...,\left.  \partial\right/  \partial\varphi_{d-l}$ define
a basis for $\operatorname{Ker}F_{\ast}$. In particular,
\[
\operatorname{Im}\left.  X\right\vert _{U}=\operatorname{Im}\frac{\partial
}{\partial\varphi_{1}}\subset\operatorname{Ker}F_{\ast}%
\]
and $\dim\left(  \operatorname{Ker}F_{\ast}\right)  =d-l$. On the other hand,
$\operatorname{Ker}\left(  \left.  \Pi\right\vert _{U}\right)  _{\ast}%
\cap\operatorname{Ker}F_{\ast}=\left\{  0\right\}  $. All that proves the
proposition.\ \ \ $\square$

\bigskip

Combining Theorem \ref{inthj} and the last proposition, the next result is immediate.

\begin{theorem}
Suppose that, on some point $m\in M$,
\begin{equation}
X\left(  m\right)  \notin\operatorname{Ker}\Pi_{\ast}. \label{nonker}%
\end{equation}
Then, there exists an open neighborhood $U$ of $m$, an $l$-manifold $\Lambda$,
with $l=\dim M-\dim N$, and a complete solution $\Sigma:\Pi\left(  U\right)
\times\Lambda\rightarrow U$ of the $\left.  \Pi\right\vert _{U}$-HJE for
$\left(  U,\left.  X\right\vert _{U}\right)  $.
\end{theorem}

Note that the condition $\left(  \ref{nonker}\right)  $ is exactly the
condition that we gave in Section \ref{hjemain} in order to ensure the
existence of a partial solution [see Eq. $\left(  \ref{nonker0}\right)  $].

\begin{example}
Consider the standard situation and assume that $H$ is simple, i.e. of
\emph{kinetic plus potential terms} form. Since $\left(  \pi_{Q}\right)
_{\ast}\circ X_{H}=\mathbb{F}H$, and $\mathbb{F}H$ is non-degenerate, the last
theorem ensures the existence of complete solutions around every point of
$T^{\ast}Q$, except for the null co-distribution.
\end{example}

\subsubsection{{Restriction} of a complete solutions to an invariant
submanifold}

\label{rest}

Let $\Sigma:N\times\Lambda\rightarrow M$ be a complete solution of the $\Pi
$-HJE for $\left(  M,X\right)  $, and let $S$ be an $X$-invariant submanifold
of $M$. We shall denote by $\left.  X\right\vert _{S}:S\rightarrow TS$ the
corresponding vector field obtained by the restriction of $X$ to $S$. To take
advantage of both $\Sigma$ and $S$ (for finding the trajectories of $X$
intersecting $S$), it would be useful to \textquotedblleft
restrict\textquotedblright\ somehow $\Sigma$ to $S$, in order to obtain a
(local) complete solution for the dynamical system $\left(  S,\left.
X\right\vert _{S}\right)  $. Assuming certain regularity conditions, we can
make such a restriction. To state that result precisely, consider the next definition.

\begin{definition}
\label{resr}Under the previous notation, we shall say that $\Sigma$
\textbf{restricts to }$S$ if for every $s_{0}\in S$ there exists an open
neighborhood $R\subset S$ of $s_{0}$ and a submanifold $\Lambda_{1}%
\subset\Lambda$ such that $\Pi\left(  R\right)  \subset N$ is a submanifold,
$\left.  \Pi\right\vert _{R}:R\rightarrow\Pi\left(  R\right)  $ is a fibration
and
\[
\left.  \Sigma\right\vert _{\Pi\left(  R\right)  \times\Lambda_{1}}:\Pi\left(
R\right)  \times\Lambda_{1}\rightarrow R
\]
is a complete solution of the $\left.  \Pi\right\vert _{R}$-HJE for $\left(
R,\left.  X\right\vert _{R}\right)  $. We shall call $\left.  \Sigma
\right\vert _{\Pi\left(  R\right)  \times\Lambda_{1}}$ a \textbf{local
restriction of }$\Sigma$ to $S$\textbf{ around }$s_{0}$.
\end{definition}

Now, the announced result.

\begin{proposition}
\label{pq}Let $\Sigma:N\times\Lambda\rightarrow M$ be a complete solution of
the $\Pi$-HJE for $\left(  M,X\right)  $, and let $S$ be an $X$-invariant
submanifold of $M$. Suppose that:

\begin{enumerate}
\item for every local momentum map $\left(  U,F\right)  $ related to $\Sigma$,
$\left.  F\right\vert _{S\cap U}:S\cap U\rightarrow F\left(  U\right)  $ has
constant rank $r_{1}$;

\item the restriction $\left.  \Pi\right\vert _{S}:S\rightarrow N$ has
constant rank $r_{2}$;

\item $r_{1}+r_{2}=\dim S$.
\end{enumerate}

Then, $\Sigma$ restricts to $S$.
\end{proposition}

\emph{Proof.} Given $s_{0}\in S$, consider a local momentum map $\left(
U,F\right)  $ related to $\Sigma$ (see Definition \ref{setlfi}) such that
$s_{0}\in U$. Recall that (see Theorem \ref{pc})
\begin{equation}
\operatorname{Im}\left.  X\right\vert _{U}\subset\operatorname{Ker}F_{\ast}
\label{hc}%
\end{equation}
and%
\begin{equation}
TU=\operatorname{Ker}\left(  \left.  \Pi\right\vert _{U}\right)  _{\ast}%
\oplus\operatorname{Ker}F_{\ast}. \label{dsum}%
\end{equation}
Since, for all $s\in S\cap U$,%
\begin{equation}
\operatorname{Ker}\left(  \left.  \Pi\right\vert _{S\cap U}\right)  _{\ast
,s}=\operatorname{Ker}\Pi_{\ast,s}\cap T_{s}S\ \ \ \text{and\ \ \ }%
\operatorname{Ker}\left(  \left.  F\right\vert _{S\cap U}\right)  _{\ast
,s}=\operatorname{Ker}F_{\ast,s}\cap T_{s}S, \label{inter}%
\end{equation}
then
\begin{equation}
\operatorname{Ker}\left(  \left.  \Pi\right\vert _{S\cap U}\right)  _{\ast
,s}\cap\operatorname{Ker}\left(  \left.  F\right\vert _{S\cap U}\right)
_{\ast,s}=0, \label{inter2}%
\end{equation}
and using the condition $3$ of the proposition, it is clear that%
\begin{equation}
T_{s}S=\operatorname{Ker}\left(  \left.  \Pi\right\vert _{S\cap U}\right)
_{\ast,s}\oplus\operatorname{Ker}\left(  \left.  F\right\vert _{S\cap
U}\right)  _{\ast,s},\ \ \ \forall s\in S\cap U. \label{ssum}%
\end{equation}
Also, using the conditions $1$ and $2$, from the Constant Rank Theorem it
follows that, for any point of $S\cap U$, and in particular for the point
$s_{0}\in S\cap U$ given above, there exists an open neighborhood $V\subset
S\cap U$ of $s_{0}$ in $S$ such that $\Pi\left(  V\right)  $ and $F\left(
V\right)  $ are submanifolds of $N$ and $\Lambda$, respectively, and $\left.
\Pi\right\vert _{V}:V\rightarrow\Pi\left(  V\right)  $ and $\left.
F\right\vert _{V}:V\rightarrow F\left(  V\right)  $ are fibrations. Moreover,
from Eq. $\left(  \ref{ssum}\right)  $ it follows that%
\[
TV=\operatorname{Ker}\left(  \left.  \Pi\right\vert _{V}\right)  _{\ast}%
\oplus\operatorname{Ker}\left(  \left.  F\right\vert _{V}\right)  _{\ast}.
\]
On the other hand, using Eq. $\left(  \ref{hc}\right)  $ and the fact that $S$
is $X$-invariant,
\[
\operatorname{Im}\left(  \left.  X\right\vert _{V}\right)  \subset
\operatorname{Ker}F_{\ast}\cap TV=\operatorname{Ker}\left(  \left.
F\right\vert _{V}\right)  _{\ast}.
\]
Then, Theorem \ref{inthj} ensures that $\left(  \left.  \Pi\right\vert
_{V},\left.  F\right\vert _{V}\right)  $ is a local diffeomorphism and that
there exists an open neighborhood $R\subset V$ of $s_{0}$ such that $\left(
\left.  \Pi\right\vert _{R},\left.  F\right\vert _{R}\right)  :R\rightarrow
\Pi\left(  R\right)  \times F\left(  R\right)  $ is a global diffeomorphism
and
\[
\Sigma_{R}:=\left(  \left.  \Pi\right\vert _{R},\left.  F\right\vert
_{R}\right)  ^{-1}:\Pi\left(  R\right)  \times F\left(  R\right)  \rightarrow
R
\]
is a solution of the $\left.  \Pi\right\vert _{R}$-HJE for $\left(  R,\left.
X\right\vert _{R}\right)  $. Of course, $\Sigma_{R}=\left.  \Sigma\right\vert
_{\Pi\left(  R\right)  \times F\left(  R\right)  }$. Choosing $\Lambda
_{1}:=F\left(  R\right)  $, it follows that $\Sigma_{R}$ is a local
restriction of $\Sigma$ to $S$ around $s_{0}$.\ \ \ $\square$

\bigskip

\begin{corollary}
Let $\Sigma:N\times\Lambda\rightarrow M$ be a complete solution of the $\Pi
$-HJE for $\left(  M,X\right)  $, and let $S$ be an $X$-invariant submanifold
of $M$. If for every momentum map $\left(  U,F\right)  $ related to $\Sigma$
we have that $\operatorname{Ker}F_{\ast,s}\subset T_{s}S\ $for all $s\in S\cap
U$, then $\Sigma$ restricts to $S$.
\end{corollary}

\emph{Proof.} It is enough to show that the items $1$, $2$ and $3$ of the last
proposition are fulfilled. In fact, if such condition above holds, the second
part of $\left(  \ref{inter}\right)  $ implies that $\operatorname{Ker}\left(
\left.  F\right\vert _{S\cap U}\right)  _{\ast,s}=\operatorname{Ker}F_{\ast
,s}$. So, $\left.  F\right\vert _{S\cap U}$ has constant rank, say $r_{1}$,
and the item $1$ follows. Using basic linear algebra and $\left(
\ref{dsum}\right)  $,
\[
\left(  \operatorname{Ker}\Pi_{\ast,s}\cap T_{s}S\right)  +\operatorname{Ker}%
F_{\ast,s}=\left(  \operatorname{Ker}\Pi_{\ast,s}+\operatorname{Ker}F_{\ast
,s}\right)  \cap T_{s}S=T_{s}S,
\]
then, from $\left(  \ref{inter2}\right)  $ and the fact that
$\operatorname{Ker}\left(  \left.  \Pi\right\vert _{S\cap U}\right)  _{\ast
,s}=\operatorname{Ker}\Pi_{\ast,s}\cap T_{s}S$, we finally have Eq. $\left(
\ref{ssum}\right)  $. This implies that the rank of $\left.  \Pi\right\vert
_{S\cap U}$ is constant, say $r_{2}$, and $r_{1}+r_{2}=\dim S$, i.e. the items
$2$ and $3$ follow.\ \ \ \ $\square$

\bigskip

Suppose now that we have an $X$-invariant foliation $\mathcal{S}$ of $M$ (i.e.
a foliation such that each one of its leaves is $X$-invariant). Denote by
$T\mathcal{S}$ the related distribution. From previous results, the next
proposition is straightforward.

\begin{proposition}
\label{resf}Let $\Sigma:N\times\Lambda\rightarrow M$ be a complete solution of
the $\Pi$-HJE for $\left(  M,X\right)  $, and let $\mathcal{S}$ be an
$X$-invariant foliation of $M$. If, for every momentum map $\left(
U,F\right)  $ related to $\Sigma$, we have that $\operatorname{Ker}F_{\ast
}\subset T\mathcal{S}$, or for every leaf $S$

\begin{enumerate}
\item $\operatorname{Ker}\Pi_{\ast}\cap TS$ (resp. $\operatorname{Ker}F_{\ast
}\cap TS$) is a constant rank distribution along $S$ (resp. along $U\cap S$),

\item $\dim\left(  \operatorname{Ker}\Pi_{\ast}\cap TS\right)  +\dim\left(
\operatorname{Ker}F_{\ast}\cap TS\right)  =\dim S$,
\end{enumerate}

then $\Sigma$ restricts to each leaf of $\mathcal{S}$.
\end{proposition}

We shall apply this result in the last section, in the context of Poisson
manifolds, where $\mathcal{S}$ will be the foliation given by the symplectic
leaves of $M$.

\subsubsection{Non-commutative integrability on Poisson manifolds}

\label{INC}

Fix a Poisson manifold $\left(  M,\Xi\right)  $ and a function $H:M\rightarrow
\mathbb{R}$. We give below a definition of the commutative and noncommutative
integrability notions (on Poisson manifolds) in terms that we find more
appropriate for this paper.\footnote{Our definition corresponds to the
so-called \emph{abstract (non)commutative integrability }of Ref. \cite{eva},
but with a \emph{momentum map}, as defined in \cite{rui}. Nevertheless, note
that we are not asking that some of the components of the momentum map are in
involution with all of them.}

\begin{definition}
\label{nci}The Hamiltonian system $\left(  M,X_{H}\right)  $ is
\textbf{non-commutative integrable (NCI)} by means of $F$ if a isotropic
fibration $F:M\rightarrow\Lambda$, such that $X_{H}\in\operatorname{Ker}%
F_{\ast}$ and $\Xi^{\sharp}\left[  \left(  \operatorname{Ker}F_{\ast}\right)
^{0}\right]  $ is integrable, can be exhibited. If in addition $F$ is
Lagrangian, $\left(  M,X_{H}\right)  $ is \textbf{commutative integrable (CI)}
by means of $F$.
\end{definition}

Under conditions above, the system $\left(  M,X_{H}\right)  $ can be
integrated by quadratures. This fact was proved in \cite{adler} for the
commutative case. A similar proof can be applied to the noncommutative
scenario, provided a number $\dim\left(  \operatorname{Ker}F_{\ast}\right)  $
of components of $F$ are in involution with all of them. In a forthcoming
paper, we shall further extend the proof to the general case.

Some authors include in the definition of NCI and CI systems one more
requirement: $F$ has compact and connected leaves. In such a case, beside
integrability by quadratures, \emph{action-angle} like coordinates can be
found for such systems (see Refs. \cite{eva} and \cite{rui}). We do not
analyze this case here.

\bigskip

Using the duality between complete solutions and first integrals (see Sections
\ref{tsfi} and \ref{fits}), a deep connection between our Hamilton-Jacobi
Theory and noncommutative integrability on Poisson manifolds can be easily
established (extending our results obtained at the end of Section \ref{scr},
in the restricted context of symplectic manifolds and commutative
integrability). We only need Proposition \ref{equi} and the next result, which
is an immediate consequence of the Remark \ref{setlfi2}. Consider a fibration
$\Pi:M\rightarrow N$ and a complete solution $\Sigma:N\times\Lambda\rightarrow
M$ of the $\Pi$-HJE for system $\left(  M,X_{H}\right)  $.

\begin{proposition}
\label{rotniso}$\left(  \Sigma^{\ast}\Xi\right)  ^{\sharp}\left(  \left(
TN\times0\right)  ^{0}\right)  $ is integrable if and only if so is
$\Xi^{\sharp}\left[  \left(  \operatorname{Ker}F_{\ast}\right)  ^{0}\right]  $
for all pair $\left(  U,F\right)  $.
\end{proposition}

Now, the mentioned connection for the general noncommutative case.

\begin{theorem}
\label{tnci}Consider a Poisson manifold $\left(  M,\Xi\right)  $, a function
$H:M\rightarrow\mathbb{R}$ and a fibration $\Pi:M\rightarrow N$.

\begin{enumerate}
\item Suppose that an isotropic complete solution $\Sigma:N\times
\Lambda\rightarrow M$ of the $\Pi$-HJE for $\left(  M,X_{H}\right)  $, with
\[
\left(  \Sigma^{\ast}\Xi\right)  ^{\sharp}\left(  \left(  TN\times0\right)
^{0}\right)
\]
integrable, can be exhibited. Then, for every $m\in M$, there exists an open
neighborhood $U$ of $m$ such that $\left(  U,\left.  X_{H}\right\vert
_{U}\right)  $ is NCI.

\item If $\left(  M,X_{H}\right)  $ is NCI by means of $F$, with $F$ and $\Pi$
transversal, then, for each $m\in M$, there exists an open neighborhood $U$ of
$m$ such that an isotropic complete solution $\Sigma$ of the $\left.
\Pi\right\vert _{U}$-HJE for $\left(  U,\left.  X_{H}\right\vert _{U}\right)
$, with
\[
\left(  \Sigma^{\ast}\Xi\right)  ^{\sharp}\left(  \left(  \Pi_{\ast}\left(
TU\right)  \times0\right)  ^{0}\right)
\]
integrable, can be exhibited.

\item If $\left(  M,X_{H}\right)  $ is NCI, then, for each $m\in M$, there
exists an open neighborhood $U$ of $m$ and a fibration $\pi:U\rightarrow
\pi\left(  U\right)  $ such that an isotropic complete solution $\Sigma$ of
the $\pi$-HJE for $\left(  U,\left.  X_{H}\right\vert _{U}\right)  $, with
\[
\left(  \Sigma^{\ast}\Xi\right)  ^{\sharp}\left(  \left(  \pi_{\ast}\left(
TU\right)  \times0\right)  ^{0}\right)
\]
integrable, can be exhibited.
\end{enumerate}
\end{theorem}

\emph{Proof.} $\left(  1\right)  $ Theorem \ref{pc} ensures the existence of
an open neighborhood $U$ of $m$ and that a fibration $F:U\rightarrow F\left(
U\right)  $, such that $\operatorname{Im}\left(  \left.  X_{H}\right\vert
_{U}\right)  \subset\operatorname{Ker}F_{\ast}$, can be constructed from
$\Sigma$. On the other hand, Propositions \ref{equi} and \ref{rotniso} ensure
that $F$ is isotropic with $\Xi^{\sharp}\left[  \left(  \operatorname{Ker}%
F_{\ast}\right)  ^{0}\right]  $ integrable. Then, $\left(  U,\left.
X_{H}\right\vert _{U}\right)  $ is NCI by means of $F$.

\smallskip

$\left(  2\right)  $ Since $\operatorname{Im}\left(  X_{H}\right)
\subset\operatorname{Ker}F_{\ast}$ and $F$ and $\Pi$ are transverse, Theorem
\ref{inthj} ensures the existence of $U$ and the possibility of constructing
$\Sigma$ from $F$; and since $F$ is isotropic and $\Xi^{\sharp}\left[  \left(
\operatorname{Ker}F_{\ast}\right)  ^{0}\right]  $ is integrable, Propositions
\ref{equi} and \ref{rotniso} ensure that the same is true for $\Sigma$ and
$\left(  \Sigma^{\ast}\Xi\right)  ^{\sharp}\left(  \left(  \Pi_{\ast}\left(
TU\right)  \times0\right)  ^{0}\right)  $.

\smallskip

$\left(  3\right)  $ We must proceed as in the previous item $(2)$, but
considering Theorem \ref{fint} instead of Theorem \ref{inthj}. \ \ \ $\square$

\bigskip

\begin{remark}
Kozlov has presented in \cite{koz} an alternative HJE (see Remark \ref{kozl})
and a connection of that equation with the non-commutative integrability in
symplectic manifolds. In contrast to our theory, the first integrals of the
system (by means of which the system is integrable) are not completely
determined by the solutions of the Kozlov's HJE.
\end{remark}

For completeness, let us express the previous result for the commutative case.

\begin{theorem}
\label{tci}Consider a Poisson manifold $\left(  M,\Xi\right)  $, a function
$H:M\rightarrow\mathbb{R}$ and a fibration $\Pi:M\rightarrow N$.

\begin{enumerate}
\item Suppose that a Lagrangian complete solution $\Sigma:N\times
\Lambda\rightarrow M$ of the $\Pi$-HJE for $\left(  M,X_{H}\right)  $ can be
exhibited. Then for every $m\in M$ there exists an open neighborhood $U$ of
$m$ such that $\left(  U,\left.  X_{H}\right\vert _{U}\right)  $ is CI.

\item If $\left(  M,X_{H}\right)  $ is CI by means of $F$, with $F$ and $\Pi$
transversal, then for each $m\in M$ there exists an open neighborhood $U$ of
$m$ such that a Lagrangian complete solution $\Sigma$ of the $\left.
\Pi\right\vert _{U}$-HJE for $\left(  U,\left.  X_{H}\right\vert _{U}\right)
$ can be exhibited.

\item If $\left(  M,X_{H}\right)  $ is CI, then for each $m\in M$ there exists
an open neighborhood $U$ of $m$ and a fibration $\pi:U\rightarrow\pi\left(
U\right)  $ such that a Lagrangian complete solution $\Sigma$ of the $\pi$-HJE
for $\left(  U,\left.  X_{H}\right\vert _{U}\right)  $ can be exhibited.
\end{enumerate}
\end{theorem}

\begin{remark}
If in Definition \ref{nci} we ask $F$ to have compact and connected leaves,
previous theorems are also true, provided we assume that $N$ is compact and
connected (see the Remark \ref{setlfi2}).
\end{remark}

\section{Weakly isotropic complete solutions and integrability by quadratures}

\label{S5}

\label{isoquad}

 Theorems \ref{tnci} and \ref{tci} say, among other things, that given a
Hamiltonian system $\left(  M,X_{H}\right)  $ on a Poisson manifold $\left(
M,\Xi\right)  $ fibered by $\Pi$, if an \textquotedblleft
appropriate\textquotedblright\ complete solution of the $\Pi$-HJE is known for
$\left(  M,X_{H}\right)  $, then $\left(  M,X_{H}\right)  $ is integrable by
quadratures. There, by \textquotedblleft appropriate\textquotedblright\ we
mean that $\Sigma$ is isotropic and the distribution $\left(  \Sigma^{\ast}%
\Xi\right)  ^{\sharp}\left(  \left(  TN\times0\right)  ^{0}\right)  $ is
integrable. We shall show in this section that only the isotropy condition is
needed to ensure integrability by quadratures.

\subsection{The case of symplectic manifolds}

\label{ghjm}

Consider a symplectic manifold $\left(  M,\omega\right)  $, a function
$H:M\rightarrow\mathbb{R}$, a fibration $\Pi:M\rightarrow N$ and a complete
solution $\Sigma:N\times\Lambda\rightarrow M$ of the $\Pi$-HJE for $\left(
M,X_{H}\right)  $. As in Section \ref{spr}, let us assume for simplicity that
$\omega=d\theta$, $N$ is simply connected (\emph{ipso facto} connected), and
$\Sigma$ is isotropic. Recall that, under these assumptions, $H\circ
\Sigma=h\circ p_{\Lambda}$ for a unique function $h:\Lambda\rightarrow
\mathbb{R}$, and there exists a function $W:N\times\Lambda\rightarrow
\mathbb{R}$ (that can be found by quadratures) satisfying%
\begin{equation}
\left\langle \left(  \Sigma^{\ast}\theta-dW\right)  \left(  n,\lambda\right)
,\left(  y,0\right)  \right\rangle =0 \label{tn}%
\end{equation}
for all $\left(  n,\lambda\right)  \in N\times\Lambda$ and $y\in T_{n}N$.
Since $d\Sigma^{\ast}\theta=d\left(  \Sigma^{\ast}\theta-dW\right)  $, using
the Cartan formula $i_{Y}\circ d+d\circ i_{Y}=L_{Y}$ for $Y:=X_{H}^{\Sigma}$,
Eq. $\left(  \ref{tn}\right)  $ and the fact $\operatorname{Im}\left(
X_{H}^{\Sigma}\right)  \subset TN\times0$, it follows that%
\[
i_{X_{H}^{\Sigma}}d\Sigma^{\ast}\theta=L_{X_{H}^{\Sigma}}\left(  \Sigma^{\ast
}\theta-dW\right)  .
\]
Combining the last equation and the identity $d\left(  H\circ\Sigma\right)
=p_{\Lambda}^{\ast}dh$, Eq. $\left(  \ref{hjes3}\right)  $ reads
\begin{equation}
p_{\Lambda}^{\ast}dh=L_{X_{H}^{\Sigma}}\left(  \Sigma^{\ast}\theta-dW\right)
. \label{qe}%
\end{equation}
Summing up, we have proved the following result.

\begin{proposition}
\label{mainglob}Consider a symplectic manifold $\left(  M,\omega\right)  $, a
function $H:M\rightarrow\mathbb{R}$, a fibration $\Pi:M\rightarrow N$ and a
complete solution $\Sigma$ of the $\Pi$-HJE for $\left(  M,X_{H}\right)  $. If
$\omega=d\theta$, $N$ is simply connected, and $\Sigma$ is isotropic, then
there exist functions $h:\Lambda\rightarrow\mathbb{R}$ and $W:N\times
\Lambda\rightarrow\mathbb{R}$ (which can be given by quadratures), such that
Eq. $\left(  \ref{qe}\right)  $ holds.
\end{proposition}

We are going to find a set of functional equations for the integral curves of
$X_{H}^{\Sigma}$. Let us consider the immersion $\varphi:N\times
\Lambda\rightarrow T^{\ast}\Lambda$ defined in (\ref{fitot}). Since the curves
of $X_{H}^{\Sigma}$ are of the form $\left(  \gamma\left(  t\right)
,\lambda\right)  $, with $\lambda\in\Lambda$ and $\gamma$ an integral curve of
$X_{H}^{\sigma_{\lambda}}$, let us fix some $\lambda$ and focus on the
immersion $\varphi_{\lambda}:N\rightarrow T_{\lambda}^{\ast}\Lambda$. Recall
that the later is given by%
\begin{equation}
\left\langle \varphi_{\lambda}\left(  n\right)  ,z\right\rangle =\left\langle
\left(  \Sigma^{\ast}\theta-dW\right)  \left(  n,\lambda\right)  ,\left(
0,z\right)  \right\rangle ,\ \ \ \forall z\in T_{\lambda}\Lambda. \label{fil}%
\end{equation}

\begin{proposition}
Given $\lambda\in\Lambda$ and an integral curve $\gamma:I\rightarrow N$ of
$X_{H}^{\sigma_{\lambda}}$, there exist unique covectors $\alpha_{\lambda
},\beta_{\lambda}\in T_{\lambda}^{\ast}\Lambda$ such that%
\[
\varphi_{\lambda}\left(  \gamma\left(  t\right)  \right)  =\alpha_{\lambda
}\,t+\beta_{\lambda},\ \ \ \forall t\in I.
\]

\end{proposition}

\emph{Proof.} Consider a vector field $z\in\mathfrak{X}\left(  \Lambda\right)
$ and the related vector field $Z=\left(  0,z\right)  \in\mathfrak{X}\left(
N\times\Lambda\right)  $. Since $\operatorname{Im}X_{H}^{\Sigma}\subset
TN\times0$, it is easy to see that $\operatorname{Im}\left[  Z,X_{H}^{\Sigma
}\right]  \subset TN\times0$. Using the identity%
\[
i_{\left[  Z,X_{H}^{\Sigma}\right]  }=i_{Z}\circ L_{X_{H}^{\Sigma}}%
-L_{X_{H}^{\Sigma}}\circ i_{Z},
\]
and the fact that $i_{Y}\left(  \Sigma^{\ast}\theta-dW\right)  =0$ for all $Y$
such that $\operatorname{Im}Y\subset TN\times0$ [see Eq. $\left(
\ref{tn}\right)  $], it follows that%
\[
i_{Z}\left[  L_{X_{H}^{\Sigma}}\left(  \Sigma^{\ast}\theta-dW\right)  \right]
=L_{X_{H}^{\Sigma}}\left[  i_{Z}\left(  \Sigma^{\ast}\theta-dW\right)
\right]  .
\]
Then, applying $i_{Z}$ on both sides of $\left(  \ref{qe}\right)  $, we obtain%
\[
\left\langle dh,z\right\rangle =L_{X_{H}^{\Sigma}}\left\langle \Sigma^{\ast
}\theta-dW,\left(  0,z\right)  \right\rangle .
\]
Evaluating on the curve $\left(  \gamma\left(  t\right)  ,\lambda\right)  $,
we have%
\begin{align*}
\left\langle dh\left(  \lambda\right)  ,z\left(  \lambda\right)
\right\rangle  &  =\frac{d}{dt}\left\langle \left(  \Sigma^{\ast}%
\theta-dW\right)  \left(  \gamma\left(  t\right)  ,\lambda\right)  ,\left(
0,z\left(  \lambda\right)  \right)  \right\rangle \\
&  =\frac{d}{dt}\left\langle \varphi_{\lambda}\left(  \gamma\left(  t\right)
\right)  ,z\left(  \lambda\right)  \right\rangle .
\end{align*}
Assume for simplicity that $0\in I$. Then, integrating the equation above
between $0$ and $t$,%
\[
\left\langle \varphi_{\lambda}\left(  \gamma\left(  t\right)  \right)
,z\left(  \lambda\right)  \right\rangle =\left\langle dh\left(  \lambda
\right)  ,z\left(  \lambda\right)  \right\rangle \,t+\left\langle
\varphi_{\lambda}\left(  \gamma\left(  0\right)  \right)  ,z\left(
\lambda\right)  \right\rangle ,
\]
i.e.%
\begin{equation}
\varphi_{\lambda}\left(  \gamma\left(  t\right)  \right)  =dh\left(
\lambda\right)  \,t+\varphi_{\lambda}\left(  \gamma\left(  0\right)  \right)
, \label{fila}%
\end{equation}
which proves the proposition.$\ \ \ \ \square$

\bigskip

Equation $\left(  \ref{fila}\right)  $ defines the functional equations for
$\gamma$ that we mentioned above. The Inverse Function Theorem for the
immersion $\varphi_{\lambda}$ ensures that such functional equations can be
solved for $\gamma\left(  t\right)  $. On the other hand, since each
$\varphi_{\lambda}$ is constructed from $\Sigma$ by quadratures (i.e. through
the function $W$), we can conclude that the dynamical system $\left(
M,X_{H}\right)  $ can be solved by quadratures. Note that this is also true in
spite of $\omega$ is not exact or $N$ is not simply connected. In fact, as we
argue at the end of Section \ref{spr}, such conditions are always true at a
local level, and consequently we can write a functional equation like $\left(
\ref{fila}\right)  $ around every point of $M$. Concluding,

\begin{theorem}
\label{ibq}Consider a symplectic manifold $\left(  M,\omega\right)  $, a
function $H:M\rightarrow\mathbb{R}$, a fibration $\Pi:M\rightarrow N$ and a
complete solution $\Sigma$ of the $\Pi$-HJE for $\left(  M,X_{H}\right)  $. If
$\Sigma$ is isotropic, then $\left(  M,X_{H}\right)  $ can be solved by quadratures.
\end{theorem}

In terms of first integrals, using Theorem \ref{fint} and Proposition
\ref{equi} (and recalling the Remark \ref{adap}), we have the next result.

\begin{theorem}
\label{ibq2}Consider a symplectic manifold $\left(  M,\omega\right)  $, a
function $H:M\rightarrow\mathbb{R}$ and a fibration $F:M\rightarrow\Lambda$
such that $\operatorname{Im}X_{H}\subset\operatorname{Ker}F_{\ast}$. If $F$ is
isotropic, then $\left(  M,X_{H}\right)  $ can be solved by quadratures.
\end{theorem}

Theorem above gives an affirmative answer to the question formulated in
\cite{zung} (page 6), if by \textquotedblleft integrable\textquotedblright\ we
mean some kind of \emph{exact solvability}.

\subsection{An illustrative example}

On the symplectic space $T^{\ast}\mathbb{R}^{4}\cong\mathbb{R}^{8}$ with its
canonical structure $\omega$, consider the Hamiltonian function $H:\mathbb{R}%
^{8}\rightarrow\mathbb{R}$ given by
\[
H(q,p)=(p_{2}-q^{2})(p_{4}-q^{4}),
\]
where $(q,p)=(q^{1},q^{2},q^{3},q^{4},p_{1},p_{2},p_{3},p_{4})$ are the
standard Darboux coordinates on $T^{\ast}\mathbb{R}^{4}$. Fixing
$N:=\mathbb{R}^{3}$,
\[
\Pi:\mathbb{R}^{8}\rightarrow\mathbb{R}^{3}:(q,p)\longmapsto(p_{2},p_{4}%
,q^{3})
\]
and $\Lambda:=\mathbb{R}^{5}$, it can be shown that the diffeomorphism
$\Sigma:\mathbb{R}^{3}\times\mathbb{R}^{5}\rightarrow\mathbb{R}^{8}$ given by%

\[
\Sigma(n,\lambda)=\left(  \lambda_{1},\lambda_{2}+n_{1},n_{3},\lambda
_{3}+n_{2},\lambda_{4},n_{1},\lambda_{5}-\lambda_{4}n_{3},n_{2}\right)  ,
\]
where $(n,\lambda):=(n_{1},n_{2},n_{3},\lambda_{1},\lambda_{2},\lambda
_{3},\lambda_{4},\lambda_{5})$, is a complete solution of the $\Pi$-HJE for
the Hamiltonian system $\left(  \mathbb{R}^{8},X_{H}\right)  $. The momentum
map related to $\Sigma$ is
\[
F:\mathbb{R}^{8}\rightarrow\mathbb{R}^{5}:(q,p)\longmapsto(q^{1},q^{2}%
-p_{2},q^{4}-p_{4},p_{1},q^{3}p_{1}+p_{3}),
\]
which is an isotropic submersion. In fact,
\[
\ker F_{\ast}=\left\langle \frac{\partial}{\partial q^{2}}+\frac{\partial
}{\partial p_{2}},\frac{\partial}{\partial q^{4}}+\frac{\partial}{\partial
p_{4}},\frac{\partial}{\partial q^{3}}-p_{1}\frac{\partial}{\partial p_{3}%
}\right\rangle \subseteq(\ker F_{\ast})^{\omega}.
\]
This implies that $\Sigma$ is isotropic too.

\begin{remark}
Since the vector fields $\displaystyle\frac{\partial}{\partial q^{3}}%
-p_{1}\displaystyle\frac{\partial}{\partial p_{3}}$ and $\displaystyle\frac
{\partial}{\partial p_{1}}$ belong to $(\ker F_{\ast})^{\omega}$, but
\end{remark}

\[
\left[  \frac{\partial}{\partial q^{3}}-p_{1}\frac{\partial}{\partial p_{3}%
},\frac{\partial}{\partial p_{1}}\right]  =\frac{\partial}{\partial p_{3}%
}\notin(\ker F_{\ast})^{\omega},
\]
we have that $(\ker F_{\ast})^{\omega}$ is \textbf{not} integrable.

Regarding the functions $h:\mathbb{R}^{5}\rightarrow\mathbb{R}$ and
$W:\mathbb{R}^{3}\times\mathbb{R}^{5}\rightarrow\mathbb{R}$, we have from the
equation $H\circ\Sigma=h\circ p_{\Lambda}$ that $h\left(  \lambda\right)
=\lambda_{2}\lambda_{3}$, and it is easy to see that
\[
W\left(  n,\lambda\right)  :=\frac{n_{1}^{2}+n_{2}^{2}}{2}+\left(  \lambda
_{3}n_{2}+\lambda_{2}n_{1}-\frac{1}{2}\lambda_{4}n_{3}^{2}\right)
\]
satisfies $\left(  \ref{tn}\right)  $. Then, the immersion $\varphi_{\lambda
}:\mathbb{R}^{3}\rightarrow T_{\lambda}^{\ast}\mathbb{R}^{5}\cong%
\mathbb{R}^{5}$ defined in $\left(  \ref{fil}\right)  $ is just
\begin{equation}
\varphi_{\lambda}(n)=\left(  0,-n_{1},-n_{2},\lambda_{1}+\frac{1}{2}n_{3}%
^{2},n_{3}\right)  , \label{fln}%
\end{equation}
and, using $\left(  \ref{fila}\right)  $ and the fact that
\[
dh\left(  \lambda\right)  =\left(  0,\lambda_{3},\lambda_{2},0,0\right)  ,
\]
it transforms the integral curve $\gamma\left(  t\right)  =\left(
n_{1}\left(  t\right)  ,n_{2}\left(  t\right)  ,n_{3}\left(  t\right)
\right)  $ of $X_{H}^{\sigma_{\lambda}}$ with initial condition $\left(
n_{1}^{0},n_{3}^{0},n_{3}^{0}\right)  $ into the expression
\begin{equation}
\varphi_{\lambda}\left(  \gamma\left(  t\right)  \right)  =\left(
0,\lambda_{3},\lambda_{2},0,0\right)  \ t+\left(  0,-n_{1}^{0},-n_{2}%
^{0},\lambda_{1}+\frac{1}{2}\left(  n_{3}^{0}\right)  ^{2},n_{3}^{0}\right)  .
\label{flg}%
\end{equation}
This implies that, combining $\left(  \ref{fln}\right)  $ and $\left(
\ref{flg}\right)  $, 
\[
n_{1}\left(  t\right)  =n_{1}^{0}-\lambda_{3}\ t,\ \ \ n_{2}\left(  t\right)
=n_{2}^{0}-\lambda_{2}\ t\ \ \ \text{and\ \ \ }n_{3}\left(  t\right)
=n_{3}^{0}.
\]
Consequently, using Eq. $\left(  \ref{Gg}\right)  $, any curve $\Gamma\left(
t\right)  =\left(  q\left(  t\right)  ,p\left(  t\right)  \right)  $ of
$X_{H}$ is given by%
\[
\left(  q\left(  t\right)  ,p\left(  t\right)  \right)  =\left(  \lambda
_{1},\lambda_{2}+n_{1}^{0}-\lambda_{3}\ t,n_{3}^{0},\lambda_{3}+n_{2}%
^{0}-\lambda_{2}\ t,\lambda_{4},n_{1}^{0}-\lambda_{3}\ t,\lambda_{5}%
-\lambda_{4}n_{3}^{0},n_{2}^{0}-\lambda_{2}\ t\right)  ,
\]
for a certain $8$-uple $\left(  n_{1}^{0},n_{3}^{0},n_{3}^{0},\lambda
_{1},\lambda_{2},\lambda_{3},\lambda_{4},\lambda_{5}\right)  \in\mathbb{R}%
^{8}$.

\bigskip

Of course, the equations of motion related to $H$ can be integrated without
any change of coordinates. We did it in this way just to illustrate the
procedure developed in the previous sections.

\subsection{The case of Poisson manifolds}

\label{poisson}

Let $\left(  M,\Xi\right)  $ be a Poisson manifold. Recall that
$\operatorname{Im}\Xi^{\sharp}$ is an integrable distribution and its maximal
connected integral manifolds $S$ are symplectic manifolds $\left(
S,\omega\right)  $: the \emph{symplectic leaves}. Denote by $\mathcal{S}$ the
related (symplectic) foliation, i.e. $\operatorname{Im}\Xi^{\sharp
}=T\mathcal{S}$. Let us fix a symplectic leaf $\left(  S,\omega\right)  $ and
indicate by $i:S\rightarrow M$ the inclusion map. For each $m\in S$,
\[
T_{m}S=\operatorname{Im}\Xi_{m}^{\sharp}=\left(  \operatorname{Ker}\Xi
_{m}^{\sharp}\right)  ^{0}%
\]
and $\omega\left(  \Xi^{\sharp}\left(  \alpha\right)  ,\Xi^{\sharp}\left(
\beta\right)  \right)  =\Xi\left(  \alpha,\beta\right)  $, $\forall
\alpha,\beta\in T_{m}^{\ast}M$. Consequently,
\begin{equation}
\Xi^{\sharp}\left(  \alpha\right)  =i_{\ast,m}\left[  \omega^{\sharp}\left(
i_{m}^{\ast}\left(  \alpha\right)  \right)  \right]  ,\ \ \ \forall\alpha\in
T_{m}^{\ast}M. \label{pw}%
\end{equation}
Fix also a function $H:M\rightarrow\mathbb{R}$ and consider its Hamiltonian
vector field $X_{H}$ w.r.t. $\Xi$. For all $m\in S$,%
\[
X_{H}\left(  m\right)  =\Xi^{\sharp}\left(  dH\left(  m\right)  \right)
=i_{\ast,m}\left[  \omega^{\sharp}\left(  i_{m}^{\ast}\left(  dH\right)
\right)  \right]  =i_{\ast,m}\left[  \left(  \omega^{\sharp}\circ i^{\ast
}dH\right)  \left(  m\right)  \right]  ,
\]
so, defining $H_{S}:=H\circ i:S\rightarrow\mathbb{R}$, we have that
$X_{H}\circ i=i_{\ast}\circ X_{H_{S}}$, where $X_{H_{S}}=\omega^{\sharp}\circ
dH_{S}\in\mathfrak{X}\left(  S\right)  $ (the Hamiltonian vector field w.r.t.
$\omega$ and with Hamiltonian function $H_{S}$). This says, in particular,
that $S$ is $X_{H}$-invariant and, moreover, the trajectories of $X_{H}$
intersecting $S$ coincide with the trajectories of $X_{H_{S}}$.

\begin{remark}
\label{sm}Thus, if we known the trajectories of each Hamiltonian system
$\left(  S,X_{H_{S}}\right)  $, for every symplectic leaf $S$, then we shall
know the trajectories of $\left(  M,X_{H}\right)  $.
\end{remark}

From now on, we shall assume that we know (each leaf $S$ of) the symplectic
foliation $\mathcal{S}$. Let $\Pi:M\rightarrow N$ be a fibration and
$\Sigma:N\times\Lambda\rightarrow M$ a complete solution of the $\Pi$-HJE for
$\left(  M,X_{H}\right)  $. According to the discussion we made in Section
\ref{rest}, in order to take advantage of both the complete solution $\Sigma$
and the $X_{H}$-invariant foliation $\mathcal{S}$, we need that $\Sigma$
restricts to each leaf $S$ of $\mathcal{S}$ (see Definition \ref{resr}). In
addition, using Theorem \ref{ibq}, if each local restriction of $\Sigma$ to a
given leaf $S$ is isotropic w.r.t. the symplectic structure $\omega$ on $S$,
the Hamiltonian system $\left(  S,X_{H_{S}}\right)  $ can be integrated by
quadratures. And from above remark, if this is true for all $S$, then $\left(
M,X_{H}\right)  $ is integrable by quadratures. Let us study sufficient
conditions for this to happen.

\begin{proposition}
\label{prop1}Consider a Poisson manifold $\left(  M,\Xi\right)  $, a function
$H:M\rightarrow\mathbb{R}$, a fibration $\Pi:M\rightarrow N$ and a complete
solution $\Sigma$ of the $\Pi$-HJE for $\left(  M,X_{H}\right)  $. If $\Sigma$
is isotropic, then $\Sigma$ restricts to each leaf of its symplectic foliation
$\mathcal{S}$.
\end{proposition}

\emph{Proof.} If $\Sigma$ is isotropic (recall Definition \ref{isodef}), all
the local momentum maps $\left(  U,F\right)  $ related to $\Sigma$ are also
isotropic (see Proposition \ref{equi}), and consequently $\operatorname{Ker}%
F_{\ast}\subset\operatorname{Im}\Xi^{\sharp}=T\mathcal{S}$. So, the hypothesis
of the Proposition \ref{resf} hold.\ \ $\square$

\bigskip

\begin{proposition}
\label{prop2}Under the same conditions of previous proposition, assume that
$\Sigma$ is weakly isotropic and restricts to a given symplectic leaf $S$.
Then, for every $s\in S$ there exists an isotropic local restriction of
$\Sigma$ to $S$ around $s$ (w.r.t. the symplectic structure $\omega$ on $S$).
\end{proposition}

\emph{Proof.} From Proposition \ref{equi}, $\Sigma$ is weakly isotropic if and
only if every local momentum map $\left(  U,F\right)  $ satisfies
\begin{equation}
\operatorname{Ker}F_{\ast}\cap\operatorname{Im}\Xi^{\sharp}\subset\Xi^{\sharp
}\left[  \left(  \operatorname{Ker}F_{\ast}\right)  ^{0}\right]  . \label{ft}%
\end{equation}
Since $\Sigma$ restricts to $S$, given $s\in S\cap U$, consider an open
neighborhood $R\subset S$ of $s$,$\ $a submanifold $\Lambda_{1}\subset\Lambda$
and a local restriction $\left.  \Sigma\right\vert _{\Pi\left(  R\right)
\times\Lambda_{1}}:\Pi\left(  R\right)  \times\Lambda_{1}\rightarrow R$ of
$\Sigma$ to $S$ around $s$. Shrink $R$ and $\Lambda_{1}$, if it is necessary,
to ensure that $R\subset S\cap U$ and $\left.  \Sigma\right\vert _{\Pi\left(
R\right)  \times\Lambda_{1}}$ is a global diffeomorphism. It is clear that
\[
\left.  F\right\vert _{R}=\left.  p_{\Lambda}\right\vert _{\Pi\left(
R\right)  \times\Lambda_{1}}\circ\left(  \left.  \Sigma\right\vert
_{\Pi\left(  R\right)  \times\Lambda_{1}}\right)  ^{-1}%
,\text{\ \ \ \ \ \ \ \ }F\left(  R\right)  =\Lambda_{1}%
\]
and $\left.  F\right\vert _{R}:R\rightarrow\Lambda_{1}$ is the global momentum
map for the system $\left(  R,\left.  X_{H}\right\vert _{R}\right)  $ related
to $\left.  \Sigma\right\vert _{\Pi\left(  R\right)  \times\Lambda_{1}}$.
Regarding Proposition \ref{equi} again (but for global momentum maps), if we
show that $\left.  F\right\vert _{R}$ is isotropic, then so is $\left.
\Sigma\right\vert _{\Pi\left(  R\right)  \times\Lambda_{1}}$. Let
$i:R\rightarrow U$ be the inclusion and write $\left.  F\right\vert
_{R}=F\circ i$. Since%
\begin{align*}
i_{\ast}\left(  \operatorname{Ker}\left(  \left.  F\right\vert _{R}\right)
_{\ast,r}\right)   &  =i_{\ast}\left(  \operatorname{Ker}\left(  F\circ
i\right)  _{\ast,r}\right)  =\operatorname{Ker}F_{\ast,r}\cap i_{\ast}\left(
T_{r}R\right)  =\operatorname{Ker}F_{\ast,r}\cap\operatorname{Im}\Xi
_{r}^{\sharp}\\
&  \subset\Xi^{\sharp}\left[  \left(  \operatorname{Ker}F_{\ast,r}\right)
^{0}\right]  ,
\end{align*}
for all $r\in R$ [see Eq. $\left(  \ref{ft}\right)  $], and%
\begin{align*}
\Xi^{\sharp}\left[  \left(  \operatorname{Ker}F_{\ast}\right)  ^{0}\right]
&  =\Xi^{\sharp}\left[  \left(  \operatorname{Ker}F_{\ast}\right)
^{0}+\operatorname{Ker}\Xi^{\sharp}\right]  =\Xi^{\sharp}\left[  \left(
\operatorname{Ker}F_{\ast}\cap\left(  \operatorname{Ker}\Xi^{\sharp}\right)
^{0}\right)  ^{0}\right] \\
&  =\Xi^{\sharp}\left[  \left(  \operatorname{Ker}F_{\ast}\cap
\operatorname{Im}\Xi^{\sharp}\right)  ^{0}\right]  ,
\end{align*}
then%
\begin{align}
i_{\ast}\left(  \operatorname{Ker}\left(  \left.  F\right\vert _{R}\right)
_{\ast,r}\right)   &  \subset\Xi^{\sharp}\left[  \left(  \operatorname{Ker}%
F_{\ast,r}\right)  ^{0}\right]  =\Xi^{\sharp}\left[  \left(
\operatorname{Ker}F_{\ast,r}\cap\operatorname{Im}\Xi_{r}^{\sharp}\right)
^{0}\right] \nonumber\\
&  =\Xi^{\sharp}\left[  \left(  i_{\ast}\left(  \operatorname{Ker}\left(
\left.  F\right\vert _{R}\right)  _{\ast,r}\right)  \right)  ^{0}\right]
,\ \ \ \ \ \forall r\in R. \label{ft1}%
\end{align}
On the other hand, from Eq. $\left(  \ref{pw}\right)  $ we have
that\footnote{Given two sets $A$ and $B$ and a function $G:A\rightarrow B$,
recall that $G\left(  G^{-1}\left(  S\right)  \right)  \subset S$ for every
subset $S\subset B$. Also, if $A$ and $B$ are vector spaces, then $\left(
G\left(  S\right)  \right)  ^{0}=\left(  G^{\ast}\right)  ^{-1}\left(
S^{0}\right)  $, for every subspace $S\subset A$.}%
\begin{align}
\Xi^{\sharp}\left[  \left(  i_{\ast,r}\left(  \operatorname{Ker}\left(
\left.  F\right\vert _{R}\right)  _{\ast,r}\right)  \right)  ^{0}\right]   &
=\Xi^{\sharp}\left[  \left(  i_{r}^{\ast}\right)  ^{-1}\left(  \left(
\operatorname{Ker}\left(  \left.  F\right\vert _{R}\right)  _{\ast,r}\right)
^{0}\right)  \right] \nonumber\\
&  =i_{\ast}\left[  \omega^{\sharp}\left[  i_{r}^{\ast}\left(  \left(
i_{r}^{\ast}\right)  ^{-1}\left(  \left(  \operatorname{Ker}\left(  \left.
F\right\vert _{R}\right)  _{\ast,r}\right)  ^{0}\right)  \right)  \right]
\right] \nonumber\\
&  \subset i_{\ast}\left[  \omega^{\sharp}\left[  \left(  \operatorname{Ker}%
\left(  \left.  F\right\vert _{R}\right)  _{\ast,r}\right)  ^{0}\right]
\right]  . \label{ft2}%
\end{align}
As a consequence, combining $\left(  \ref{ft1}\right)  $ and $\left(
\ref{ft2}\right)  $,%
\[
i_{\ast}\left(  \operatorname{Ker}\left(  \left.  F\right\vert _{R}\right)
_{\ast}\right)  \subset i_{\ast}\left[  \omega^{\sharp}\left(  \left(
\operatorname{Ker}\left(  \left.  F\right\vert _{R}\right)  _{\ast}\right)
^{0}\right)  \right]  ,
\]
and from the injectivity of $i_{\ast}$ we have that $\left.  F\right\vert
_{R}$ is isotropic w.r.t. $\omega$. \ \ \ $\square$

\bigskip

Now, we are ready to the main result of this section.

\begin{theorem}
\label{ibqp}Consider a Poisson manifold $\left(  M,\Xi\right)  $, its
symplectic foliation $\mathcal{S}$, a function $H:M\rightarrow\mathbb{R}$, a
fibration $\Pi:M\rightarrow N$ and a complete solution $\Sigma$ of the $\Pi
$-HJE for $\left(  M,X_{H}\right)  $. If $\Sigma$ is weakly isotropic and for
every leaf $S$ of $\mathcal{S}$ we have that

\begin{enumerate}
\item $\operatorname{Ker}\Pi_{\ast}\cap TS$ (resp. $\operatorname{Ker}F_{\ast
}\cap TS$) is a constant rank distribution along $S$ (resp. $S\cap U$),

\item $\dim\left(  \operatorname{Ker}\Pi_{\ast}\cap TS\right)  +\dim\left(
\operatorname{Ker}F_{\ast}\cap TS\right)  =\dim S$,
\end{enumerate}

then $\left(  M,X_{H}\right)  $ can be solved by quadratures.
\end{theorem}

\emph{Proof.} According to the Propositions \ref{resf} and \ref{prop2}, the
conditions above ensure that there exists an isotropic local restriction of
$\Sigma$ to every $S$ and around every $s\in S$. Then, by Theorem
\ref{ibq}, each Hamiltonian system $\left(  S,X_{H_{S}}\right)  $ can be
integrated by quadratures, and from Remark \ref{sm} the same is true for
$\left(  M,X_{H}\right)  $. \ \ \ $\square$

\bigskip

Combining the Propositions \ref{prop1} and \ref{prop2} and the last theorem,
we have the following corollary.

\begin{corollary}
Consider a Poisson manifold $\left(  M,\Xi\right)  $, its symplectic foliation
$\mathcal{S}$, a function $H:M\rightarrow\mathbb{R}$, a fibration
$\Pi:M\rightarrow N$ and a complete solution $\Sigma$ of the $\Pi$-HJE for
$\left(  M,X_{H}\right)  $. If $\Sigma$ is isotropic, then $\left(
M,X_{H}\right)  $ can be solved by quadratures.
\end{corollary}

\bigskip

In terms of first integrals, we have the next result and its immediate corollary.

\begin{theorem}
Consider a Poisson manifold $\left(  M,\Xi\right)  $, its symplectic foliation
$\mathcal{S}$, a function $H:M\rightarrow\mathbb{R}$ and a fibration
$F:M\rightarrow\Lambda$ such that $\operatorname{Im}X_{H}\subset
\operatorname{Ker}F_{\ast}$. If $F$ is weakly isotropic and for every leaf $S$
of $\mathcal{S}$ we have that $\operatorname{Ker}F_{\ast}\cap TS$ is a
constant rank distribution along $S$, then $\left(  M,X_{H}\right)  $ can be
solved by quadratures.
\end{theorem}

\emph{Proof.} Fix a leaf $S$ and a point $s\in S$. Since $\operatorname{Ker}%
F_{\ast}\cap TS$ is a constant rank distribution along $S$, then $\left.
F\right\vert _{S}$ has constant rank. As a consequence, as we saw in the proof
of Proposition \ref{pq}, there exists an open neighborhood $V\subset S$ of $s$
such that $F\left(  V\right)  \subset\Lambda$ is a submanifold, $\left.
F\right\vert _{V}:V\rightarrow F\left(  V\right)  $ is a fibration and
\[
\operatorname{Im}\left(  \left.  X_{H_{S}}\right\vert _{V}\right)
=\operatorname{Im}\left(  \left.  X_{H}\right\vert _{V}\right)  \subset
\operatorname{Ker}\left(  \left.  F\right\vert _{V}\right)  _{\ast}.
\]
On the other hand, since $F$ is weakly isotropic, we can prove as in
Proposition \ref{prop2} [see Eqs. $\left(  \ref{ft1}\right)  $ and $\left(
\ref{ft2}\right)  $] that $\left.  F\right\vert _{V}$ is isotropic w.r.t. the
symplectic structure $\omega$ of $S$. Then, the present theorem follows from
Theorem \ref{ibq2} and the Remark \ref{sm}.\ \ \ $\square$

\bigskip

\begin{corollary}
Consider a Poisson manifold $\left(  M,\Xi\right)  $, its symplectic foliation
$\mathcal{S}$, a function $H:M\rightarrow\mathbb{R}$ and a fibration
$F:M\rightarrow\Lambda$ such that $\operatorname{Im}X_{H}\subset
\operatorname{Ker}F_{\ast}$. If $F$ is isotropic, then $\left(  M,X_{H}%
\right)  $ can be solved by quadratures.
\end{corollary}

Concluding, given a Hamiltonian system on a Poisson manifold, in order to
ensure integrability by quadratures it is enough to have a weakly isotropic
set of first integrals (and an additional regularity condition), provided the
symplectic foliation of the manifold is known. Note that we are not asking for
the Poisson structure to have constant rank.

\section{Final comments and future work}

We have defined a Hamilton-Jacobi equation (HJE) for general dynamical systems
on fibered phase spaces. We have focused on the consequences of having a
partial or a complete solution of such an equation, for a given dynamical
system, and, in particular, on the possibility of integrating by quadratures
the equations of motion of the system in question. We have shown that the
resulting Hamilton-Jacobi Theory is an extension of those developed for
Hamiltonian systems on symplectic, Poisson and almost-Poisson manifolds. We
have also shown, in the context of general dynamical systems, that there
exists a deep connection between first integrals and the complete solutions of
our HJE. This enabled us, in the context of Poisson manifolds, to identify
properties of the complete solutions that give rise to non-commutative
integrability. Also, a new method to integrate by quadratures a Hamiltonian
system on a Poisson manifold have been developed. It is worth mentioned that,
such a method, can be applied under weaker assumptions than those considered
in the notion of non-commutative integrability.

This paper contains our first steps in the study of this new Hamilton-Jacobi
Theory. In future works, we are planning to address the following problems.

\textbf{Action-angle coordinates. }We want to study the relationship between
the generalized action-angle coordinates, which appears in a non-commutative
integrable system on a Poisson manifold, and the immersion $\varphi
:N\times\Lambda\rightarrow T^{\ast}\Lambda$ defined in Eq. $\left(
\ref{fitot}\right)  $ for isotropic complete solutions.

\textbf{Hamiltonian systems with external forces and constrained systems.}
Given a manifold $Q$, a Hamiltonian system with external forces on $Q$ is a
dynamical system $\left(  T^{\ast}Q,X\right)  $ with $X=X_{H}+Y$ for some
vertical vector field $Y\in\mathfrak{X}\left(  T^{\ast}Q\right)  $ (i.e.
$\operatorname{Im}Y\subset\operatorname{Ker}\left(  \pi_{Q}\right)  _{\ast}$):
\emph{the external force}. Such a force can represent a constraint force. So,
the nonholonomic and generalized nonholonomic systems are particular examples.
The complete version of the $\Pi$-HJE for $\left(  T^{\ast}Q,X\right)  $ [see
the second part of Eq. $\left(  \ref{Srel}\right)  $], and for any fibration
$\Pi:T^{\ast}Q\rightarrow N$, can be written
\[
i_{X^{\Sigma}}\Sigma^{\ast}\omega=\Sigma^{\ast}\left(  dH+\beta\right)
,\ \ \ \text{with\ \ }\beta:=\omega^{\flat}\left(  Y\right)  .
\]
If we ask, for instance, that $\Sigma$ is isotropic and $\left(  \ast\right)
$ $\operatorname{Ker}\left(  \Sigma^{\ast}\beta\right)  ,\operatorname{Ker}%
\left(  \Sigma^{\ast}d\beta\right)  \subset TN\times0$, it can be shown that a
procedure analogous to that realized in Section \ref{ghjm} can be done for
$\left(  T^{\ast}Q,X\right)  $. Consequently, the system $\left(  T^{\ast
}Q,X\right)  $ can be integrated by quadratures. We are studying weaker
conditions than $\left(  \ast\right)  $, and trying to apply them to
particular examples of nonholonomic systems.

\textbf{Time-dependent systems, time-dependent HJE, co-symplectic and contact
manifolds.} Cosymplectic and contact manifolds are the counterpart of the
symplectic manifolds in odd dimension. The cosymplectic geometry plays an
important role in the geometric description of the time-dependent Mechanics.
On the other hand, the Thermodynamics admits a geometric formulation in terms
of a certain contact structure on the phase space of the system. We shall
study the relationship between the integration of the equations which
determine the dynamics in both cases (cosymplectic and contact) and the
Hamilton-Jacobi Theory. In the case of the cosymplectic structures, since they
are Poisson structures, part of the theory developed in this paper can be considered.

\textbf{Symmetries and reconstruction. }Consider a Lie group $G$ and a
$G$-invariant dynamical\footnote{This means that there exists an action
$\rho:G\times M\rightarrow M$ such that $\left(  \rho_{g}\right)  _{\ast}\circ
X=X\circ\rho_{g}$, for all $g\in G$.} system $\left(  M,X\right)  $ such that
$N:=M/G$ is a manifold and the canonical projection $\Pi:M\rightarrow N$ is a
fibration. Let $Y\in\mathfrak{X}\left(  N\right)  $ be the reduced field, i.e.
the unique vector field on $N$ such that $\Pi_{\ast}\circ X=Y\circ\Pi$. It can
be shown that, for every complete solution $\Sigma:N\times\Lambda\rightarrow
M$ of the $\Pi$-HJE for $\left(  M,X\right)  $, we have that $X^{\Sigma
}\left(  n,\lambda\right)  =\left(  Y\left(  n\right)  ,0\right)  $ for all
$n,\lambda$. Then, each integral curve $\Gamma:I\rightarrow M$ of $X$ can be
obtained from an integral curve $\gamma:I\rightarrow N$ of the reduced field
$Y$ by the formula $\Gamma\left(  t\right)  =\Sigma\left(  \gamma\left(
t\right)  ,\lambda\right)  $, for some $\lambda\in\Lambda$. This means that
the $\Pi$-HJE gives rise to a \emph{reconstruction process}. We want to
further study this fact.

\section{Acknowledgements}

This work was partially supported by Ministerio de Econom{\'\i}a y Competitividad  (Spain) grants MTM2012-34478 (E.
Padr\'{o}n) and the European Community IRSES-project GEOMECH-246981 (S. Grillo
and E. Padr\'{o}n). S. Grillo thanks CONICET for its financial support.


\bigskip

\begin{thebibliography}{99}                                                                                               %


\bibitem {am}R. Abraham, J.E. Marsden, \textit{Foundation of Mechanics}, New
York, Benjaming Cummings (1985).

\bibitem {adler}M. Adler, P. van Moerbeke, and P. Vanhaecke. \textit{Algebraic
integrability, Painlev e geometry and Lie algebras}, volume 47 of Results in
Mathematics and Related Areas. 3rd Series. A Series of Modern Surveys in
Mathematics. Springer-Verlag, Berlin, 2004.

\bibitem {ar2}V.I. Arnold, Ordinary differential equations,\text{ MIT Press, }
(1973) Translated from the Russian by Richard A. Silverman.

\bibitem {ar}V.I. Arnold, \textit{Mathematical Models in Classical Mechanics},
Berlin, Springer-Verlag (1978).

\bibitem {bmmp}P. Balseiro, J.C. Marrero, D. Mart{\'{\i}}n de Diego, E.
Padr\'{o}n, \textit{A unified framework for Mechanics, Hamilton-Jacobi
equation and applications}, {Nonlinearity}, \textbf{23}, 8 (2010), 1887--1918.

\bibitem {bfs}L. Bates, S. Fasso, N. Sansonetto, \textit{The Hamilton-Jacobi
equation, integrability, and nonholonomic systems}, Journal of Geometric
Mechanics \textbf{6}, 4 (2014), 441-449.

\bibitem {bates}L. Bates, J. Sniatycki, \textit{Nonholonomic reduction},
Reports on Math. Phys. \textbf{32} (1993), 99-115.

\bibitem {boot}W.M. Boothby, \textit{An Introduction to Differentiable
Manifolds and Riemannian Geometry}, New York, Academic Press (1985).

\bibitem {can}F. Cantrijn, M. de Le\'{o}n, D. Mart\'{\i}n de Diego, \textit{On
almost-Poisson structures in nonholonomic mechanics}, Nonlinearity \textbf{12}
(1999), 721-737.

\bibitem {pepin-holo}J. Cari\~{n}ena, X. Gr\`{a}cia, G. Marmo, E.
Mart\'{\i}nez, M. Mu\~{n}oz-Lecanda, N. Roman-Roy, \textit{Geometric
Hamilton-Jacobi theory}, International Journal of Geometric Methods in Modern
Physics \textbf{3}, 7 (2006), 1417-1458.

\bibitem {pepin-noholo}J. Cari\~{n}ena, X. Gr\`{a}cia, G. Marmo, E.
Mart\'{\i}nez, M. Mu\~{n}oz-Lecanda, N. Roman-Roy, \textit{Geometric
Hamilton-Jacobi theory for nonholonomic dynamical systems}, International
Journal of Geometric Methods in Modern Physics \textbf{7}, 3 (2010), 431-454.

\bibitem {cg}H. Cendra, S. Grillo, \textit{Generalized nonholonomic mechanics,
servomechanisms and related brackets}, J. Math. Phys. \textbf{47} (2006), 2209-38.

\bibitem {crampin}M. Crampin, \textit{On the differential geometry of the
Euler--Lagrange equations and the inverse problem of Lagrangian dynamics}, J.
Phys. A \textbf{14} (1981), 2567--2575.

\bibitem {D}J.J. Duistermaat: \textit{On global action-angle coordinates,}
Comm. Pure Appl. Math. \textbf{33} (1980), 687--706.

\bibitem {fasso}F. Fass\`{o}, \textit{Superintegrable Hamiltonian Systems:
Geometry and Perturbations}, Acta Applicandae Mathematicae \textbf{87} (2005), 93--121.

\bibitem {rui}R. Fern\'{a}ndez, C. Laurent-Gengoux, P. Vanhaecke,
\textit{Global action-angle variables for non-commutative integrable systems},
arXiv:1503.00084 (2015).

\bibitem {gold}H. Goldstein,\textit{\ Classical Mechanics}, Addison-Wesley (1951).

\bibitem {gu}J. Grabowski, P. Urbanski, \textit{Algebroid-general differential
calculus on vector bundles}, J. Geom. Phys. \textbf{31} (1999), 111-141.

\bibitem {klein}J. Klein, \textit{Op\'{e}rateurs diff\'{e}rentielles sur les
variet\'{e}s presque tangentes}, C.R. Math. Acad. Sci. Paris \textbf{257}
(1963), 2392--2394.

\bibitem {kn}S. Kobayashi, K. Nomizu, \textit{Foundations of Differential
Geometry}, New York, John Wiley \& Son (1963).

\bibitem {koz}V. Kozlov, \textit{An extension of the Hamilton-Jacobi method},
J. Appl. Maths Mechs \textbf{60}, 6 (1996), 911-920.

\bibitem {eva}C. Laurent-Gengoux, E. Miranda, P. Vanhaecke, Action-angle
coordinates for integrable systems on Poisson manifolds, IMRN \textbf{13}, 8 (2011)1839-1869.

\bibitem {lmm}M. de Le\'on, J.C. Marrero, D. Mart{\'\i}n de Diego:
\textit{Linear almost Poisson structures and Hamilton-Jacobi equation.
Applications to nonholonomic Mechanics}, J. Geom. Mech. \textbf{2} (2010) 159--198.

\bibitem {hjp}M. de Le\'{o}n, D. Mart\'{\i}n de Diego, M. Vaquero: \textit{A
Hamilton-Jacobi theory on Poisson manifolds}, Journal of Geometric Mechanics
\textbf{6}, 1 (2014) 121-140.

\bibitem {mr}J.E. Marsden, T.S. Ratiu,\textit{ Introduction to Mechanics and
Symmetry, }New York, Springer-Verlag (1994).

\bibitem {mf}A.S. Mischenko, A.T. Fomenko: \textit{Generalized Liouville
Methods of integration of Hamiltonian systems, } Funct. Anal. Appl., \textbf{
12}, 2 (1978) 113--1978.

\bibitem {mrgm}J.E. Marsden, T.S. Ratiu,\textit{ Manifolds, Tensor Analysis
and Applications}, New York, Springer-Verlag (2001).

\bibitem {ortega}J.-P. Ortega, V. Planas-Bielsa, \textit{Dynamics on Leibniz
manifolds}, J. Geom. Phys. \textbf{52} (2004), 1-27.

\bibitem {van}A.J. van der Schaft, B.M. Maschke, \textit{On the Hamiltonian
formulation of nonholonomic mechanical systems}, Rep. Math. Phys. \textbf{34}
(1994), 225-233.

\bibitem {zung}N. T. Zung, \textit{Reduction and Integrability},
arXiv:math/0201087v2 [math.DS] (2002).
\end{thebibliography}
\end{document}